  \documentclass{cambridge6A}
  \usepackage[authoryear,round]{natbib}
  
  \usepackage[all]{xy}

\usepackage{array}
\usepackage{units}
\usepackage{textcomp}
\usepackage{multirow}

  \usepackage{rotating}
  \usepackage{floatpag}
  \rotfloatpagestyle{empty}

  \usepackage{amsthm}
\usepackage{amsmath}
\usepackage{amssymb}
  \usepackage{graphicx}

\newenvironment{lyxlist}[1]
{\begin{list}{}
{\settowidth{\labelwidth}{#1}
 \setlength{\leftmargin}{\labelwidth}
 \addtolength{\leftmargin}{\labelsep}
 }}
{\end{list}}

  \usepackage{multind}
  \makeindex{authors}
  \makeindex{subject}

  \theoremstyle{plain}% default
  \newtheorem{theorem}{Theorem}[chapter]
  \newtheorem{thm}{Theorem}[chapter] % lyx uses thm, rem

  \theoremstyle{definition}
  \newtheorem{rem}{Remark}[chapter]
  
  \newtheorem{example}[theorem]{Example}

  \theoremstyle{remark}

\renewcommand{\theenumii}{(\alph{enumii})}
\begin{document}

  \title[]
    {to appear in\\The New Handbook of Mathematical Psychology,\\to be published by Cambridge University Press.\\This version may differ from the published chapter.}

  \frontmatter
  \maketitle
  \tableofcontents
 
  \mainmatter

      % chapD.tex
% 2010/09/09, v2.10

% for multi-contributor books,  use \author
% for single-contributor books, though not required, use \chapterauthor

% uncomment \begin{abstract}...\end{abstract} for the Abstract to apppear

  \alphafootnotes

  \chapterauthor{Ehtibar Dzhafarov\footnotemark\ and Janne Kujala\footnotemark}
  
  \chapter{Probability, Random Variables, and Selectivity}

  \footnotetext[1]{Purdue University, USA}
    \footnotetext[2]{University of Jyv\"askyl\"a, Finland}

  \arabicfootnotes
    
     \contributor{Ehtibar Dzhafarov
    \affiliation{Purdue University}}

\section{What is it about?}

This chapter is about systems with several random outputs whose joint
distribution depends on several inputs. More specifically, it is about
selectiveness in the dependence of random outputs on inputs. That
is, we are concerned with the question of which of the several outputs
are influenced by which of the several inputs. A system can be anything:
a person, animal, group of people, neural network, technical gadget,
two entangled electrons running away from each other. Outputs are
responses of the system or outcomes of measurements performed on it.
Inputs are entities upon whose values the outputs of the system are
conditioned. Even if inputs are random variable in their own right,
the outputs are being conditioned upon every specific stimulus. Inputs
therefore are always deterministic (not random) entities insofar as
their relationship to random outputs is concerned.
\begin{example}
In a double-detection experiment, the stimulus presented in each trial
may consist of two flashes, say, right one and left one, separated
by some distance in visual field. Suppose that each flash can have
one of two contrast levels, one zero and one (slightly) above zero.
These contrasts play the role of two binary inputs, that we can call
$\lambda^{\text{left}}$ and $\lambda^{\text{right}}$ (each one with
values present/absent). The inputs are used in a completely crossed
experimental design: that is, the stimulus in each trial is described
by one of four combinations of the two inputs: $\left(\lambda^{\text{left}}=\textnormal{present, }\lambda^{\text{right}}=\textnormal{present}\right)$,
$\left(\lambda^{\text{left}}=\textnormal{present, }\lambda^{\text{right}}=\textnormal{absent}\right)$,
etc. In response to each such a combination (called a treatment),
the participant is asked to say whether the left flash was present
(yes/no) and whether the right flash was present (yes/no). These are
the two binary outputs, we can denote them $A^{\text{left}}$ and
$A^{\text{right}}$ (each with two possible values, yes/no). The outputs
are random variables. Theoretically, they are characterized by joint
distributions tied to each of four treatments:

\[%
\begin{tabular}{r|cc}
$\left(\lambda^{\text{left}}=i,\lambda^{\text{right}}=j\right)$ & $A^{\text{right}}=\text{yes}$ & $A^{\text{right}}=\text{no}$\tabularnewline
\hline 
$A^{\text{left}}=\text{yes}$ & $p_{\text{yes},\text{yes}}$ & $p_{\text{yes},\text{no}}$\tabularnewline
$A^{\text{left}}=\text{no}$ & $p_{\text{no},\text{yes}}$ & $p_{\text{no},\text{no}}$\tabularnewline
\end{tabular} \]where $i,j$ stand for ``present'' or ``absent'' each. Suppose
now that the experimenter hypothesizes that the response to the left
stimulus depends only on the contrast of the left stimulus, and the
response to the right stimulus depends only on the contrast of the
right stimulus, 
\[
\xymatrix{\lambda^{\text{left}}\ar[d] & \lambda^{\text{right}}\ar[d]\\
A^{\text{left}} & A^{\text{right}}
}
\]
This hypothesis can be justified, for example, by one's knowledge
that the separation between the locations of the flashes is too large
to allow for interference, and that subjectively, nothing seems to
change in the appearance of the left stimulus as the right one is
switched on and off, and vice versa. The meaning of this hypothesis
is easy to understand if the two random outputs are known to be stochastically
independent, which in this case means that, for every one of the four
treatments,
\[
p_{\text{yes},\text{yes}}=\Pr\left(A^{\text{left}}=\text{yes},A^{\text{right}}=\text{yes}\right)=\Pr\left(A^{\text{left}}=\text{yes}\right)\Pr\left(A^{\text{right}}=\text{yes}\right).
\]
In this case the test of the selectiveness consists in finding out
if the distribution of $A^{\text{left}}$, in this case defined by
$\Pr\left(A^{\text{left}}=\text{yes}\right)$, remains unchanged as
one changes the value of $\lambda^{\text{right}}$ while keeping $\lambda^{\text{left}}$
fixed, and analogously for $A^{\text{right}}$. The experimenter,
however, is likely to find out that stochastic independence in such
an experiment does not hold: for some, if not all of the four treatments,
\[
p_{\text{yes},\text{yes}}\not=\Pr\left(A^{\text{left}}=\text{yes}\right)\Pr\left(A^{\text{right}}=\text{yes}\right).
\]
Now the conceptual clarity may be lost. Does the lack of stochastic
independence invalidate the hypothesis that the outputs are selectively
influenced by the corresponding inputs? Indeed, one might reason that
it does, because if $A^{\text{left}}$ and $A^{\text{right}}$ are
not independent, then $A^{\text{left}}$ certainly ``depends on''
$A^{\text{right}}$, whence $A^{\text{left}}$ should also depend
on anything $A^{\text{right}}$ depends on (and this includes $\lambda^{\text{right}}$).
But one might also reason that stochastic relationship between the
two outputs can be ignored altogether. Cannot one declare that the
hypothesis in question holds if one establishes that the marginal
distributions (i.e., $\Pr\left(A^{\text{left}}=\text{yes}\right)$
and $\Pr\left(A^{\text{right}}=\text{yes}\right)$, taken separately)
are invariant with respect to changes in the non-corresponding inputs
(here, $\lambda^{\text{right}}$ and $\lambda^{\text{left}}$, respectively)?
We will see in this chapter that stochastic relationship must not
be ignored, but that lack of stochastic independence does not by itself
rule out selectiveness in the dependence of random outputs on inputs.\hfill$\square$
\end{example}
It is easy to generate formally equivalent examples by trivial modifications.
For instance, one can replace the two responses of a participant with
activity levels of two neurons, determining whether each of them is
above or below its background level. The two locations can be replaced
with two stimulus features (say, orientation and spatial frequency
of a grating pattern) that are hypothesized to selectively trigger
the responses from the two neurons. 

One can also easily modify any of such examples by increasing the
number of inputs and outputs involved, or increasing the number of
possible values per input or output. Thus, in the example with double-detection,
one can think of several levels of contrast for each of the flashes.
Or one can think of responses being multi-level confidence rating
instead of the binary yes/no. 

Let us consider a few more examples, however, to appreciate the variety
in the nature of inputs and outputs falling within the score of our
analysis.
\begin{example}
Let a very large group of students have to take three exams, in physics,
geometry, and French. Each student prepares for each of the exams,
and the preparation times are classified as ``short'' or ``long''
by some criteria (which may be different for different exams). The
three preparation times serve as the inputs in this example. We denote
them by $\lambda^{\text{physics}}$, $\lambda^{\text{geometry}}$,
and $\lambda^{\text{French}}$ (each with possible values short/long).
The outputs are scores the students eventually receive: $A^{\text{physics}}$,
$A^{\text{geometry}}$, and $A^{\text{French}}$ (say, from 0 to 100\%
each). The hypothesis to be tested is that preparation time for a
given subject selectively affects the score in that subject,
\[
\xymatrix{\lambda^{\text{physics}}\ar[d] & \lambda^{\text{geometry}}\ar[d] & \lambda^{\text{French}}\ar[d]\\
A^{\text{physics}} & A^{\text{geometry}} & A^{\text{French}}
}
\]
To see if this is the case we subdivide the group of students into
eight subgroups, corresponding to the eight combinations of the three
preparation times, 
\[
\left(\lambda^{\text{physics}}=\textnormal{short/long, }\lambda^{\text{geometry}}=\textnormal{short/long },\lambda=\textnormal{short/long }\right).
\]
Assuming each group is very large, we look at the joint distribution
of scores within each of them. The conceptual difficulty here stems
from the fact that, for any given treatment, test scores are typically
positively correlated rather than stochastically independent.\hfill$\square$
\end{example}

\begin{example}
Let us modify the previous example by assigning to each student in
each subject a binary grade, ``high'' or ``low,'' according as
the student is, respectively, above or below the median score in this
subject received by all student in the same preparation group. Thus,
in the preparation group $\left(\lambda^{\text{physics}}=\textnormal{long },\lambda^{\text{geometry}}=\textnormal{short },\lambda^{\text{French}}=\textnormal{short }\right)$,
if the median scores in physics is $m$, a student gets the grade
``high'' if her score is above $m$ and ``low'' if it is not.
This defines three outputs that we can call $B^{\text{physics}},B^{\text{geometry}},B^{\text{French}}$.
The hypothesis represented by the diagram
\[
\xymatrix{\lambda^{\text{physics}}\ar[d] & \lambda^{\text{geometry}}\ar[d] & \lambda^{\text{French}}\ar[d]\\
B^{\text{physics}} & B^{\text{geometry}} & B^{\text{French}}
}
\]
is more subtle than in the previous example. It says that if one factors
out the possible dependence of the median score in physics on all
three preparation times (with no selectiveness assumed in this dependence),
then whether a student's physics score will or will not fall above
the median may only depend on the preparation time for physics, and
not on the preparation times for two other subjects. And analogously
for geometry and French. Since the grades assigned to students are
binary, their theoretical distribution for each of the eight treatments
is given by eight joint probabilities
\[
\Pr\left(B^{\text{physics}}=\textnormal{high/low, }B^{\text{geometry}}=\textnormal{high/low, }B^{\text{French}}=\textnormal{high/low}\right).
\]
Again, the conceptual difficulty is in that this probability is not
typically equal to $\nicefrac{1}{8}$ for all combinations of the
high/low values, as it would have to be if the three random variables
were independent. Indeed, the marginal (separately taken) probabilities
here are, by the definition of median,
\[
\Pr\left(B^{\text{physics}}=\textnormal{high}\right)=\Pr\left(B^{\text{geometry}}=\textnormal{high}\right)=\Pr\left(B^{\text{French}}=\textnormal{high}\right)=\frac{1}{2}.
\]
This example also shows why it is not wise to ignore the joint distributions
and look at the marginal ones only. If we did this, none of the random
outputs $B^{\text{physics}},B^{\text{geometry}},B^{\text{French}}$
would be viewed as influenced by any of the inputs $\lambda^{\text{physics}},\lambda^{\text{geometry}},\lambda^{\text{French}}$.
But this view would clash with the fact that in different preparation
groups the corresponding joint probabilities will typically be different.\hfill$\square$ 
\end{example}

\begin{example}
This example is not from behavioral sciences but from quantum physics.
It is not as strange as it may appear to the reader. The fact is,
the mathematical formalisms independently developed to study selective
influences in psychology turn out to be identical to those developed
in quantum physics to study the types of determinism involved in the
behavior of so-called entangled particles. Two entangled particles
can be thought of as being created as a single particles and then
split into two mirror-images running away from each other. Particles
possess a property called spin, something that can be measured along
differently oriented spatial axes. In the case of so-called spin$\textnormal{-}\nicefrac{1}{2}$
particles, such as electrons, once an axis is chosen the spin can
attain one of only two possible values, referred to as ``spin-up''
and ``spin-down.'' Suppose that two entangled electrons run away
from each other towards two observers, Alice and Bob (a traditional
way of referring to them in quantum physics), with previously synchronized
clocks. At one and the same moment by these clocks Alice and Bob measure
spins of their respective electrons along axes they previously chose.
The nature of the entanglement is such that if the axes chosen by
the two observers are precisely the same, then the spin values recorded
will necessarily have opposite values: if Bob records spin-down, Alice
will record spin-up. Suppose that Bob always chooses one of two axes,
which we will denote $\lambda^{\text{Bob}}=\beta_{1}$ and $\lambda^{\text{Bob}}=\beta_{2}$.
We view $\lambda^{\text{Bob}}$, therefore, as one of the two inputs
of the system. The other input is the axis chosen by Alice, $\lambda^{\text{Alice}}$.
Let it also have two possible values, $\lambda_{1}^{\text{Alice}}=\alpha_{1}$
and $\lambda_{2}^{\text{Alice}}=\alpha_{2}$. The outcome of Bob's
measurement is the first of two outputs of the system. We denote it
by $A^{\text{Bob}}$, with the possible values ``spin-up'' and ``spin-down''.
The random output $A^{\text{Alice}}$, with the same two values, is
defined analogously. The theoretical representation of this situation
is given by the joint probabilities\[%
\begin{tabular}{r|cc}
$\left(\lambda^{\text{Alice}}=\alpha_{1},\lambda^{\text{Bob}}=\beta_{j},\right)$ & $A^{\text{Bob}}={\uparrow}$ & $A^{\text{Bob}}={\downarrow}$\tabularnewline
\hline 
$A^{\text{Alice}}={\uparrow}$ & $p_{\uparrow\uparrow}$ & $p_{\uparrow\downarrow}$\tabularnewline
$A^{\text{Alice}}={\downarrow}$ & $p_{\downarrow\uparrow}$ & $p_{\downarrow\downarrow}$\tabularnewline
\end{tabular} \]where $i$ and $j$ stand for 1 or 2 each. It is reasonable to
hypothesize that 
\[
\xymatrix{\lambda^{\text{Alice}}\ar[d] & \lambda^{\text{Bob}}\ar[d]\\
A^{\text{Alice}} & A^{\text{Bob}}
}
\]
In other words, the spin recorded by Alice may depend on which axes
she chose, but not on the axis chosen by Bob. And vice versa. But
the two outcomes here, for any of the four possible combinations of
Alice's and Bob's axes, are not stochastically independent. This makes
this situation formally identical to that described in the example
with double detection, except that in the entanglement paradigm the
invariance of the marginal distributions is guaranteed: $\Pr\left(A^{\text{Bob}}={\uparrow}\right)$
is the same no matter what axis was chosen by Alice, and vice versa.
In fact, it may very well be the case that these probabilities always
remain equal to $\nicefrac{1}{2}$, as in the second example with
the three exams.\hfill$\square$ 
\end{example}
Behavioral sciences abound with cases when selective influences are
assumed with respect to random variables whose realizations are not
directly observable. Rather these random variables are hypothetical
entities from which random variables with observable realizations
can be derived theoretically. Thus, one may posit the existence of
certain unobservable processes selectively influenced by certain experimental
manipulations and manifested by their contribution to observable response
times. For instance, one may assume the existence of processes called
perception and response choice with respective durations $A^{\text{percept}}$
and $A^{\text{response}}$, and assume that the observed response
time is $A^{\text{percept}}+A^{\text{response}}$. One can further
assume that stimulus characteristics selectively influence $A^{\text{percept}}$
and instruction versions (such as speed emphasis versus accuracy emphasis)
selectively influence $A^{\text{response}}.$ The conceptual problem
mentioned in the previous examples arises here if the two durations
are not assumed to be stochastically independent. 

In analyzing ``same-different'' judgments for pairs of sounds, the
observable entities are sounds $\lambda^{\text{first}}$ and $\lambda^{\text{second}}$,
each varying on several levels, and responses ``same'' or ``different''
for each pair of these sounds' levels. It is typically postulated,
however, that the response is a function (in the mathematical sense
of the word) of two unobservable random variables, $A^{\text{first}}$
and $A^{\text{second}}$, interpreted as internal representations
of the two sounds, their images. For instance, a model may postulate
that the response ``same'' is given if and only if the distance
between $A^{\text{first}}$ and $A^{\text{second}}$ in some metric
is less than some epsilon. It is reasonable to hypothesize then that
\[
\xymatrix{\lambda^{\text{first}}\ar[d] & \lambda^{\text{second}}\ar[d]\\
A^{\text{first}} & A^{\text{second}}
}
\]
Otherwise, why would one interpret $A^{\text{first}}$ and $A^{\text{second}}$
as ``separate'' respective images of $\lambda^{\text{first}}$ and
$\lambda^{\text{second}}$, rather than speaking of $A=\left(A^{\text{first}},A^{\text{second}}\right)$
as one image of the compound stimulus $\left(\lambda^{\text{first}},\lambda^{\text{second}}\right)$?

Stochastic independence of random outputs is, of course, a special
case of stochastic relationship. It is clear from our opening examples
that this is one case when the issue of defining and testing for selective
influences is conceptually transparent. Deterministic outputs are
a special case of random outputs, moreover, they can be formally considered
stochastically independent. To see that a deterministic output $a$
is influenced by an input $\lambda$ but not input $\lambda'$, see
if its value changes in response to changes in $\lambda$ but remains
constant if $\lambda'$ changes with $\lambda$ fixed. The only reason
for mentioning here this obvious consideration is this: there is a
wide class of theoretical models which deal with deterministic inputs
and and random outputs, but in which selectiveness of influences is
formulated as a relationship between deterministic entities, namely,
between the inputs and some parameters of the distributions of the
random outputs. Parameters of distributions are, by definition, deterministic
quantities. Such models require no special theory of selective influences.
\begin{example}
In multinomial processing tree models we see simple examples of random
variables related to inputs through parameters describing these variables'
distributions. A prototypical example is provided by R. Duncan Luce's
(1959) two-state low threshold model of detection,

\emph{
\[
\xymatrix{\lambda^{\textnormal{stimulus}}\ar[rr] &  & \bullet\ar@{=>}[dl]_{p}\ar@{=>}[dr]\\
\lambda^{\textnormal{payoff}}\ar[r]\ar@/_{1pc}/[rrr] & \bullet(\textnormal{detected})\ar@{=>}[drr]\ar@{=>}[d]_{q} &  & \bullet\left(\textnormal{not detected}\right)\ar@{=>}[d]^{r}\ar@{=>}[dll]\\
 & \bullet\left(\textnormal{No}\right) &  & \bullet\left(\textnormal{Yes}\right)
}
\]
}The processing flow is shown by the double-line arrows: from the
root of the tree to the root's children nodes, labeled ``detected''
and ``not detected,'' and from each of those to their children nodes,
labeled ``Yes'' and ``No.'' The labels $p$, $q$, and $r$ are
probabilities. The information shown in the processing tree is sufficient
for computations, except for one additional constraint: the model
stipulates that $qr=0$ (i.e., when one of the $q$ and $r$ is nonzero
the other one must be zero). The inputs $\lambda^{\text{stimulus}}$
and $\lambda^{\text{payoff}}$ are shown on the margins. A single-line
arrow pointing at a node of the tree indicates influence on the random
variable whose possible values are the children of this node. Stimulus
influences the distribution of the (unobservable) binary random variable
called ``detection state.'' It has two values occurring with probabilities
$p$ and $1-p$. Payoff is any procedure involving feedback and designed
to bias to various degrees the participants towards or against saying
``Yes.'' This input influences the (observable) random variable
``response.'' The point to note here is this: there is no reason
to consider the joint distributions of detection state and response
for different combinations of stimuli and payoffs; all we need is
to declare which of the three parameters of the model, $p,q,r$ depends
on which input, 

\[
p=p\left(\lambda^{\text{stimulus}}\right),q=q\left(\lambda^{\text{payoff}}\right),r=r\left(\lambda^{\text{payoff}}\right).
\]
This is simple and clear, even though the outputs ``detection state''
and ``response'' are not stochastically independent. \hfill$\square$ 
\end{example}
As it turns out, it is impossible to answer the questions posed in
this introductory section without getting ``back to basics,'' to
the foundational concepts of probability, random variable, joint distribution,
and dependence of joint distributions on deterministic variables.
It is simply impossible not to make mistakes and not to get hopelessly
confused in dealing with the issues of selective influences if one
is only guided by intuitive and informal understanding of these notions.
This applies even if the random variables involved are as simple as
binary responses. The first part of this chapter (Sections \ref{sec:What-is-a}-\ref{sec:Random-outputs-depending})
is dedicated to these foundational issues. The reader should be especially
attentive when we discuss the fact that not all random variables are
jointly distributed, that a set of random variables can always be
assigned a joint distribution in the absence of any constraints, but
that this may not be possible if the joint distribution should agree
with the known distributions of some subsets of this set of random
variables. Essentially, the issue of selective influences boils down
to establishing whether this is or is not possible in specific cases.
We deal with this issue beginning with Section \ref{sec:Selectiveness-in-the},
as well as the issue of methods by which one can determine whether
a particular pattern of selective influences holds. In Section \ref{Conditional determinism} we
show how the theory of selective influences applies to a classical problem of cognitive psychology,
the problem of determining, based on the overall response time, whether certain hypothetical 
processes involved in the formation of the response are concurrent or serial. The chapter concludes
with a brief guide to the relevant literature.

\section{\label{sec:What-is-a}What is a random variable?}

Let us begin with the notion of a\emph{ distribution} of a random
variable. The formal definition of this notion is as follows: the
distribution of a random variable $A$ is a triple 
\[
\overline{A}=\left(S,\Sigma,p\right),
\]
where
\begin{enumerate}
\item $S$ is some nonempty set, called the \emph{set of possible values}
of $A$;
\item $\Sigma$ is a \emph{sigma-algebra} \emph{over} $S$, which means
a collection of subsets of $S$, each called an\emph{ event} or a\emph{
measurable set}, such that

\begin{enumerate}[(a)]
\item $S\in\Sigma$,
\item if $S'\in\Sigma$, then $S-S'\in\Sigma$,
\item if $S_{1},S_{2},\ldots\in\Sigma$ (a finite or countably infinite
sequence), then 
\[
\bigcup_{i=1,2,\ldots}S_{i}\in\Sigma;
\]

\end{enumerate}
\item $p$ is some function (called \emph{probability measure}) from $\Sigma$
to $\left[0,1\right]$, such that $p\left(S'\right)$ for $S'\in\Sigma$
is interpreted as the probability with which a value of $A$ falls
in (belongs to) event $S'$; it is assumed that 

\begin{enumerate}[(a)]
\item $p\left(S\right)=1$,
\item (\emph{sigma-additivity}) if $S_{1},S_{2},\ldots\in\Sigma$ (a finite
or countably infinite sequence), and if in this sequence $S_{i}\cap S_{j}=\emptyset$
whenever $i\not=j$ (i.e., the subsets in the sequence are pairwise
disjoint), then 
\[
p\left(\bigcup_{i=1,2,\ldots}S_{i}\right)=\sum_{i=1,2,\ldots}p\left(S_{i}\right).
\]

\end{enumerate}
\end{enumerate}
The following consequences of this definition are easily derived:
\begin{enumerate}
\item $\emptyset\in\Sigma$ and $p\left(\emptyset\right)=0$;
\item if $S_{1},S_{2},\ldots\in\Sigma$, then $\bigcap_{i=1}^{\infty}S_{i}\in\Sigma$;
\item if $S_{1},S_{2},\ldots\in\Sigma$ and $S_{1}\subset S_{2}\subset\ldots$,
then 
\[
\lim_{i\rightarrow\infty}p\left(S_{i}\right)=p\left(\bigcup_{i=1}^{\infty}S_{i}\right);
\]

\item if $S_{1},S_{2},\ldots\in\Sigma$ and $S_{1}\supset S_{2}\supset\ldots$,
then 
\[
\lim_{i\rightarrow\infty}p\left(S_{i}\right)=p\left(\bigcap_{i=1}^{\infty}S_{i}\right);
\]

\item if $S_{1},S_{2}\in\Sigma$ and $S_{1}\subset S_{2}$, then $S_{2}-S_{1}\in\Sigma$
and 
\[
p\left(S_{1}\right)+p\left(S_{2}-S_{1}\right)=p\left(S_{2}\right);
\]

\item if $S_{1},S_{2}\in\Sigma$, then 
\[
p\left(S_{1}\cap S_{2}\right)\leq\min\left(p\left(S_{1}\right),p\left(S_{2}\right)\right)\leq\max\left(p\left(S_{1}\right),p\left(S_{2}\right)\right)\leq p\left(S_{1}\cup S_{2}\right).
\]

\end{enumerate}
Most of these consequences are known as \emph{elementary properties
of probability}. It is customary to write $p\left(S'\right)$ for
$S'\in\Sigma$ as $\Pr\left(A\in S'\right)$, if the distribution
of $A$ is known from the context.

We see that in order to know the distribution of a random variable
$A$ we have to know its set of possible values $S$ and a set of
specially chosen subsets of $S$, called events. And we should have
a procedure ``measuring'' each event, that is, assigning to it a
probability with which a value of $A$ (an element of $S$) falls
within this event (which is also described by saying that the event
in question ``occurs''). 
\begin{example}
\label{Example: Finite-.}For a finite $S$, the sigma-algebra is
usually defined as the power set, i.e., the set of all subsets of
$S$. For example, the distribution of the outcome $A$ of a roll
of a fair die can be represented by the distribution 
\[
\overline{A}=\left(S=\text{\{1,2,3,4,5,6\},\ensuremath{\Sigma}=\ensuremath{\mathcal{P}}(S)},p\right),
\]
where $\mathcal{P}(S)$ denotes the power set of $S$ and $p(\{s_{1},\dots,s_{k}\})=\nicefrac{k}{6}$
for any set $\{s_{1},\dots,s_{k}\}\in\Sigma$ of $k$ elements in
$S$. Similarly, the sum of two dice can be represented by the distribution
$\overline{A}=\left(S=\text{\{2,\ensuremath{\dots},12\},\ensuremath{\Sigma}=\ensuremath{\mathcal{P}}(S)},p\right)$,
where
\[
p(\{s_{1},\dots,s_{k}\})=\sum_{i=1}^{k}p(\{s_{i}\})
\]
and $p(\{s\})=\frac{1}{36}(6-|7-s|)$ gives the probability of each singleton (one-element subset) $\{s\}$. \hfill$\square$ 
\end{example}

\begin{example}
\label{Example: Real-line-.}Let $S$ be an interval of real numbers,
finite or infinite, perhaps the entire set $\mathbb{R}$ of real numbers.
For continuous distributions defined on $S$, at the very least we
want to be able to measure the probability of all intervals $\left(a,b\right)\subset S$.
This requirement implies that our sigma-algebra $\Sigma$ of events
must contain all so-called \emph{Borel subsets} of $S$. The Borel
sets form the smallest sigma-algebra $\Sigma$ over $S$ that contains
all open (or, equivalently, all closed) intervals. One can construct
this sigma algebra by the following \emph{recursive procedure}: (1)
include in $\Sigma$ all intervals in $S$; (2) add to this set of
intervals all countable unions of these intervals and of their complements;
(3) add to the previously obtained sets all countable unions of these
sets of their complements; (4) and so on. Clearly, these steps are
recursive applications of the operations (b) and (c) in the definition
of a sigma-algebra. Every Borel set will be obtained at some step
of this procedure. 

The Borel sigma-algebra is sufficient for most purposes, but often
the sigma-algebra is further enlarged by adding to all Borel sets
all \emph{null sets}. The latter are sets that can be covered by a
countable sequence of intervals with arbitrarily small total length
(see Section \ref{sec:Random-variables-narrow}). The motivation for
this extension is that anything that can be covered by an arbitrarily
small length should have its measure equal to zero (and for this it
should be measurable). The smallest sigma-algebra containing intervals
and null sets is called the \emph{Lebesgue sigma-algebra}. 

A continuous distribution on the real line can be defined using a
density function $f(a).$ The distribution is given by $\overline{A}=\left(S,\Sigma,p\right)$,
where $\Sigma$ is the Lebesgue sigma-algebra, and the probability
measure of a set $S_{A}\in\Sigma$ is given by the integral of the
density function $f$ over the set $S_{A}$, 
\[
p(S_{A})=\int_{S_{A}}f(a)\mathrm{d}a.
\]
(To be well defined for all Lebesgue-measurable sets $S_{A}$, the
integral here should be understood in the Lebesgue sense, but we need
not go into this.)\hfill$\square$ 
\end{example}
We see that measurability of a subset of $S$ is not a property of
the subset itself, but of this subset taken in conjunction with a
sigma-algebra $\Sigma$. Examples of \emph{non-measurable subsets}
of $S$ therefore are easily constructed: choose $\Sigma$ which is
not the entire power set of $S$, and choose a subset of $S$ which
is not in $\Sigma$. For instance, if $\Sigma=\left\{ \emptyset,\left\{ 1\right\} ,\left\{ 2,3\right\} ,\left\{ 1,2,3\right\} \right\} $
over the set $S=\left\{ 1,2,3\right\} $, then the single-element
subset $\left\{ 3\right\} $ is non-measurable. This means that if
$A$ is distributed as $\left(S,\Sigma,p\right)$, the probability
$p\left(\left\{ 3\right\} \right)$ with which $A$ falls in $\left\{ 3\right\} $
(or, simply, equals 3) is undefined. This example may seem artificial,
as nothing prevents one from complementing $\Sigma$ with all other
subsets of $S=\left\{ 1,2,3\right\} $ (i.e., to assume that $p$
is defined for all of them even if it is only known for some). If
$S$ is an interval of reals, however, then there are deeper reasons
for not including in $\Sigma$ all subsets of $S$.

It is obvious that different random variables can have one and the
same distribution. For instance, Peter and Paul can flip a fair coin
each, and describe the outcomes by one and the same distribution
\[
\overline{A}=\left(S=\left\{ 0,1\right\} ,\Sigma=\left\{ \emptyset,\left\{ 0\right\} ,\left\{ 1\right\} ,\left\{ 0,1\right\} \right\} ,p\left(\Sigma\right)=\left\{ 0,\frac{1}{2},\frac{1}{2},1\right\} \right).
\]
To distinguish one random variable from another, therefore, it is
not sufficient to know its distribution. We should, in addition, have
a \emph{label} or \emph{name} for the random variable: for instance,
we can identify one random variable as $\text{coin}_{1}$, distributed
as $\overline{A}$, and another as $\text{coin}_{2}$, also distributed
as $\overline{A}$.

Generally speaking, a random variable $A$ can be viewed as a quadruple
$\left(\iota_{A},S,\Sigma,p\right)$, where $\iota_{A}$ is its unique
name and $\overline{A}=\left(S,\Sigma,p\right)$ is its distribution.
We do not need to be that formal, however, as the notation for a random
variable, $A$, also serves as its name. (The reader familiar with
the conventional definition of a random variable as a measurable function
on a sample space should wait patiently until Sections \ref{sec:Random-variables-as}
and \ref{sec:Stochastically-unrelated-and}. A function may serve
as an identifying label too.)
\begin{rem}
Alternatively, one can assume that the name of a random variable is
always (implicitly) part of the elements of its domain $S$. For instance,
the domain for one of the two coins mentioned above may be defined
as $S^{1}=\left\{ \left(0,\text{coin}_{1}\right),\left(1,\text{coin}_{1}\right)\right\} $
and for another as $S^{2}=\left\{ \left(0,\text{coin}_{2}\right),\left(1,\text{coin}_{2}\right)\right\} $.
The sigma-algebras $\Sigma^{1}$ and $\Sigma^{2}$ then have to be
(re)defined accordingly. If this approach is followed consistently,
every random variable is uniquely determined by its distribution.
We do not follow this route in this chapter.
\end{rem}

\section{\label{sec:Jointly-distributed-random}Jointly distributed random
variables}

Let $A$, $B$, and $C$ be random variables with distributions $\overline{A}=\left(S^{1},\Sigma^{1},p_{1}\right)$,
$\overline{B}=\left(S^{2},\Sigma^{2},p_{2}\right)$, and $\overline{C}=\left(S^{3},\Sigma^{3},p_{3}\right)$. 
\begin{rem}
We will consistently use numerical superscripts to refer to the domain
sets for random variables, to sigma-algebras over these sets, and
later to random variables and inputs. Notation $S^{3}$, for example,
always refers to a domain set of some random variable, \emph{not}
to the Cartesian product $S\times S\times S$. This should not cause
any difficulties, as we use numerical exponents in this chapter only
twice, and both times this is explicitly indicated. 
\end{rem}
Let $S_{A}\in\Sigma^{1}$, $S_{B}\in\Sigma^{2}$, and $S{}_{C}\in\Sigma^{3}$
be three events. We know that $p_{1}\left(S_{A}\right)$ is interpreted
as the probability with which a value of $A$ falls in $S_{A}$ (or,
the probability that the event $S_{A}$ ``occurs''); and analogously
for $p_{2}\left(S_{B}\right)$ and $p_{3}\left(S_{C}\right)$. We
also speak of events \emph{occurring jointly}, or \emph{co-occurring},
a concept whose substantive meaning we will discuss in Section \ref{sec:Stochastically-unrelated-and}.
For now we will take it formally. In order to speak of $S_{A},S_{B},S_{C}$
co-occurring and to ask of the probabilities with which they co-occur,
we have to introduce a new random variable, denoted $D_{ABC}$. As
any random variable, it is defined by some unique name (e.g., ``$D_{ABC}$'')
and a distribution 
\[
\overline{D_{ABC}}=\left(S^{123},\Sigma^{123},p_{123}\right).
\]
The set $S^{123}$ of possible values of $D_{ABC}$ is the Cartesian
product $S^{1}\times S^{2}\times S^{3}$ (the set of all ordered triples
with the first components chosen from $S^{1}$, the second from $S^{2}$,
the third from $S^{3}$). The sigma-algebra $\Sigma^{123}$ is denoted
$\Sigma^{1}\otimes\Sigma^{2}\otimes\Sigma^{3}$ and defined as the \emph{smallest
sigma-algebra} containing the Cartesian products $S_{A}\times S_{B}\times S_{C}$
for all $S_{A}\in\Sigma^{1}$, $S_{B}\in\Sigma^{2}$ and $S_{C}\in\Sigma^{3}$.
This means that $\Sigma^{123}=\Sigma^{1}\otimes\Sigma^{2}\otimes\Sigma^{3}$
is a set of subsets of $S^{1}\times S^{2}\times S^{3}$, such that
\begin{enumerate}
\item it contains all the Cartesian products $S_{A}\times S_{B}\times S_{C}$
just mentioned;
\item with every subset $S'$ it contains, it also contains the complement
$S^{123}-S'$;
\item with every sequence of subsets $S_{1},S_{2}\ldots$ it contains, it
also contains their union, $\bigcup_{i=1,w,\ldots}S_{i}$;
\item it is included in any other set of subsets of $S^{1}\times S^{2}\times S^{3}$
satisfying 1-2-3 above.
\end{enumerate}
The probability measure $p_{123}$ is called a \emph{joint probability
measure}. It should satisfy the general requirements of a probability
measure, namely:
\[
p_{123}\left(S^{1}\times S^{2}\times S^{3}\right)=1,
\]
and
\[
p_{123}\left(\bigcup_{i=1,2,\ldots}S_{i}\right)=\sum_{i=1,2,\ldots}p\left(S_{i}\right)
\]
for any sequence of pairwise disjoint elements $S_{1},S_{2},\ldots$
of $\Sigma^{123}$. In addition, $p_{123}$ should satisfy the following
\emph{1-marginal probability} \emph{equations}: for any $S_{A}\in\Sigma^{1}$,
$S_{B}\in\Sigma^{2}$ and $S_{C}\in\Sigma^{3}$,

\[
\begin{array}{c}
p_{123}\left(S_{A}\times S^{2}\times S^{3}\right)=p_{1}\left(S_{A}\right),\\
\\
p_{123}\left(S^{1}\times S_{B}\times S^{3}\right)=p_{2}\left(S_{B}\right),\\
\\
p_{123}\left(S^{1}\times S^{2}\times S_{C}\right)=p_{3}\left(S_{C}\right).
\end{array}
\]

\begin{example}
\label{Example:joint-distribution-of-three}Let 
\[
S=\{0,1\},\;\Sigma=\{\emptyset,\{0\},\{1\},\{0,1\}\},
\]
and let the random variables $A$, $B$, and $C$ be distributed as
\[
\overline{A}=\left(S,\Sigma,p_{1}\right),\;\overline{B}=\left(S,\Sigma,p_{2}\right),\;\overline{C}=\left(S,\Sigma,p_{3}\right),
\]
where 
\[
p_{1}(\Sigma)=\{0,\nicefrac{1}{2},\nicefrac{1}{2},1\},\; p_{2}(\Sigma)=\{0,\nicefrac{1}{4},\nicefrac{3}{4},1\},\; p_{3}(\Sigma)=\{0,1,0,1\}.
\]
A joint distribution of $A,B,C$ is defined on the product sigma-algebra
$\Sigma^{123}=\Sigma\otimes\Sigma\otimes\Sigma$, which is the smallest
sigma-algebra containing all Cartesian products $S_{A}\times S_{B}\times S_{C}$
such that $S_{A},S_{B},S_{C}\in\Sigma.$ As the Cartesian products
include those of all singletons (one-element subsets) $\{(a,b,c)\}=\{a\}\times\{b\}\times\{c\}$,
and all subsets of $S\times S\times S$ can be formed by finite unions
of these, the product sigma algebra $\Sigma\otimes\Sigma\otimes\Sigma$
is the full power set of $S\times S\times S$. One possible joint
distribution for $A,B,C$ is given by 
\[
\overline{D_{ABC}}=\left(S^{123}=S\times S\times S,\Sigma^{123}=\Sigma\otimes\Sigma\otimes\Sigma,p_{123}\right),
\]
where 
\[
p_{123}(S_{ABC})=\sum_{(a,b,c)\in S_{ABC}}p_{123}(\{(a,b,c)\})
\]
and $p_{123}(\{(a,b,c)\})$ is given by the table\[\small%
\begin{tabular}{ccc|ccccc|c}
$a$ & $b$ & $c$ & $p_{123}(\{(a,b,c)\})$ &  & $a$ & $b$ & $c$ & $p_{123}(\{(a,b,c)\})$\tabularnewline
\cline{1-4} \cline{6-9} 
0 & 0 & 0 & $\nicefrac{1}{16}$ &  & 1 & 0 & 0 & $\nicefrac{3}{16}$\tabularnewline
0 & 0 & 1 & 0 &  & 1 & 0 & 1 & 0\tabularnewline
0 & 1 & 0 & $\nicefrac{7}{16}$ &  & 1 & 1 & 0 & $\nicefrac{5}{16}$\tabularnewline
0 & 1 & 1 & 0 &  & 1 & 1 & 1 & 0\tabularnewline
\end{tabular}\]Let us verify that this distribution satisfies the 1-marginal probability
equations and is thus a proper joint distribution of $A,B,C:$
\begin{align*}
p_{123}(\{0\}\times S\times S) & =\nicefrac{1}{16}+0+\nicefrac{7}{16}+0=\nicefrac{1}{2}=p_{1}(\{0\}),\\
p_{123}(\{1\}\times S\times S) & =\nicefrac{3}{16}+0+\nicefrac{5}{16}+0=\nicefrac{1}{2}=p_{1}(\{1\}),\\
p_{123}(S\times\{0\}\times S) & =\nicefrac{1}{16}+0+\nicefrac{3}{16}+0=\nicefrac{1}{4}=p_{2}(\{0\}),\\
p_{123}(S\times\{1\}\times S) & =\nicefrac{7}{16}+0+\nicefrac{5}{16}+0=\nicefrac{3}{4}=p_{2}(\{1\}),\\
p_{123}(S\times S\times\{0\}) & =\nicefrac{1}{16}+\nicefrac{7}{16}+\nicefrac{3}{16}+\nicefrac{5}{16}=1=p_{3}(\{0\}),\\
p_{123}(S\times S\times\{1\}) & =0+0+0+0=0=p_{3}(\{1\}).
\end{align*}
For each 1-marginal, it suffices to verify the probabilities of the
points $0$ and $1$ as the probability values for singletons fully
determine the discrete distributions.\hfill$\square$ 
\end{example}
The random variable $D_{ABC}$ is commonly called a \emph{vector of
the (jointly distributed) random variables} $A$, $B$, and $C$,
and it is denoted $\left(A,B,C\right)$. We will use this vectorial
notation in the sequel. One should keep in mind, however, that any
such a vector is a random variable in its own right. Furthermore,
one should keep in mind that the distribution $\overline{\left(A,B,C\right)}$,
called the \emph{joint distribution} with respect to the individual
random variables $A,B,C$, is not uniquely determined by these $A,B,C$.
Specifically, although the set $S^{123}=S^{1}\times S^{2}\times S^{3}$
and the sigma-algebra $\Sigma^{123}=\Sigma^{1}\otimes\Sigma^{2}\otimes\Sigma^{3}$
are uniquely determined by the sets and sigma-algebras in the distributions
$\overline{A}$, $\overline{B},$ and $\overline{C}$, there can generally
be more than one\emph{ }joint probability measure\emph{ }$p_{123}$.
The\emph{ }individual $p_{1}$, $p_{2}$, and $p_{3}$ only serve
as constraints, in the form of the 1-marginal probability equations
above.

$A$, $B$, and $C$ in $\left(A,B,C\right)$ are called \emph{stochastically
independent} if, for any $S_{A}\in\Sigma^{1}$, $S_{B}\in\Sigma^{2}$
and $S_{C}\in\Sigma^{3}$,
\[
p_{123}\left(S_{A}\times S_{B}\times S_{C}\right)=p_{1}\left(S_{A}\right)p_{2}\left(S_{B}\right)p\left(S_{C}\right).
\]
This joint probability measure always satisfies the 1-marginal probability
equations. 
\begin{example}
Let $A$ and $B$ be standard normally distributed random variables.
A bivariate normal joint distribution $(A,B)(\rho)$ can be defined
with the density function
\[
f_{12}(a,b;\rho)=\frac{1}{2\pi\sqrt{1-\rho^{2}}}\exp\left(-\frac{a^{2}+b^{2}-2\rho ab}{2(1-\rho^{2})}\right),
\]
where $-1<\rho<1$ denotes the correlation coefficient. The sigma
algebra $\Sigma_{12}=\Sigma_{1}\otimes\Sigma_{2}$ of the joint distribution
is the product of two Lebesgue sigma-algebras (called a Lebesgue sigma-algebra
itself). The 1-marginal probability equations can be verified by checking
that integrating out either $a$ or $b$ yields the standard normal
density function with respect to the remaining variable. The probability
measure for $C=(A,B)(\rho)$ is obtained as
\[
p_{12}(S_{C})=\int_{(a,b)\in S_{c}}f_{12}(a,b;\rho)\mathrm{d}(a,b).
\]
Do $C=\left(A,B\right)\left(\rho_{1}\right)$ and $D=\left(A,B\right)\left(\rho_{2}\right)$
with $\rho_{1}\ne\rho_{2}$ exclude each other? Not in the sense that
defining one of them makes the other meaningless. They both can be
defined as variables of interest. But $C$ and $D$ cannot be jointly
distributed. \hfill$\square$ 
\end{example}
The reverse relationship between joint and marginal distributions
is more straightforward: the distribution $\overline{\left(A,B,C\right)}$
uniquely determines the distributions and identity of $A$, $B$,
$C$, called the \emph{1-marginal random variables} with respect to
$\left(A,B,C\right)$, as well as the distributions and identity of
$\left(A,B\right)$, $\left(B,C\right)$, and $\left(A,C\right)$,
called the \emph{2-marginal random variables} with respect to $\left(A,B,C\right)$.
Thus, in the distribution $\overline{A}$ the set $S^{1}$ is the
projection $\mathrm{Proj}_{1}$ of the set $S^{123}=S^{1}\times S^{2}\times S^{3}$,
defined by
\[
\mathrm{Proj}_{1}\left(a,b,c\right)=a.
\]
The sigma-algebra $\Sigma^{1}$ consists of the projections $\mathrm{Proj}_{1}$
of the elements of the sigma-algebra $\Sigma^{123}=\Sigma^{1}\otimes\Sigma^{2}\otimes\Sigma^{3}$
having the form $S_{A}\times S^{2}\times S^{3}$. And the probability
measure $p_{1}$ is determined by the 1-marginal probability equations.
The \emph{2-marginal distributions} $\overline{\left(A,B\right)}$,
$\overline{\left(B,C\right)}$, and $\overline{\left(A,C\right)}$
are found analogously. For example, if one defines function $\mathrm{Proj}_{23}$
by
\[
\mathrm{Proj}_{23}\left(a,b,c\right)=\left(b,c\right),
\]
we have 
\[
\overline{\left(B,C\right)}=\left(S^{23},\Sigma^{23},p_{23}\right),
\]
where
\[
S^{23}=\mathrm{Proj}_{23}\left(S^{1}\times S^{2}\times S^{3}\right),
\]
$\Sigma^{23}$ consists of the sets of the form
\[
\mathrm{Proj}_{23}\left(S^{1}\times S_{BC}\right),\: S_{BC}\in\Sigma^{2}\otimes\Sigma^{3},
\]
and
\[
p_{23}\left(S_{BC}\right)=p_{123}\left(S^{1}\times S_{BC}\right).
\]
The last equality is one of the three \emph{2-marginal probability
equations }(the remaining two being for $p_{12}$ and $p_{13}$).\emph{ }

One can check that 
\[
S^{23}=S^{2}\times S^{3},
\]
and
\[
\Sigma^{23}=\Sigma^{2}\otimes\Sigma^{3},
\]
which is the smallest sigma-algebra containing the Cartesian products
$S_{B}\times S_{C}$ for all $S_{B}\in\Sigma^{2}$ and $S_{C}\in\Sigma^{3}$.
In other words, the set $S^{23}$ and the sigma-algebra $\Sigma^{23}$
over it in the 2-marginal distribution are precisely the same as if
they were formed for a joint distribution $\overline{\left(B,C\right)}$
with respect to the 1-marginal distributions $\overline{B}$ and $\overline{C}$.
Moreover, the 2-marginal probability $p_{23}$ is a joint probability
satisfying the 1-marginal probability equations
\[
\begin{array}{c}
p_{23}\left(S_{B}\times S^{3}\right)=p_{2}\left(S_{B}\right),\\
\\
p_{23}\left(S^{2}\times S_{C}\right)=p_{3}\left(S_{C}\right).
\end{array}
\]

\begin{example}
Continuing from Example~\ref{Example:joint-distribution-of-three},
we can derive the following 2-marginals (and 1-marginals shown at
the sides of the 2-marginals):\[\small%
\begin{tabular}{c|cc|c}
\!\!\!$p_{12}(\{(a,b)\})$\!\!\! & \!$b=0$\!\! & \!\!$b=1$\! & \tabularnewline
\hline 
$a=0$ & $\nicefrac{1}{16}$ & $\nicefrac{7}{16}$ & $\nicefrac{1}{2}$\tabularnewline
$a=1$ & $\nicefrac{3}{16}$ & $\nicefrac{5}{16}$ & $\nicefrac{1}{2}$\tabularnewline
\hline 
 & $\nicefrac{1}{4}$ & $\nicefrac{3}{4}$ & \tabularnewline
\end{tabular}\quad{} %
\begin{tabular}{c|cc|c}
\!\!\!$p_{12}(\{(b,c)\})$\!\!\! & \!$c=0$\!\! & \!\!$c=1$\! & \tabularnewline
\hline 
$b=0$ & $\nicefrac{1}{4}$ & $0$ & $\nicefrac{1}{4}$\tabularnewline
$b=1$ & $\nicefrac{3}{4}$ & $0$ & $\nicefrac{3}{4}$\tabularnewline
\hline 
 & $1$ & $0$ & \tabularnewline
\end{tabular}\]\[\small %
\begin{tabular}{c|cc|c}
\!\!\!$p_{12}(\{(a,c)\})$\!\!\! & \!$c=0$\!\! & \!\!$c=1$\! & \tabularnewline
\hline 
$a=0$ & $\nicefrac{1}{2}$ & $0$ & $\nicefrac{1}{2}$\tabularnewline
$a=1$ & $\nicefrac{1}{2}$ & $0$ & $\nicefrac{1}{2}$\tabularnewline
\hline 
 & $1$ & $0$ & \tabularnewline
\end{tabular}\]

\hfill$\square$ 
\end{example}
{\small It} should be clear now how one should generalize the notion
of a joint distribution to an arbitrary number $n$ of random variables,
$A^{1},\ldots,A^{n}$, and how to define $k\textnormal{-}$marginal
distributions for $k=1,\ldots,n$ ($n\textnormal{-}$marginal distributions
being permutations of the joint one, including itself). 
\begin{rem}
\label{REM: For-an-infinite}For an infinite set of random variables
(countable or not) the definition of a joint distribution is less
obvious. We will not deal with this notion in this chapter except
for mentioning it occasionally, for completeness sake. With little
elaboration, let $\left(A^{k}:k\in K\right)$ be an indexed family
of random variables (with an arbitrary indexing set $K$), each distributed
as $\left(S^{k},\Sigma^{k},p_{k}\right)$. We say that the random
variables in $\left(A^{k}:k\in K\right)$ are jointly distributed
if $A=\left(A^{k}:k\in K\right)$ is a random variable with the distribution
\[
\overline{A}=\left(\prod_{k\in K}S^{k},\bigotimes_{k\in K}\Sigma^{k},p\right),
\]
where \end{rem}
\begin{enumerate}
\item $\prod_{k\in K}S^{k}$ is the Cartesian product of the sets $S^{k}$
(its elements are functions choosing for each element of $K$ an element
of $S^{k}$); 
\item $\bigotimes_{k\in K}\Sigma^{k}$ is the smallest sigma-algebra containing
sets of the form $S'\times\prod_{k\in K-\left\{ k_{0}\right\} }S^{k}$,
for all $k_{0}\in K$ and $S'\in\Sigma^{k_{0}}$;
\item $p$ is a probability measure on $\bigotimes_{k\in K}\Sigma^{k}$
such that $p\left(S'\times\prod_{k\in K-\left\{ k_{0}\right\} }S^{k}\right)=p_{k_{0}}\left(S'\right)$,
for all $k_{0}\in K$ and $S'\in\Sigma^{k_{0}}$.
\end{enumerate}
The random variables $A^{k}$ in $A=\left(A^{k}:k\in K\right)$ are
said to be stochastically independent if any finite subset of them
consists of stochastically independent elements.

\begin{rem}
Marginal random variables sometimes have to be defined hierarchically.
Consider, for example, $A'=\left(A,B\right)$ and $B'=\left(C,D\right)$.
Then $C'=\left(A',B'\right)$ has the 1-marginal distributions $\overline{A'}=\overline{\left(A,B\right)}$
and $\overline{B'}=\overline{\left(C,D\right)}$. And $A'=\left(A,B\right)$,
in turn, has 1-marginal distributions $\overline{A}$ and $\overline{B}$.
It may sometimes be convenient to speak of all of $\left(A,B\right)$,
$\left(C,D\right)$, $A$, $B$, $C$, $D$ as marginal random variables
with respect to a random variable $C'=\left(\left(A,B\right),\left(C,D\right)\right)$.
Note that $\left(\left(A,B\right),\left(C,D\right)\right)$, $\left(\left(A,B,C\right),D\right)$,
$\left(A,\left(B,\left(C,D\right)\right)\right)$, etc. are all distributed
as $\left(A,B,C,D\right)$, because the Cartesian product $S^{1}\times S^{2}\times S^{3}\times S^{4}$
and the product sigma algebra $\Sigma^{1}\otimes\Sigma^{2}\otimes\Sigma^{3}\otimes\Sigma^{3}$
are associative. The random variables $\left(\left(A,B\right),\left(C,D\right)\right)$,
$\left(\left(A,B,C\right),D\right)$, $\left(A,\left(B,\left(C,D\right)\right)\right)$,
etc. differ in their labeling only. (In the infinite case (Remark
\ref{REM: For-an-infinite}) the formal definition is rather straightforward,
but it involves potentially more than a finite number of hierarchical
steps. We will assume that the notion is clear and a formal definition
may be skipped\label{REM: Marginal-random-variables}.) 
\end{rem}

\section{\textmd{\normalsize \label{sec:Random-variables-narrow}}Random variables
in the narrow sense}

The concept of a random variable used in this chapter is very general,
with no restrictions imposed on the sets and sigma-algebras in their
distributions. Sometimes such random variables are referred to as
\emph{random entit}ies, \emph{random element}s, or \emph{random variables
in the broad sense}, to distinguish them from \emph{random variables
in the narrow sense}. The latter are most important in applications.
In particular, all our example involve random variables in the narrow
sense. They can be defined as follows. Let $A$ be distributed as
$\overline{A}=\left(S,\Sigma,p\right)$. 

(i) If $S$ is countable, $\Sigma$ is the power set of $S$ (the
set of all its subsets), then $A$ is a random variable in the narrow
sense; 

(ii) if $S$ is an interval of real numbers, $\Sigma$ is the \emph{Lebesgue
sigma-algebra} over $S$ (as defined in Example~\ref{Example: Real-line-.}),
then $A$ is a random variable in the narrow sense; 

(iii) if $A_{1},\ldots,A_{n}$ are random variables in the narrow
sense, then any jointly distributed vector $\left(A_{1},\ldots,A_{n}\right)$
is a random variable (also referred to as a \emph{random vector}) in the
narrow sense.

Random variables satisfying (i) are called \emph{discrete}. The distribution
of such a random variable is uniquely determined by the probabilities
assigned to its singleton (one-element) subsets. These probabilities
can also be viewed as assigned to the elements themselves, in
which case they form a \emph{probability mass function}. An example
of a discrete random variable is given in Example \ref{Example: Finite-.}.
But $S$ may also be countably infinite. 
\begin{example}
Let $S$ be the set of positive integers $\left\{ 1,2,\ldots,n,\ldots\right\} $,
and let $p\left(\left\{ n\right\} \right)=\alpha^{n-1}\left(1-\alpha\right)$,
where $\alpha$ is a constant in $\left[0,1\right]$. This defines
a discrete random variables interpreted as the number of independent
trials $n$ with binary outcomes (success/failure) until the first
failure. It is customary to replace (or even confuse) $p\left(\left\{ n\right\} \right)$
with the probability mass function function $p^{*}\left(n\right)=p\left(\left\{ n\right\} \right)$.\hfill$\square$ 
\end{example}
Random variables satisfying (ii) are called \emph{continuous} (see
Example \ref{Example: Real-line-.}). Any such a variable can be viewed
as having $S$ extended to the entire set of reals, and its distribution
is uniquely determined by the \emph{distribution function} 
\[
F\left(x\right)=p\left(\left(-\infty,x\right]\right),
\]
for every real $x$. The function $F\left(x\right)$ has the following
properties: 
\begin{enumerate}
\item it is nondecreasing; 
\item as $x\rightarrow-\infty$, $F\left(x\right)\rightarrow0$; 
\item as $x\rightarrow\infty$, $F\left(x\right)\rightarrow1$; 
\item for any real $x_{0}$, as $x\rightarrow x_{0}+$, $F\left(x\right)\rightarrow F\left(x_{0}\right)$
(right-continuity);
\item for any real $x_{0}$, as $x\rightarrow x_{0}+$, $F\left(x\right)$
tends to a limit.
\end{enumerate}
$F\left(x\right)$ generally is not left-continuous: as $x\rightarrow x_{0}-$,
the limit of $F\left(x\right)$ need not coincide with $F\left(x_{0}\right)$, the function may
instead ``jump'' from the value of $\lim_{x\rightarrow x_{0}-}F\left(x\right)$
to $F\left(x_{0}\right)$. The difference $F\left(x_{0}\right)-\lim_{x\rightarrow x_{0}-}F\left(x\right)$
equals $p\left(\left\{ x_{0}\right\} \right)$, so the jumps
occur if and only if $p\left(\left\{ x_{0}\right\} \right)>0$. A
distribution function cannot have more than a countable set of jump points.
For any two reals $x_{1}\leq x_{2}$, 
\[
F\left(x_{2}\right)-F\left(x_{1}\right)=p\left(\left(x_{1},x_{2}\right]\right).
\]

\begin{example}
\label{A-discrete-continuous}A discrete random variable can always
be redefined as a continuous one. Thus, the variable in the previous
example can be redefined into a random variable $X$ whose distribution
is given by 
\[
F\left(x\right)=\left\{ \begin{array}{ccc}
0 & \textnormal{for} & x<1\\
\alpha^{n-1}\left(1-\alpha\right) & \textnormal{for} & \left\lfloor x\right\rfloor =n\geq1,
\end{array}\right.
\]
where $\left\lfloor x\right\rfloor $ is the floor function (the largest
integer not exceeding $x$).\hfill$\square$
\end{example}
The Lebesgue sigma-algebra over the reals, as defined in Example~\ref{Example: Real-line-.},
is the smallest sigma-algebra including all intervals and all \emph{null
sets}. A subset $S'$ of reals is a null set if, for any $\varepsilon>0$,
however small, $S'$ is contained within a union of open intervals
$S_{1},S_{2},\ldots$ whose overall length is less than $\varepsilon$.
An empty set is, obviously a null set, and so is a single point, and
a countable set of points. 
\begin{rem}
Let us prove that a countable set of points is a null set, to better
understand the definition. Enumerate this set as $x_{1},x_{2},\ldots$,
choose an $\varepsilon>0$, and enclose each $x_{i}$ into interval
$\left]x-\frac{\varepsilon}{2^{i+1}},x+\frac{\varepsilon}{2^{i+1}}\right[$
. The length of this interval is $\frac{\varepsilon}{2^{i}}$, whence
the overall length of the system of such intervals cannot exceed 
\[
\sum_{i=1,2,\ldots}\frac{\varepsilon}{2^{i}}\leq\varepsilon.
\]
We conclude that a countable subset of $S$ is a null set. There are
uncountable null sets.
\end{rem}
As should be clear from our discussion of jumps and Example
\ref{A-discrete-continuous}, a null set may have a nonzero probability.
If this does not happen, i.e., if $F\left(x\right)$ has no jumps,
the distribution of the random variable is called \emph{absolutely
continuous}.

Finally, the combination rule (iii) allows one to form vectors of
discrete, continuous, and mixed jointly distributed random variables
using the construction discussed in Section \ref{sec:Jointly-distributed-random}.

\section{\label{sec:Functions-of-random}Functions of random variables}

Let $A$ be a random variable with distribution $\overline{A}=\left(S^{1},\Sigma^{1},p_{1}\right)$,
let $S^{2}$ be some set, and let $f:S^{1}\rightarrow S^{2}$ be some
function. Consider some sigma algebra $\Sigma^{2}$ of events over
$S^{2}$. For every $S_{B}\in\Sigma^{2}$ one can determine the subset
of all elements of $S^{1}$ that are mapped by $f$ into $S_{B}$,
\[
f^{-1}\left(S_{B}\right)=\left\{ a\in S^{1}:f\left(a\right)\in S_{B}\right\} .
\]
This subset, $f^{-1}\left(S_{B}\right)$, does not have to be an event
in $\Sigma^{1}$. If it is, for every $S_{B}\in\Sigma^{2}$, then
$f$ is said to be a \emph{measurable function} (or $\Sigma^{1}\rightarrow\Sigma^{2}\textnormal{-}$\emph{measurabl}e\emph{
function}, to be specific). Measurability of a function therefore
is not a property of the function itself, but of the function taken
in conjunction with two sigma-algebras. In particular, given $S^{1}$
and $\Sigma^{1}$, any \emph{onto} function $f:S^{1}\rightarrow S^{2}$
(one with $f\left(S^{1}\right)=S^{2}$) will be measurable if we agree
to define $\Sigma^{2}=f\left(\Sigma^{1}\right)$, the set of all $f$-images
of the elements of $\Sigma^{1}$; it is easy to prove that $f\left(\Sigma^{1}\right)$
is a sigma-algebra over $f\left(S^{1}\right)$, for any $f$.
\begin{example}
Let $S^{1}=S^{2}=\{1,2,3\}$, 
\[
\Sigma^{1}=\{\emptyset,\{1\},\{2,3\},\{1,2,3\}\},
\]
and 
\[
\Sigma^{2}=\{\emptyset,\{1\},\{2\},\{3\},\{1,2\},\{2,3\},\{1,3\},\{1,2,3\}\}.
\]
Then the function $f:\Sigma^{1}\to\Sigma^{2}$ defined by $f(a)=a$
is not measurable, because $\{2\}\in\Sigma^{2}$ but $f^{-1}(\{2\})=\{2\}\notin\Sigma^{1}$.
However, one can easily verify that $f(a)=\min(a,2)$ is a $\Sigma^{1}\to\Sigma^{2}$-measurable
function. \hfill$\square$ 
\end{example}
Of course, with finite $S^{1},S^{2}$, one can always define the sigma-algebras
as full power sets and then all functions between these sets will
be measurable.

Why is the notion of a measurable function important? Because measurable
functions can be used to obtain new random variables from existing
ones. Given a random variable $A$ and a $\Sigma^{1}\rightarrow\Sigma^{2}\textnormal{-}$measurable
function $f:S^{1}\rightarrow S^{2}$, one can define a random variable
$B=f\left(A\right)$ distributed as $\overline{B}=\left(S^{2},\Sigma^{2},p_{2}\right)$
by putting, for any $S'\in\Sigma^{2}$,
\[
p_{2}\left(S'\right)=p_{1}\left(f^{-1}\left(S'\right)\right).
\]
In other words, the probability with which the new variable $B$ falls
in an event belonging to $\Sigma^{2}$ is defined as the probability
with which $A$ falls in the $f\textnormal{-}$preimage of this event
in $\Sigma^{1}$ (which probability is well defined because $f$ is
measurable). Of course, the notation $B=f\left(A\right)$ serves as
a unique identification of $B$ once we agree that $A$ is uniquely
identified. 
\begin{example}
Let $S^{1}$ and $S^{2}$ be two intervals of reals, and let $\Sigma^{1}$
and $\Sigma^{2}$ be the Borel sigma-algebras over them (see Example
\ref{Example: Real-line-.}). A function $f:S^{1}\rightarrow S^{2}$
which is $\Sigma^{1}\rightarrow\Sigma^{2}$-measurable is called a
\emph{Borel-measurable function}. If in this definition $\Sigma^{1}$
is the Lebesgue sigma algebra over $S^{1}$ while $\Sigma^{2}$ continues
to be the Borel sigma-algebra over $S^{2}$ (note the asymmetry),
then $f$ is a \emph{Lebesgue-measurable function}. It is sufficient
to require in these two definitions that for any interval $\left(a,b\right)\subset S^{2}$,
its preimage $f^{-1}\left(\left(a,b\right)\right)$ be a Borel-measurable
(respectively, Lebesgue-measurable) subset of $S^{1}$. It is easy
to prove that if $f$ is monotone or continuous, then it is Borel-measurable
(hence also Lebesgue-measurable).

Let now $A$ be a random variable with distribution $\overline{A}=\left(\mathbb{R},\Sigma^{1},p\right)$,
where $\Sigma^{1}$ is the Lebesgue sigma-algebra over $\mathbb{R}$.
The function $F\left(x\right)=p\left(\left(-\infty,x\right]\right)$
is called the \emph{distribution function} for $A$. It is monotonically
non-decreasing and maps into $S^{2}=\left[0,1\right]$. If we define
$\Sigma^{2}$ to be the Borel sigma-algebra over $\left[0,1\right]$,
then $F$ (being nondecreasing) is Lebesgue-measurable. If we apply $F$ to $A$, the resulting
random variable $B=F\left(A\right)$ is distributed on $\left[0,1\right]$. If furthermore $F$ is a continuous function, then the distribution of $B=F\left(A\right)$ on $\left[0,1\right]$ is uniform.
That is, its distribution is $\overline{B}=\left(\left[0,1\right],\Sigma^{2},q\right)$,
where $q\left(\left(a,b\right)\right)=b-a$ for any $\left(a,b\right)\subset\left[0,1\right]$.
\hfill$\square$ 
\end{example}
Let $A$ be distributed as $\overline{A}=\left(S^{1},\Sigma^{1},p_{1}\right)$,
and let $B=f\left(A\right)$ and $C=g\left(A\right)$ be two random
variables with distributions $\overline{B}=\left(S^{2},\Sigma^{2},p_{2}\right)$
and $\overline{C}=\left(S^{3},\Sigma^{3},p_{3}\right)$. This implies
that both $f$ and $g$ are measurable functions in the sense of,
respectively, $\Sigma^{1}\rightarrow\Sigma^{2}$ and $\Sigma^{1}\rightarrow\Sigma^{3}$.
For every $S_{B}\in\Sigma^{2}$ and every $S_{C}\in\Sigma^{3}$ we
have 
\[
p_{2}\left(S_{B}\right)=p_{1}\left(f^{-1}\left(S_{B}\right)\right),\textnormal{ and }p_{3}\left(S_{C}\right)=p_{1}\left(g^{-1}\left(S_{C}\right)\right).
\]
A value $b$ of $B$ falls in $S_{B}$ if and only if $b=f\left(a\right)$
for some $a\in f^{-1}\left(S_{B}\right)$. A value $c$ of $C$ falls
in $S_{C}$ if and only if $c=g\left(a\right)$ for some $a\in g^{-1}\left(S_{C}\right)$.
This suggests a way of defining the notion of a \emph{joint occurrence}
of these events, $S_{B}$ and $S_{C}$: they occur jointly if and
only if $a$ in the previous two sentences is one and the same. In
other words, a value $b$ of $B$ falls in $S_{B}$ and, jointly,
a value $c$ of $C$ falls in $S_{C}$ if and only if, for some $a\in f^{-1}\left(S_{B}\right)\cap g^{-1}\left(S_{C}\right)$,
$b=f\left(a\right)$ and $c=g\left(a\right)$. Since $f^{-1}\left(S_{B}\right)\cap g^{-1}\left(S_{C}\right)$
is $\Sigma^{1}\textnormal{-}$measurable in (belongs to $\Sigma^{1}$),
the probability
\[
p_{23}\left(S_{B}\times S_{C}\right)=p_{1}\left(f^{-1}\left(S_{B}\right)\cap g^{-1}\left(S_{C}\right)\right)
\]
is well defined, and we can take it as the joint probability of $S_{B}$
and $S_{C}$. 

We now can construct the joint distribution of $\left(B,C\right)$,
\[
\overline{\left(B,C\right)}=\left(S^{2}\times S^{3},\Sigma^{2}\otimes\Sigma^{3},p_{23}\right),
\]
where the set and the sigma-algebra are defined as required by the
general notion of a joint distribution (Section \ref{sec:Jointly-distributed-random}).
The joint probability measure $p_{23}$ defined above for $S_{B}\times S_{C}$-type
sets is extended to all other members of $\Sigma^{2}\otimes\Sigma^{3}$
by using the basic properties of a probability measure (Section \ref{sec:What-is-a}).
Equivalently, the joint probability measure $p_{23}$ can be defined
by
\[
p_{23}\left(S'\right)=p\left(\left(f,g\right)^{-1}\left(S'\right)\right),
\]
for any $S'\in\Sigma^{2}\otimes\Sigma^{3}$. The notation $\left(f,g\right)^{-1}\left(S'\right)$
designates the set $S_{A}$ of all $a\in S$, such that $\left(f\left(a\right),g\left(a\right)\right)\in S'$.
It can be shown that $S_{A}\in\Sigma^{1}$, that is, $\left(f,g\right)$
is a measurable function.

It can easily be checked that $p_{23}$ satisfies the 1-marginal probability
equations,
\[
\begin{array}{c}
p_{23}\left(S_{B}\times S^{3}\right)=p_{1}\left(f^{-1}\left(S_{B}\right)\cap g^{-1}\left(S^{3}\right)\right)=p_{1}\left(f^{-1}\left(S_{B}\right)\right)=p_{2}\left(S_{B}\right),\\
\\
p_{23}\left(S^{2}\times S_{C}\right)=p_{1}\left(f^{-1}\left(S^{2}\right)\cap g^{-1}\left(S_{C}\right)\right)=p_{1}\left(g^{-1}\left(S_{C}\right)\right)=p_{3}\left(S_{C}\right),
\end{array}
\]
where we used the fact that 
\[
g^{-1}\left(S^{3}\right)=f^{-1}\left(S^{2}\right)=A.
\]

We see that if two random variables are formed as functions of another
random variable, their joint distribution is uniquely determined. 
\begin{example}
\label{Example: identical}A simple but instructive example is the
joint distribution of a random variable $A$ and itself. Let $A$
be distributed as $\left(S,\Sigma,p\right)$. $\left(A,A\right)$
is a random variable both components of which are functions of one
and the same random variable, $A=\mathrm{id}\left(A\right)$, where
$\mathrm{id}$ is the identity function defined by $\mathrm{id}$$\left(a\right)=a$.
Let the distribution of $\left(A,A\right)$ be $\left(S\times S,\Sigma\otimes\Sigma,p_{2}\right)$.
By the general theory, for any $S'\in\Sigma$ we have $S'\times S'\in\Sigma\otimes\Sigma$
and
\[
p_{2}\left(S'\times S'\right)=p\left(\mathrm{id}^{-1}\left(S'\right)\cap\mathrm{id}^{-1}\left(S'\right)\right)=p\left(S'\right),
\]
as it should be. It is not always true, however, that the probability
measure $p_{2}$ of the set of pairs 
\[
\mathrm{diag}_{S\times S}=\left\{ \left(a,a\right):a\in S\right\} 
\]
equals 1, because this set is not necessarily an event in $\Sigma\otimes\Sigma$.
As an example, $\left\{ \left(1,1\right),\left(2,2\right),\left(3,3\right),\left(4,4\right)\right\} $
is not such an event if $\Sigma=\left\{ \emptyset,\left\{ 1,2\right\} ,\left\{ 3,4\right\} ,\left\{ 1,2,3,4\right\} \right\} $.
If, however, $\mathrm{diag}_{S\times S}\in\Sigma\otimes\Sigma$, then
\[
p_{2}\left(\mathrm{diag}_{S\times S}\right)=p\left(\left(\mathrm{id},\mathrm{id}\right)^{-1}\left(\mathrm{diag}_{S\times S}\right)\right)=p\left(S\right)=1.
\]
\hfill$\square$ 
\end{example}
The generalization to several functions of a random variable $A$
is trivial. Thus, we can form a joint distribution not just of $B,C$
but of $A,B,C$ (for symmetry, we can consider $A$ the identity function
of $A$). In particular, the joint probability of $S_{B}\in\Sigma^{2}$,
$S_{C}\in\Sigma^{3}$, and $S_{A}\in\Sigma^{1}$ is defined here as
\[
p_{23}\left(S_{A}\times S_{B}\times S_{C}\right)=p_{1}\left(S_{A}\cap f^{-1}\left(S_{B}\right)\cap g^{-1}\left(S_{C}\right)\right).
\]

One of the important classes of measurable functions of random variables
are projections. We have already dealt with them in Section \ref{sec:Jointly-distributed-random},
when discussing marginal distributions. More generally, a vector of
jointly distributed random variables $A^{1},A^{2},\ldots,A^{n}$ is
a random variable with a distribution
\[
\left(S^{1}\times\ldots\times S^{n},\Sigma^{1}\otimes\ldots\otimes\Sigma^{n},p_{1\ldots n}\right),
\]
where the notation should be clear from the foregoing. A projection
function $\mathrm{Proj}_{i_{1}\ldots i_{k}}$, where $k\leq n$ and
$i_{1},\ldots,i_{k}$ is a set of $k$ distinct numbers chosen from
$\left(1,\ldots,n\right)$, is defined by
\[
\mathrm{Proj}_{i_{1}\ldots i_{k}}\left(a_{1},\ldots,a_{n}\right)=\left(a_{i_{1}},\ldots,a_{i_{k}}\right).
\]
Without loss of generality, let $\left(i_{1},\ldots,i_{k}\right)=\left(1,\ldots,k\right)$;
if this is not the case, one can always make it so by renumbering
the original set of $n$ random variables. The function $\mathrm{Proj}_{1\ldots k}$
creates a $k\textnormal{-}$\emph{marginal random variable} 
\[
\mathrm{Proj}_{1\ldots k}\left(A^{1},\ldots,A^{n}\right)=\left(A^{1},\ldots,A^{k}\right),
\]
with the $k\textnormal{-}$marginal distributions
\[
\left(S^{1}\times\ldots\times S^{k},\Sigma^{1}\otimes\ldots\otimes\Sigma^{k},p_{1\ldots k}\right).
\]
where, for any measurable even $S^{'}$ in $\Sigma^{1}\otimes\ldots\otimes\Sigma^{k}$,
\[
p_{1\ldots k}\left(S'\right)=p_{1\ldots n}\left(S'\times S^{k+1}\times\ldots\times S^{n}\right).
\]

\section{\label{sec:Random-variables-as}Random variables as measurable functions}

We have seen that if $A^{1},\ldots,A^{n}$ are all functions of one
and the same random variable $R$, then they posses a joint distribution.
To recapitulate, if 
\[
A^{1}=f_{1}\left(R\right),\ldots,A^{n}=f_{n}\left(R\right),
\]
\[
\overline{R}=\left(S^{*},\Sigma^{*},p_{*}\right),
\]
and
\[
\overline{A^{i}}=\left(S^{i},\Sigma^{i},p_{i}\right),\; i=1,\ldots,n,
\]
then
\[
\overline{\left(A^{1},\ldots,A^{n}\right)}=\left(S^{1}\times\ldots\times S^{n},\Sigma^{1}\otimes\ldots\otimes\Sigma^{n},p_{1\ldots n}\right),
\]
where
\[
p_{1\ldots n}\left(S'\right)=p_{*}\left(\left(f_{1},\ldots,f_{n}\right)^{-1}\left(S'\right)\right),
\]
for any $S'\in\Sigma^{1}\otimes\ldots\otimes\Sigma^{n}$. In particular,
\[
p_{1\ldots n}\left(S_{1}\times\ldots\times S_{n}\right)=p_{*}\left(\bigcap f_{i}^{-1}\left(S_{i}\right)\right),
\]
for all
\[
S_{1}\in\Sigma^{1},\ldots,S_{n}\in\Sigma^{n}.
\]

It is easy to see that the reverse of this statement is also true:
if $A^{1},\ldots,A^{n}$ have a joint distribution, they can be presented
as functions of one and the same random variable. Indeed, denoting
the random variable $\left(A^{1},\ldots,A^{n}\right)$ by $R$, we
have 
\[
A^{1}=f_{1}\left(R\right),\ldots,A^{n}=f_{n}\left(R\right),
\]
where
\[
f_{i}\equiv\mathrm{Proj}_{i}.
\]
These two simple observations constitute a proof of an important theorem.
\begin{thm}
\label{TH: main}A vector $\left(A^{1},\ldots,A^{n}\right)$ of random
variables possesses a joint distribution if and only if there is a
random variable $R$ and a vector of functions $\left\{ f_{1},\ldots,f_{n}\right\} $,
such that $A^{1}=f_{1}\left(R\right),\ldots,A^{n}=f_{n}\left(R\right)$.
\end{thm}
Note that we need not specify here that the functions are measurable,
because both $A^{i}$ and $R$ in $A^{i}=f_{i}\left(R\right)$ are
random variables (implying that $f_{i}$ is measurable). 

Although we do not deal in this chapter with infinite sets of jointly
distributed random variables, it must be mentioned that Theorem \ref{TH: main}
has the following \emph{generalized formulation} (see Remark \ref{REM: For-an-infinite}). 
\begin{thm}
\label{TH: main generalized}A family $\left(A^{k}:k\in K\right)$
of random variables possesses a joint distribution if and only if
there is a random variable $R$ and a family of functions $\left(f_{k}:k\in K\right)$
such that $A^{k}=f_{k}\left(R\right)$ for all $k\in K$.
\end{thm}
In probability textbooks, consideration is almost always confined
to random variables that are jointly distributed. This enables what
we may call the \emph{traditional conceptualization of random variables}.
It consists in choosing some distribution 
\[
\overline{R}=\left(S^{*},\Sigma^{*},p_{*}\right),
\]
calling it a \emph{sample (probability) space}, and identifying any
random variable $A$ as a ($\Sigma^{*}\rightarrow\Sigma^{1}$)-measurable
function $f:S^{*}\rightarrow S^{1}$. The set and sigma-algebra pair
$\left(S^{1},\Sigma^{1}\right)$ being chosen, the probability measure
$p_{1}$ satisfying, for every $S'\in\Sigma^{1}$,
\[
p_{1}\left(S'\right)=p_{*}\left(f^{-1}\left(S'\right)\right),
\]
is referred to as an i\emph{nduced probability measure}, and the distribution
$\overline{A}=\left(S^{1},\Sigma^{1},p_{1}\right)$ as an \emph{induced
(probability) space}. 

The sample space $\overline{R}$ is the distribution of some random
variable $R$; in the language just presented $R$ should be defined
as the identity function $\mathrm{id}:S^{*}\rightarrow S^{*}$ (one
that maps each element into itself) on the sample space $\overline{R}$;
its induced probability space is, obviously, also $\overline{R}$.
In our conceptual framework we simply define $R$ by its distribution
$\overline{R}$ and some unique identifying label (such as ``$R$'').
Note that the traditional language, too, requires an identifying label
and a distribution (using our terminology) in order to define the
sample space itself. 
\begin{rem}
The traditional language does not constitute a different approach.
It is a terminological variant of the conceptual set-up adopted in
this chapter and applied to a special object of study: a class $\mathcal{A}$
of random variables that can be defined as functions of some ``primary''
random variable $R$. In accordance with Theorem \ref{TH: main generalized},
$\mathcal{A}$ can also be described without mentioning $R$, as a
class of random variables such that, for any indexed family of random
variables $\left(A^{k}:k\in K\right)$ with $A^{k}\in\mathcal{A}\left(R\right)$
for all $k\in K$, there is a random variable $A=\left(A^{k}:k\in K\right)$
that also belongs to $\mathcal{A}$. 
\end{rem}

\section{\label{sec:Stochastically-unrelated-and}Unrelated random variables
and coupling schemes}

There are two considerations to keep in mind when using the traditional
language of random variables as measurable functions on sample spaces. 

One of them is that sample spaces $\overline{R}$ (or ``primary''
random variables $R$) are more often than not nebulous: they need
not be and usually are not explicitly introduced when dealing with
collections of jointly distributed random variables, and they often
have no substantive interpretation if introduced. Consider an experiment
in which a participant is shown one of two stimuli, randomly chosen,
and is asked to identify them by pressing one of two keys as soon
as possible. In each trial we record two random variables: stimulus
presented and response time observed, RT. The joint distribution of
stimuli and response times is well defined by virtue of pairing them
trial-wise. But what would the ``primary'' random variable $R$
be of which stimulus and RT would be functions? No one would normally
attempt determining one, and it is difficult if one tries, except
for the trivial choice $R=$(stimulus, RT) or some one-to-one function
thereof. The stimulus and RT then would be projections (i.e., functions)
of $R$, but this hardly adds insights to our understanding of the
situation. Moreover, as soon as one introduces a new random variable
in the experimental design, say, ``response key,'' indicating which
of the two keys was pressed, the ``primary'' random variable $R$
has to be redefined. It may now be the jointly distributed triple
$R=$(stimulus, response key, RT). 

The second consideration is that there can be no such thing as a single
``primary'' random variable $R$ allowing one to define all conceivable
random variables as its functions. This is obvious from the cardinality
considerations alone: the set $S^{*}$ in $\overline{R}$ would have
to be ``larger'' than the set of possible values for any conceivable
random variable (which can, of course, be chosen arbitrarily large).
It is a mathematical impossibility. The universe of all conceivable
random variables should necessarily include random variables that
are not functions of a common ``primary'' one. In view of Theorem
\ref{TH: main generalized}, this means that there must be random
variables that do not possess a joint distribution. The situation
should look like in the diagram below, with $A^{1},A^{2},\ldots$
being functions of some $R^{1}$, $B^{1},B^{2},\ldots$ being functions
of some $R^{2}$, but $R^{1}$ and $R^{2}$ being \emph{stochastically
unrelated}, with no joint distribution.

\medskip{}

\begin{footnotesize}

\[
\xymatrix{ &  & R^{1}\ar[d]\ar[dl]\ar[dr]\ar[dll]\ar[drr]\\
\ldots & A^{1} & A^{2} & A^{3} & \ldots
}
\quad\xymatrix{ &  & R^{2}\ar[d]\ar[dl]\ar[dr]\ar[dll]\ar[drr]\\
\ldots & B^{1} & B^{2} & B^{3} & \ldots
}
\]
\end{footnotesize}

\medskip{}

\protect{\noindent}It is true that, as explained below, once $R^{1}$
and $R^{2}$ are introduced (by their distributions and identifying
labels), there is always a way to introduce a new random variable
$\left(H^{1},H^{2}\right)$ (whose components are functions of some
random variables) such that $H^{1}$ has the same distribution as
$R^{1}$ and $H^{2}$ has the same distribution as $R^{2}$. But there
is no way of conceiving all random variables in the form of functions
of a single ``primary'' one.

Examples of random variables that normally are not introduced as jointly
distributed are easy to find. If RTs in an experiment with two stimuli
(say, ``green'' and ``red'') are considered separately for stimulus
``green'' and stimulus ``red'', we have two random variables:
RT$^{green}$ and RT$^{red}$. What ``natural'' stochastic relationship
they might have? The answer is, none: the two random variables occur
in mutually exclusive conditions, so there is no privileged way of
coupling realizations of RT$^{green}$ and RT$^{red}$ and declaring
them co-occurring. Once these random variables are introduced, one
can impose a joint distribution on them. For example, one may consider
them stochastically independent, essentially forcing on them the coupling
scheme in which each realization of RT$^{green}$ considered as if
it co-occurred with every realization RT$^{red}$. But it is also
possible to couple them differently, for instance, by the common quantile
ranks, so that the $q$th quantile of RT$^{red}$ is paired with and
only with the $q$th quantile of RT$^{green}$. The two random variables
then are functions of the quantile rank, which is a random variable
uniformly distributed between 0 and 1. The point is, neither of these
nor any of the infinity of other coupling schemes for the realizations
of RT$^{green}$ and RT$^{red}$ is privileged, and none is necessary:
one need not impose any joint distribution on RT$^{green}$ and RT$^{red}$.

It can be shown that stochastic independence can be imposed on any
set of pairwise \emph{stochastically unrelated} random variables.
\begin{thm}
\label{TH: universality of ind coupling}For any vector $\left(R^{1},\ldots,R^{n}\right)$
(more generally, any family $\left(R^{k}:k\in K\right)$) of random
variables that are pairwise stochastically unrelated there is a random
variable $H=\left(H^{1},\ldots,H^{n}\right)$ (generally, $H=\left(H^{k}:k\in K\right)$)
with stochastically independent $H^{k}$, such that $\overline{H^{k}}=\overline{R^{k}}$
for all $k\in K$. 
\end{thm}
$H$ is called the \emph{independent coupling} for $\left(R^{k}:k\in K\right)$.
In general, a \emph{coupling} for a family of random variables $\left(R^{k}:k\in K\right)$,
is any random variable $H=\left(H^{k}:k\in K\right)$ whose every
1-marginal random variable $H^{k}$ is distributed as $R^{k}$. 

Theorem \ref{TH: universality of ind coupling} must not be interpreted
to mean that one can take all pairwise stochastically unrelated random
variables and consider them stochastically independent. The reason
for this is that this class is not a well defined set, and cannot
be therefore indexed by any set. Indeed, if it were possible to present
it as $\left(R^{k}:k\in K\right)$, then one could form a new random
variable $R=\left(R^{k}:k\in K\right)$ whose distribution is the
same as $\overline{\left(H^{k}:k\in K\right)}$ in Theorem \ref{TH: universality of ind coupling},
and it would follow that the set contains itself as an element (which
is impossible for a set).

Summarizing, in practice random variables are often well defined without
their joint distribution being well defined. There is nothing wrong
in dealing with stochastically unrelated random variables without
trying to embed them in jointly distributed system. When such an embedding
is desirable, the joint distribution is ``in the eyes of the beholder,''
in the sense of depending on how one wishes to couple the realizations
of the variables being interrelated.

\section{On sameness, equality, and equal distributions}

We have to distinguish two different meanings in which one can understand
the equality of random variables, $A=B$. 

One meaning is that $A$ and $B$ are different notations for one
and the same variable, that is, that $A$ and $B$ have the same identifying
label and the same distribution. This meaning of equality is implicit
when we say ``let $D$ be $\left(A,B,C\right)$, jointly distributed''
or ``there is a random variable $A=\left(A^{k}:k\in K\right)$.'' 

The other meaning of $A=B$ is that 
\begin{enumerate}
\item these random variables have (or may have) different identifying labels
(i.e., they are not or may not be the same);
\item they are identically distributed, $\overline{A}=\overline{B}=\left(S,\Sigma,p_{1}\right)$;
\item they are jointly distributed, and their joint distribution has the
form $\left(S\times S,\Sigma\otimes\Sigma,p_{2}\right)$;
\item for any $S'\in\Sigma$, 
\[
p_{2}\left(S'\times S'\right)=p_{1}\left(S'\right).
\]

\end{enumerate}
In some cases, if $\mathrm{diag}_{S}=\left\{ \left(a,a\right):a\in S\right\} $
is a measurable set (i.e., it belongs to $\Sigma\otimes\Sigma$), one
can replace the last property with 
\[
p_{2}\left(\mathrm{diag}_{S}\right)=1,
\]
which can also be presented as
\[
\Pr\left(A=B\right)=1.
\]
If $A$ and $B$ about which we know that $A=B$ are represented as
functions of some random variable $R$, then it is usually assumed
that $\mathrm{diag}_{S}\in\Sigma\otimes\Sigma$, and the two functions
representing $A$ and $B$ are called \emph{equal with probability
1} (or \emph{almost surely}). Of course, if $A$ and $B$ are merely
different notations for one and the same random variable, they are
always jointly distributed and equal in the second sense of the term
(see Example \ref{Example: identical}).

The equality of random variables, in either sense, should not be confused
with the equality of distributions, $\overline{A}=\overline{B}$.
The random variables $A$ and $B$ here may but do not have to be
jointly distributed. They may very well be stochastically unrelated.
We will use the symbol $\sim$ in the meaning of ``has the distribution''
or ``has the same distribution as.'' Thus, $A\sim\overline{A}$
always, $A\sim B$ if and only if $\overline{A}=\overline{B}$, and
$A=B$ always implies $A\sim B$.

An important notational consideration applies to random variables
with imposed on them or redefined joint distributions. One may write
$\left(A,B\right)$ either as indicating a pair of stochastically
unrelated random variables, or some random variable $C=\left(A,B\right)$.
The two meanings are distinguished by context. Nothing prevents one,
in principle, from considering the same $A$ and $B$ as components
of two differently distributed pairs, $C=\left(A,B\right)$ and $C'=\left(A,B\right)$,
or as components of a $C=\left(A,B\right)$ possessing a joint distribution
and a pair $\left(A,B\right)$ of stochastically unrelated random
variables. Doing this within the same context, however, will create
conceptual difficulties. For one thing, we would lose the ability
of presenting $A$ and $B$ as functions of some $R$ (based on their
joint distribution in $C$). 

There is a simple and principled way of avoiding this inconvenience:
use different symbols for random variables comprising different pairs
(more generally, vectors or indexed families), considering them across
the pairs (vectors, families) as equally distributed stochastically
unrelated random variables. In our example, we can write $C=\left(A,B\right)$
and $C'=\left(A',B'\right)$, where $A\sim A'$ and $B\sim B'$, with
$C$ and $C'$ being stochastically unrelated. The same principle
was applied in the formulation of Theorem \ref{TH: universality of ind coupling}
and more generally, in the definition of a coupling: rather than saying
that given a family of stochastically unrelated $\left(R^{k}:k\in K\right)$,
its coupling is any random variable $H=\left(R^{k}:k\in K\right)$
whose components are jointly distributed (e.g., independent), the
definition says that a coupling is a random variable $H=\left(H^{k}:k\in K\right)$
such that $H^{k}\sim R^{k}$ for all $k\in K$. This means, in particular,
that every vector of random variables is stochastically unrelated
to any of its couplings.

\section{\label{sec:Random-outputs-depending}Random outputs depending on
inputs}

Let a random variable be distributed as $\left(S,\Sigma,p_{\phi}\right)$,
where $\phi$ stands for some deterministic variable taking values
in a set $\Phi$. This means that the probability measure on $\Sigma$
(the entire function) is generally different for different values
of $\Phi$. One could also write $p\left(\phi\right)$ instead of
$p_{\phi}$, but one should keep in mind that this is not a function
from $\Phi$ to a set of values of $p$ (real numbers between 0 and
1) but rather a function from $\Phi$ to the set of all possible probability
measures on $\Sigma$. The dependence of $p_{\phi}$ on $\phi$ means
that the distribution $\left(S,\Sigma,p_{\phi}\right)$ of the random
variable in question depends on $\phi$. We can present it as $\overline{A_{\phi}}$,
and the random variable itself as $A_{\phi}$. One can say that the
random variable $A$ \emph{depends on} $\phi$, which is equivalent
to saying that there is an indexed family of random variables $\left(A_{\phi}:\phi\in\Phi\right)$.

Let $\phi_{1}$ and $\phi_{2}$ be two different elements of $\Phi$.
We will assume throughout the rest of the chapter that the corresponding
random variables $A_{\phi_{1}}$ and $A_{\phi_{2}}$ always have different
identifying labels (such as ``$A$ at $\phi=\phi_{1}$'' and ``$A$
at $\phi=\phi_{2}$''), that is, they are never one and the same
variable. But they may have one and the same distribution function,
if $p_{\phi_{1}}\equiv p_{\phi_{2}}$. If $A$ is a vector of jointly
distributed random variables $\left(A^{1},\ldots,A^{n}\right)$, then
its dependence on $\phi$ can be shown as $A_{\phi}=\left(A^{1},\ldots,A^{n}\right)_{\phi}$
or $A_{\phi}=\left(A_{\phi}^{1},\ldots,A_{\phi}^{n}\right)$.

In the following, $\phi$ always represents mutually exclusive conditions
under which $A$ is observed, and the indexed family $\left(A_{\phi}:\phi\in\Phi\right)$
abbreviated by $A$ consists of pairwise stochastically unrelated
random variables. The elements of $\Phi$ are referred to as \emph{treatments},
the term being used in the same way as in the analysis of variance:
a combination of values of different \emph{factors}, or \emph{inputs}.
We will use the latter term. An input is simply a variable $\lambda$
with a set of possible values $\Lambda$. If the number of inputs
considered is $m$, a treatment is a vector 
\[
\phi=\left(\lambda^{1},\ldots,\lambda^{m}\right),
\]
with $\lambda^{1}\in\Lambda^{1}$, ..., $\lambda^{m}\in\Lambda^{m}$.
The set of treatments is therefore 
\[
\Phi\subset\Lambda^{1}\times\ldots\times\Lambda^{m}.
\]
 
\begin{rem}
\label{REM: As-it-is}As it is commonly done in mathematics, we will
use the same symbol to denote a variable and its specific values.
For example, in $\lambda^{1}\in\Lambda^{1}$ the symbol $\lambda^{1}$
refers to a value of $\lambda^{1}$, whereas in the sentence ``$A^{1}$
depends on $\lambda^{1}$'' the same symbol refers to the variable
as a whole. This ambiguity is possible to avoid by using $\Lambda^{1}$
in place of $\lambda^{1}$ when referring to the entire variable,
and using a pair $\left(\lambda^{1},\Lambda^{1}\right)$ when referring
to an input value as that of a given input. We do not use this rigorous
notation here, assuming context will be sufficient for disambiguation.\end{rem}
\begin{example}
Let $\phi$ describe a stimulus presented to a participant. Let it
attain eight possible values formed by combinations of three binary
attributes, such as 
\[
\lambda^{1}\in\Lambda^{1}=\left\{ \text{large},\text{small}\right\} ,\lambda^{2}\in\Lambda^{2}=\left\{ \text{bright},\text{dim}\right\} ,\lambda^{3}\in\Lambda^{3}=\left\{ \text{round},\text{square}\right\} .
\]
Let the participant respond by identifying (correctly or incorrectly)
these attributes, by saying $A^{1}={}$``large'' or ``small'',
$A^{2}={}$``bright'' or ``dim'', and $A^{3}={}$``round'' or
``square''. The response therefore is a vector of three binary random
variables $\left(A^{1},A^{2},A^{3}\right)_{\phi}$ that depends on
stimuli $\phi=\left(\lambda^{1},\lambda^{2},\lambda^{3}\right)$.
Equivalently, we can say that there are eight triples of random variables,
one for each treatment, $\left(A^{1},A^{2},A^{3}\right)_{\phi_{1}},\ldots$,
$\left(A^{1},A^{2},A^{3}\right)_{\phi_{8}}$. \hfill$\square$ 
\end{example}
The set of all treatments $\Phi$ may be equal to $\Lambda^{1}\times\ldots\times\Lambda^{m}$,
but it need not be. Some of the logically possible combinations of
input values may not be physically realizable or simply may not be
of interest. The elements of $\Phi$ therefore are referred to as
\emph{allowable treatments}. We will see later that this notion is
important in pairing inputs with random outputs. 
\begin{example}
Suppose $\Lambda^{1}$ and $\Lambda^{2}$ denote the sets of possible
lengths of two line segments presented side by side in the visual
field of an observer. Let $A^{1}$ and $A^{2}$ denote the observer's
numerical estimates of the two lengths. If the goal of the experiment
is to study perceptual discrimination, it may be reasonable (and time-saving)
to exclude the pairs with large values of $\left|\lambda^{1}-\lambda^{2}\right|$.
For example, if $\Lambda^{1}=\Lambda^{2}=\{5,6,7,8,9\}$, the set
of allowable treatments may be defined as 
\[
\Phi=\{(\lambda^{1},\lambda^{2})\in\Lambda^{1}\times\Lambda^{2}:|\lambda^{1}-\lambda^{2}|\le2\}.
\]
This set contains only 19 treatments of the 25 logically possible
combinations.\hfill$\square$ 
\end{example}
As explained in the introductory section, inputs may very well be
random variables themselves, but only their possible values rather
than their distributions are relevant in our analysis: the distributions
of random outputs are always \emph{conditioned} upon particular treatments.
All inputs therefore are always treated as deterministic quantities.

\section{\label{sec:Selectiveness-in-the}Selectiveness in the dependence
of outputs on inputs}

We are interested in the relationship between (deterministic) inputs
and random outputs. Specifically, we are interested in the selectiveness
in this relationship: which input may and which may not influence
a given output. Such selectiveness can be presented in the form of
a \emph{diagram} \emph{of influences}, where an arrow from an input
$\lambda$ to a random output $A$ means that $\lambda$ \emph{influences}
$A$ (note that the meaning of ``influence'' has not been as yet
defined). The absence of an arrow from an input $\lambda$ to a random
output $A$ excludes $\lambda$ from the set of inputs that influence
$A$.

Consider, for example the following arrow diagram

\emph{
\[
\xymatrix{\alpha\ar[d]\ar[dr] & \beta\ar[d]\ar[dl]\ar[dr]\ar[drr] & \gamma\ar[dr]\ar[dl]\\
A & B & C & D
}
\]
}This diagram can be interpreted by saying that: 
\begin{enumerate}
\item the random outputs $\left(A,B,C,D\right)$ are jointly distributed,
and their joint distribution (specifically, joint probability measure)
depends on the inputs $\left(\alpha,\beta,\gamma\right)$; in other
words, $\left(A,B,C,D\right)$ is in fact $\left(A,B,C,D\right)_{\alpha\beta\gamma}$,
or $\left(A_{\alpha\beta\gamma},B_{\alpha\beta\gamma},C_{\alpha\beta\gamma},D_{\alpha\beta\gamma}\right)$.
\item output $A$ is influenced by inputs $\alpha,\beta$ but not by input
$\gamma$;
\item output $B$ is influenced by all inputs, $\alpha,\beta,\gamma$;
\item output $C$ is influenced by input $\beta$ but not by inputs $\alpha,\gamma$;
\item output $D$ is influenced by inputs $\beta$ and $\gamma$, but not
by $\alpha$.
\end{enumerate}
The first thing to do here is to ask the question we asked in the
introductory section: does this even make sense? It certainly does
if $\left(A,B,C,D\right)_{\alpha\beta\gamma}$, for every treatment
$\left(\alpha,\beta,\gamma\right)$, is a vector of \emph{independent}
random variables. Then the points 2,3, and 4, above simply translate
into the statements: the marginal distribution of $A$ depends on
$\alpha,\beta$ but not on $\gamma$; the marginal distribution of
$B$ depends on $\alpha,\beta,\gamma$; etc. But does the selectiveness
make sense if the random outputs are not stochastically independent?
Look at the diagram below, the same as above, but with added point
lines indicating stochastic interdependences.

\emph{
\[
\xymatrix{\alpha\ar[d]\ar[ddr] & \beta\ar[dd]\ar[dl]\ar[ddr]\ar[drr] & \gamma\ar[dr]\ar[ddl]\\
A\ar@{:}[dr]\ar@{:}[drr]\ar@{:}[rrr] &  &  & D\\
 & B\ar@{:}[r]\ar@{:}[urr] & C\ar@{:}[ur]
}
\]
}We see, for instance, that output $A$ is influenced by $\alpha$,
and output $C$ is stochastically dependent on $A$. In what sense
then one can say that $\alpha$ does not influence $C$? The output
$B$ is influenced by all inputs, and every other output is stochastically
dependent on $B$. Does not this mean that every output is influenced
by every input?

This seemingly compelling line of reasoning is a conceptual confusion.
It confuses two types of relations, both of which can be described
using the word ``dependence.'' Stochastic dependence and dependence
of outputs on inputs are different in nature. This is easy to understand
if we consider the following diagram:

\[
\xymatrix{\alpha\ar[d]\ar[dr] & \beta\ar[d]\ar[dl]\ar[drr]\ar[drrr] &  & \gamma\ar[dr]\ar[dll]\\
A' & B' &  & C' & D'\\
 &  & R\ar[ull]\ar[ul]\ar[ur]\ar[urr]
}
\]
In this diagram, every random variable is a function of all the arguments
from which the arrows leading to this random variable initiate:
\[
\begin{array}{l}
A'_{\alpha\beta\gamma}=f_{1}\left(\alpha,\beta,R\right),\\
B'_{\alpha\beta\gamma}=f_{2}\left(\alpha,\beta,\gamma,R\right),\\
C'_{\alpha\beta\gamma}=f_{3}\left(\beta,R\right),\\
D'_{\alpha\beta\gamma}=f_{4}\left(\alpha,\beta,R\right).
\end{array}
\]
For every value of $R$ and for every treatment $\left(\alpha,\beta,\gamma\right)$,
the values of $\left(A',B',C',D'\right)_{\alpha\beta\gamma}$ are
determined uniquely. Suppose now that we have, for every treatment,
\[
\left(A',B',C',D'\right)_{\alpha\beta\gamma}\sim\left(A,B,C,D\right)_{\alpha\beta\gamma}.
\]
This assumption explains the coexistence of the stochastic relationship
between the random outputs and the selectiveness in their dependence
on the inputs. \emph{For} \emph{any given treatment}, the components
of $\left(A,B,C,D\right)_{\alpha\beta\gamma}$ are generally stochastically
interdependent because they are distributed as functions of one and
the same random variable $R$ (of course, as a special case, they
may also be stochastically independent). At the same time, \emph{for
any fixed value $r$ of $R$}, the value $a=f_{1}\left(\alpha,\beta,r\right)$
of the output $A'_{\alpha\beta\gamma}$ cannot depend on $\gamma$,
the value $c=f_{3}\left(\beta,r\right)$ of the output $C'_{\alpha\beta\gamma}$
cannot depend on anything but $\beta$, etc. And since the distributions
of $\left(A',B',C',D'\right)_{\alpha\beta\gamma}$ and $\left(A,B,C,D\right)_{\alpha\beta\gamma}$
are the same, they share the same selectiveness pattern. 

This consideration leads us to a rigorous definition of what it means
for a vector of random outputs $\left(A,B,C,D\right)_{\alpha\beta\gamma}$
to satisfy the pattern of selective influences represented in the
opening diagram of this section: this pattern is satisfied if and
only if the equations above are satisfied for some choice of a random
variable $R$ and function $f_{1},f_{2},f_{3},f_{4}$. This definition
can be generalized to an arbitrary family of random outputs and an
arbitrary family of inputs. However, we will confine our attention
to the case when these families are finite vectors. And we will use
a special (re-)arrangement of the inputs to make the definition especially
simple. 
\begin{rem}
It should be kept in mind that the meaning of ``$\lambda$ influences
$A$'' includes, as a special case the possibility of $\lambda$
not influencing $A$. There is an asymmetry in saying that, in the
example used in this section, $C$ depend on $\beta$, and saying
that $C$ does not depend on $\alpha$. The latter is a definitive
statement: $\alpha$ is not within the list of arguments in the function
$c=f_{3}\left(\beta,r\right)$. The dependence on $\beta$ means that
$\beta$ is within this list. But a constant function is a special
case of a function. So $c=f_{3}\left(\beta,r\right)$ may, as a special
case, be constant at all values of $R$, or at all values of $R$
except on a subset of measure zero. For instance, if $R$ is uniformly
distributed between 0 and 1 (we will see below that this choice is
possible in a wide class of cases) and $c=f_{3}\left(\beta,r\right)$
is a non-constant function of $\beta$ only at rational $r$, then
$C$ does not depend on $\beta$ with probability 1 (because the set
of all rational points is countable, hence its Lebesgue measure is
zero). This shows that the terms ``depends on'' and ``influences''
should generally be understood as ``may depend on'' and ``may influence.''
\end{rem}

\section{\label{sec:Canonical-form}Selective Influences in a canonical form}

Continuing with the same example, let us consider the random outputs
one by one, and for each of them group together all inputs that influence
it. We get 
\[
\xymatrix{\lambda^{1}=\left(\alpha,\beta\right)\ar[d] & \lambda^{2}=\left(\alpha,\beta,\gamma\right)\ar[d] & \lambda^{3}=\left(\beta\right)\ar[d] & \lambda^{4}=\left(\beta,\gamma\right)\ar[d]\\
A & B & C & D
}
\]
Let us assume that each of the inputs $\alpha,\beta,\gamma$ has three
possible values, crossed in all possible ways to form 27 treatments.
Each of the newly formed groups of inputs can be viewed as a new input
in its own right. Thus, $\lambda^{1}$ and $\lambda^{4}$ are inputs
whose sets of possible values $\Lambda^{1}$ and $\Lambda^{4}$ have
nine possible values each, $\lambda^{2}$ is an input with 27 possible
values in $\Lambda^{2}$, and $\lambda^{3}$ is an input with three
values in $\Lambda^{3}$. 

Such a rearrangement is always possible, whatever the original pattern
of influences, and it achieves a one-to-one correspondence between
random outputs and inputs. We call a diagram with such one-to-one
correspondence a \emph{canonical diagram of influences}. (The term
``canonical'' is used in mathematics to refer to a standard representation
into which a variety of other representations can be transformed.)
The problem of selectiveness with a canonical diagram acquires a simple
form: is every random output selectively influenced by its corresponding
input?

When dealing with canonical diagrams it is especially important to
keep in mind that \emph{allowable treatments} are generally just a
subset of the Cartesian product of the sets of input values. In our
example, this Cartesian product is $\Lambda^{1}\times\Lambda^{2}\times\Lambda^{3}\times\Lambda^{4}$
and it consists of $9\times27\times3\times9$ elements. But, obviously,
only 27 combinations of new inputs' values are allowable, corresponding
to the 27 treatments formed by the completely crossed original inputs.
Thus, if $\lambda^{2}=\left(\alpha,\beta,\gamma\right)$, then the
only allowable treatment containing this value of $\lambda^{2}$ also
contains $\lambda^{1}=\left(\alpha,\beta\right)$, $\lambda^{3}=\left(\beta\right)$,
and $\lambda^{4}=\left(\beta,\gamma\right)$. 

Another consideration related to the canonical diagrams of influences
is that in order to ensure one-to-one correspondence between inputs
and random outputs, we may need to allow for ``dummy'' inputs, with
a single possible value. Consider the following example: 

\emph{
\[
\xymatrix{\alpha\ar[d] & \beta\ar[d]\ar[dl] & \gamma\ar[dl]\\
A & B & C
}
\]
}Not being influenced by any inputs (as it is the case with the output
$C$) is a special case of selectiveness, so this situation falls
within the scope of our analysis. Presented in the canonical form,
this diagram becomes
\[
\xymatrix{\lambda^{1}=\left(\alpha,\beta\right)\ar[d] & \lambda^{2}=\left(\beta,\gamma\right)\ar[d] & \lambda^{3}=\left(\right)\ar[d]\\
A & B & C
}
\]
The new input $\lambda^{3}$ represents an empty subset of original
inputs. Therefore $\lambda^{3}$ does not change, and should formally
viewed as an input whose set of possible values $\Lambda^{3}$ contains
a single element, that we may denote arbitrarily.

We are ready now to give a formal definition of selective influences.
Let $\left(\lambda^{1},\ldots,\lambda^{n}\right)$ be a vector of
inputs, with values belonging to nonempty sets $\left(\Lambda^{1},\ldots,\Lambda^{n}\right)$,
respectively. Let $\Phi\subset\Lambda^{1}\times\ldots\times\Lambda^{n}$
be a nonempty set of allowable treatments. Let $\left(A_{\phi}^{1},\ldots,A_{\phi}^{n}\right)$
be a vector of random variables jointly distributed for every $\phi\in\Phi$.
(Recall that for $\phi\not=\phi'$, the random variables $\left(A_{\phi}^{1},\ldots,A_{\phi}^{n}\right)$
and $\left(A_{\phi'}^{1},\ldots,A_{\phi'}^{n}\right)$ are stochastically
unrelated.) We say that the dependence of $\left(A_{\phi}^{1},\ldots,A_{\phi}^{n}\right)$
on $\phi$ satisfies the (canonical) diagram of influences 
\[
\xymatrix{\lambda^{1}\ar[d] & \ldots & \lambda^{n}\ar[d]\\
A^{1} & \ldots & A^{n}
}
\]
if and only if one can find a random variable $R$ and functions $f_{1},\ldots,f_{n}$
such that
\[
\left(A_{\phi}^{1},\ldots,A_{\phi}^{n}\right)\sim\left(f_{1}\left(\lambda^{1},R\right),\ldots,f_{n}\left(\lambda^{n},R\right)\right)
\]
for every $\left(\lambda^{1},\ldots,\lambda^{n}\right)=\phi\in\Phi.$ 
\begin{rem}
There is no implication of uniqueness in this definition: below, in
the discussion of the linear feasibility test, we will reconstruct
$R$ explicitly, and we will see that it can, as a rule, be chosen
in infinitely many ways. Theorem \ref{TH: Under-the-conditions} below
shows the non-uniqueness of $R$ by another argument.
\end{rem}
Instead of drawing diagrams, in the sequel we will present the same
pattern of selective influences as
\[
\left(A^{1},\ldots,A^{n}\right)\looparrowleft\left(\lambda^{1},\ldots,\lambda^{n}\right),
\]
and say that $A^{1},\ldots,A^{n}$ are \emph{selectively influenced}
by $\lambda^{1},\ldots,\lambda^{n}$ (respectively). If it is known
that for a given vector of input-output pairs the definition above
is not satisfied whatever $R$ and $f_{1},\ldots,f_{n}$ one chooses,
then we write
\[
\left(A^{1},\ldots,A^{n}\right)\not\looparrowleft\left(\lambda^{1},\ldots,\lambda^{n}\right).
\]
Note that for this schematic notation to make sense, context in which
it is used should specify the sets of input values, the distributions
of $\left(A_{\phi}^{1},\ldots,A_{\phi}^{n}\right)$, and the set of
allowable treatments. 
\begin{example}
Let $R=(R_{1},R_{2},R_{3})$ denote a vector of independent standard
normal random variables, and suppose the input factors $\Lambda^{1}$
and $\Lambda^{2}$ are some subsets of $\mathbb{R}$. Then, the binary
random variables
\begin{align*}
A_{(\lambda^{1},\lambda^{2})}^{1} & =\begin{cases}
1 & \textnormal{if }R_{1}<\lambda\text{\textonesuperior}+R_{3},\\
0 & \text{otherwise},
\end{cases}\\
A_{(\lambda^{1},\lambda^{2})}^{2} & =\begin{cases}
1 & \textnormal{if }R_{2}<\lambda^{2}+R_{3},\\
0 & \text{otherwise},
\end{cases}
\end{align*}
are selectively influenced by respectively $\lambda^{1}\in\Lambda^{1}$
and $\lambda^{2}\in\Lambda^{2}$, because $A^{1}$ depends only on
$(\lambda^{1},R)$ and $A^{2}$ depends only on $(\lambda^{2},R)$.
For any given $\left(\lambda^{1},\lambda^{2}\right)$, the random
variables $A_{(\lambda^{1},\lambda^{2})}^{1}$ and $A_{(\lambda^{1},\lambda^{2})}^{2}$
are not stochastically independent because $R_{1}-R_{3}$ and $R_{2}-R_{3}$
have a nonzero correlation.\hfill$\square$ 
\end{example}

\section{\label{sec:Joint-Distribution-Criterion}Joint Distribution Criterion}

Let us begin by making sure that the simplest special case, when $\left(A_{\phi}^{1},\ldots,A_{\phi}^{n}\right)$
are mutually independent random variables at every allowable treatment
$\phi$, falls within the scope of the general definition. We expect,
if our general definition is well constructed, that in this case selectiveness
of influences, $\left(A^{1},\ldots,A^{n}\right)\looparrowleft\left(\lambda^{1},\ldots,\lambda^{n}\right)$,
follows from the fact that the distribution of $A_{\phi}^{k}$ (for
$k=1,\ldots,n$) depends only on $\lambda^{k}$. In order not to deal
with infinite indexed families, let us assume that $\lambda^{k}$
has a finite number of values, enumerated as $1,\dots,m_{k}$. Consider
the random variable
\[
H=\left(H_{1}^{1},\ldots,H_{m_{1}}^{1},\ldots,H_{1}^{k},\ldots,H_{m_{k}}^{k},\ldots,H_{1}^{n},\ldots,H_{m_{n}}^{n}\right)
\]
with stochastically independent components, such that, for all $i=1,\dots,m_{k}$
and $k=1,\ldots,n$,
\[
H_{i}^{k}\sim A_{\phi}^{k}
\]
whenever $\lambda^{k}=i$ is in $\phi$. In other words, if the treatment
$\phi$ contains the $i$th value of the input $\lambda^{k}$, then
we pick $A_{\phi}^{k}$, and change its identifying label with its
distribution intact to create $H_{i}^{k}$. Clearly, the $H_{i}^{k}$
will be the same (provided we always use the same label) irrespective
of which $\phi$ contains $\lambda^{k}=i$. The variable $H$ above
always exists by Theorem \ref{TH: universality of ind coupling}.
Let us define function $f_{k}$ for $k=1,\ldots,n$ by
\[
f_{k}\left(i,h_{1}^{1},\ldots,h_{m_{k}}^{1},\ldots,h_{1}^{k},\ldots,h_{m_{k}}^{k},\ldots,h_{1}^{n},\ldots,h_{m_{k}}^{n}\right)=h_{i}^{k}.
\]
This can be understood as the ``first-level'' $k$th projection
that selects from the range of the arguments the subrange $h_{1}^{k},\ldots,h_{m_{k}}^{k}$,
followed by the ``second-level'' $i$th projection that selects
from this subrange the argument $h_{i}^{k}$. It is obvious then that,
for every $\phi\in\Phi$,
\[
A_{\phi}^{k}\sim f_{k}\left(i,H\right)
\]
whenever $\phi$ contains $\lambda^{k}=i$. But then 
\[
\left(A_{\phi}^{1},\ldots,A_{\phi}^{n}\right)\sim\left(f_{1}\left(\lambda^{1},H\right),\ldots,f_{n}\left(\lambda^{n},H\right)\right)
\]
whenever $\left(\lambda^{1},\ldots,\lambda^{n}\right)=\phi\in\Phi$,
as it is required by the general definition. 

The vector $H$ constructed in this analysis is a special case of
the \emph{reduced coupling vector} introduced next. As it turns out,
the existence of such a vector, with one random variable per each
value of each input is the general criterion for selective influences.
A criterion for a statement is another statement which is equivalent
to it. Put differently, a criterion is a condition which is both necessary
and sufficient for a given statement. 

Consider the statement that $A^{1},\ldots,A^{n}$ are selectively
influenced by $\lambda^{1},\ldots,\lambda^{n}$, respectively. By
definition, for this to be true, there should exist functions $f_{1},\ldots,f_{n}$
and a random variable $R$ such that 
\[
\left(A_{\phi}^{1},\ldots,A_{\phi}^{n}\right)\sim\left(f_{1}\left(\lambda^{1},R\right),\ldots,f_{n}\left(\lambda^{n},R\right)\right)
\]
for every $\left(\lambda^{1},\ldots,\lambda^{n}\right)=\phi\in\Phi.$
We continue to assume that every input $\lambda^{k}$ has a finite
number of values, enumerated $1,\ldots,m_{k}$. (Recall, from the
discussion of dummy inputs, that $m_{k}=1$ is allowed.)

For each $k$ and every value of $\lambda^{k}$, denote 
\[
H_{\lambda^{k}}^{k}=f_{k}\left(\lambda^{k},R\right).
\]
As $\lambda^{k}$ runs from $1$ to $m_{k}$ and $k$ runs from $1$
to $n$, this creates $m_{1}+\ldots+m_{n}$ random variables, one
random variable per each value of each input, jointly distributed
due to being functions of one and the same $R$. We have therefore
a random variable
\[
H=\left(H_{1}^{1},\ldots,H_{m_{1}}^{1},\ldots,H_{1}^{k},\ldots,H_{m_{k}}^{k},\ldots,H_{1}^{n},\ldots,H_{m_{n}}^{n}\right).
\]
If follows from the definition of selective influences that if $\left(A^{1},\ldots,A^{n}\right)\looparrowleft\left(\lambda^{1},\ldots,\lambda^{n}\right),$
then, for every allowable treatment $\phi=\left(\lambda^{1},\ldots,\lambda^{n}\right)$,
\[
\left(A_{\phi}^{1},\ldots,A_{\phi}^{n}\right)\sim\left(H_{\lambda^{1}}^{1},\ldots,H_{\lambda^{n}}^{n}\right).
\]
In other words, the existence of a jointly distributed vector of random
variables $H$ with this property is a necessary condition for $\left(A^{1},\ldots,A^{n}\right)\looparrowleft\left(\lambda^{1},\ldots,\lambda^{n}\right).$ 

Let us now assume that a vector $H$ with the above property exists.
Let us define functions as we did it in the case with stochastic independence,
\[
f_{k}\left(i,h_{1}^{1},\ldots,h_{m_{1}}^{1},\ldots,h_{1}^{k},\ldots,h_{m_{k}}^{k},\ldots,h_{1}^{n},\ldots,h_{m_{n}}^{n}\right)=h_{i}^{k}.
\]
Then
\[
\left(A_{\phi}^{1},\ldots,A_{\phi}^{n}\right)\sim\left(H_{\lambda^{1}}^{1},\ldots,H_{\lambda^{n}}^{n}\right)=\left(f_{1}\left(\lambda^{1},H\right),\ldots,f_{n}\left(\lambda^{n},H\right)\right)
\]
for every $\left(\lambda^{1},\ldots,\lambda^{n}\right)=\phi\in\Phi.$
This means that the existence of $H$ is a sufficient condition for
$\left(A^{1},\ldots,A^{n}\right)\looparrowleft\left(\lambda^{1},\ldots,\lambda^{n}\right).$

Summarizing, we have proved the following theorem.
\begin{thm}[Joint Distribution Criterion]
\label{TH: JDC}Let $\left(\lambda^{1},\ldots,\lambda^{n}\right)$
be a vector of inputs, with $\lambda^{k}\in\Lambda^{k}=\left\{ 1,\ldots,m_{k}\right\} $
($m_{k}\geq1$, $k=1,\ldots,n$). Let $\Phi\subset\Lambda^{1}\times\ldots\times\Lambda^{n}$
be a nonempty set of allowable treatments. Let $\left(A_{\phi}^{1},\ldots,A_{\phi}^{n}\right)$
be a set of random variables jointly distributed for every $\phi\in\Phi$.
Then
\[
\left(A^{1},\ldots,A^{n}\right)\looparrowleft\left(\lambda^{1},\ldots,\lambda^{n}\right)
\]
if and only if there exists a vector of jointly distributed random
variables
\[
H=\left(\overbrace{H_{1}^{1},\ldots,H_{m_{1}}^{1}},\ldots,\overbrace{H_{1}^{k},\ldots,H_{m_{k}}^{k}},\ldots,\overbrace{H_{1}^{n},\ldots,H_{m_{n}}^{n}}\right),
\]
(one variable per each value of each input) such that
\[
\left(A_{\phi}^{1},\ldots,A_{\phi}^{n}\right)\sim\left(H_{\lambda^{1}}^{1},\ldots,H_{\lambda^{n}}^{n}\right)
\]
for every $\left(\lambda^{1},\ldots,\lambda^{n}\right)=\phi\in\Phi$.
\end{thm}
The vector $H$ in this theorem is called a \emph{reduced coupling
vector} for the family $\left(\left(A_{\phi}^{1},\ldots,A_{\phi}^{n}\right):\phi\in\Phi\right)$
(or for a given pattern of selective influences). 
\begin{rem}
According to the general definition of a coupling (Section \ref{sec:Stochastically-unrelated-and}),
a coupling for the family $\left(\left(A_{\phi}^{1},\ldots,A_{\phi}^{n}\right):\phi\in\Phi\right)$
is any random variable 
\[
H^{*}=\left(\left(H_{\phi}^{1},\ldots,H_{\phi}^{n}\right):\phi\in\Phi\right)
\]
such that, for all $\phi\in\Phi$, 
\[
\left(A_{\phi}^{1},\ldots,A_{\phi}^{n}\right)\sim\left(H_{\phi}^{1},\ldots,H_{\phi}^{n}\right).
\]
The vector $H$ of Theorem \ref{TH: JDC} is obtained from such a
coupling by imposing on it additional constraints: for any $k=1,\ldots,n$
and any $\phi,\phi'\in\Phi$ sharing the same value of input $\lambda^{k}$,
\[
H_{\phi}^{k}=H_{\phi'}^{k}.
\]
These constraints allow one to reduce all different occurrences of
$H^{k}$ in $H$ to one occurrence per each value of factor $\lambda^{k}$.
Hence the adjective ``reduced'' in the name for this special coupling.
(In the literature on selective influences the reduced coupling was
also called a \emph{joint distribution vector}, and a \emph{Joint
Distribution Criterion vector}. We will not use these terms here.) 
\end{rem}
Theorem \ref{TH: JDC} is much more important than it may be suggested
by its simple proof (essentially, by means of renaming functions of
a random variable into random variables and vice versa). The reasons
for its importance are two:
\begin{enumerate}
\item it is often easier to determine whether a coupling vector exists than
whether one can find certain functions of a single random variable
(unless the latter is taken to be the reduced coupling vector and
the functions to be its projections);
\item even when a reduced coupling vector is not explicitly constructed,
its existence provides insights into the nature of the random variable
$R$ in the definition of selective influences.
\end{enumerate}
The first of these reasons is yet another illustration of the fact
that jointly distributed random variables are not, as a rule, introduced
as functions of a single random variable (see Section \ref{sec:Stochastically-unrelated-and}).
Take a simple example, when there are two binary inputs $\lambda^{1},\lambda^{2}$
(with values 1,2 each) paired with two binary outputs (with values
1,2 each). Let the set of allowable treatments consist of all four
combinations,
\[
\left(\lambda^{1}=1,\lambda^{2}=1\right),\left(\lambda^{1}=1,\lambda^{2}=2\right),\left(\lambda^{1}=2,\lambda^{2}=1\right),\left(\lambda^{1}=2,\lambda^{2}=2\right).
\]
Note that 1 and 2 as values for the inputs are chosen merely for convenience.
We could replace them by any numbers or distinct symbols (say, $\boxtimes,\boxplus$
for $\lambda^{1}$, and $ $ $\leftthreetimes,\rightthreetimes$ for
$\lambda^{2}$). The existence of the jointly distributed vectors
$\left(A_{\phi}^{1},A_{\phi}^{2}\right)$ means that for each of the
four treatments $\phi$ we are given four probabilities of the form
\begin{align*}
 & \Pr(A_{\phi}^{1}=1,A_{\phi}^{2}=1),\quad\Pr(A_{\phi}^{1}=1,A_{\phi}^{2}=2),\\
 & \Pr(A_{\phi}^{1}=2,A_{\phi}^{2}=1),\quad\Pr(A_{\phi}^{1}=2,A_{\phi}^{2}=2).
\end{align*}
Of course, the four probabilities sum to 1. Again, the use of 1 and
2 for values here is arbitrary, other symbols, generally different
for $A_{\phi}^{1}$ and $A_{\phi}^{2}$, would do as well. According
to the Joint Distribution Criterion, $\left(A^{1},A^{2}\right)\looparrowleft\left(\lambda^{1},\lambda^{2}\right)$
means the existence of four jointly distributed random variables
\[
H=\left(H_{1}^{1},H_{2}^{1},H_{1}^{2},H_{2}^{2}\right),
\]
with $H_{1}^{1}$ corresponding to the first value of input $\lambda^{1}$,
$H_{2}^{1}$ to the second value of input $\lambda^{1}$, etc., such
that
\begin{align*}
 & \left(A^{1},A^{2}\right)_{\lambda^{1}=1,\lambda^{2}=1}\sim\left(H_{1}^{1},H_{1}^{2}\right),\quad\left(A^{1},A^{2}\right)_{\lambda^{1}=1,\lambda^{2}=2}\sim\left(H_{1}^{1},H_{2}^{2}\right),\\
 & \left(A^{1},A^{2}\right)_{\lambda^{1}=2,\lambda^{2}=1}\sim\left(H_{2}^{1},H_{1}^{2}\right),\quad\left(A^{1},A^{2}\right)_{\lambda^{1}=2,\lambda^{2}=2}\sim\left(H_{2}^{1},H_{2}^{2}\right).
\end{align*}
This implies, of course, that $H_{1}^{1},H_{2}^{1},H_{1}^{2},H_{2}^{2}$
are all binary random variables, with values $1$ and 2 each. 

What is the meaning of saying that they are jointly distributed? The
meaning is that for any of the $2\times2\times2\times2$ possible
combinations of values for $H_{1}^{1},H_{2}^{1},H_{1}^{2},H_{2}^{2}$
we can find a probability,
\[
\Pr\left(H_{1}^{1}=i,H_{2}^{1}=i',H_{1}^{2}=j,H_{2}^{2}=j'\right)=p_{ii'jj'},
\]
where $i,j,i',j'\in\left\{ 1,2\right\} $. It does not matter what
these probabilities $p_{ii'jj'}$ are, insofar as they
\begin{description}
\item [{\textmd{(i)}}] are legitimate probabilities, that is, they are
nonnegative and sum to 1 across the 16 values of $H$;
\end{description}
(ii) satisfy the 2-marginal constraints
\[
\left(A^{1},A^{2}\right)_{\lambda^{1}=i,\lambda^{2}=j}\sim\left(H_{i}^{1},H_{j}^{2}\right),
\]
for all $i,j\in$$\left\{ 1,2\right\} $.

The latter translates into
\[
\small\begin{array}{c}
p_{i1j1}+p_{i1j2}+p_{i2j1}+p_{i2j2}=\Pr\left(H_{1}^{1}=i,H_{1}^{2}=j\right)=\Pr\left(A^{1}=i,A^{2}=j\right)_{\lambda^{1}=1,\lambda^{2}=1},\\
p_{i11j}+p_{i12j}+p_{i21j}+p_{i22j}=\Pr\left(H_{1}^{1}=i,H_{2}^{2}=j\right)=\Pr\left(A^{1}=i,A^{2}=j\right)_{\lambda^{1}=1,\lambda^{2}=2},\\
p_{1ij1}+p_{1ij2}+p_{2ij1}+p_{2ij2}=\Pr\left(H_{2}^{1}=i,H_{1}^{2}=j\right)=\Pr\left(A^{1}=i,A^{2}=j\right)_{\lambda^{1}=2,\lambda^{2}=1},\\
p_{1i1j}+p_{1i2j}+p_{2i1j}+p_{2i2j}=\Pr\left(H_{1}^{1}=i,H_{1}^{2}=j\right)=\Pr\left(A^{1}=i,A^{2}=j\right)_{\lambda^{1}=2,\lambda^{2}=2}.
\end{array}
\]
This is a simple system of four linear equations with 16 unknowns,
subject to being legitimate probabilities (i.e., being non-negative
and summing to 1). We will discuss this algebraic structure in the
next section, but it should be clear that this is a much more transparent
task than the one of finding a random variable $R$ and some functions,
or proving that they cannot be found.
\begin{example}
Let $A^{1},A^{2}$ have values in $\text{\{1,2\}}$ and depend on
the factors $\lambda^{1}\in\Lambda^{1}=\{1,2\}$ and $\lambda^{2}\in\Lambda^{2}=\{1,2\}$.
Let all four possible treatments be allowable. Suppose we observe
the following joint distributions of $A^{1},A^{2}$ for these treatments:\[\small%
\begin{tabular}{cc|cc|c}
$\lambda^{1}$ & $\lambda^{2}$ & $A^{1}$ & $A^{2}$ & $\Pr$\tabularnewline
\hline 
1 & 1 & 1 & 1 & $.140$\tabularnewline
 &  & 1 & 2 & $.360$\tabularnewline
 &  & 2 & 1 & $.360$\tabularnewline
 &  & 2 & 2 & $.140$\tabularnewline
\end{tabular}\qquad%
\begin{tabular}{cc|cc|c}
$\lambda^{1}$ & $\lambda^{2}$ & $A^{1}$ & $A^{2}$ & $\Pr$\tabularnewline
\hline 
1 & 2 & 1 & 1 & $.198$\tabularnewline
 &  & 1 & 2 & $.302$\tabularnewline
 &  & 2 & 1 & $.302$\tabularnewline
 &  & 2 & 2 & $.198$\tabularnewline
\end{tabular}\]\[\small%
\begin{tabular}{cc|cc|c}
$\lambda^{1}$ & $\lambda^{2}$ & $A^{1}$ & $A^{2}$ & $\Pr$\tabularnewline
\hline 
2 & 1 & 1 & 1 & $.189$\tabularnewline
 &  & 1 & 2 & $.311$\tabularnewline
 &  & 2 & 1 & $.311$\tabularnewline
 &  & 2 & 2 & $.189$\tabularnewline
\end{tabular}\qquad%
\begin{tabular}{cc|cc|c}
$\lambda^{1}$ & $\lambda^{2}$ & $A^{1}$ & $A^{2}$ & $\Pr$\tabularnewline
\hline 
2 & 2 & 1 & 1 & $.460$\tabularnewline
 &  & 1 & 2 & $.040$\tabularnewline
 &  & 2 & 1 & $.040$\tabularnewline
 &  & 2 & 2 & $.460$\tabularnewline
\end{tabular}\]The question of whether $\left(A^{1},A^{2}\right)\looparrowleft\left(\lambda^{1},\lambda^{2}\right)$
now reduces to finding a solution for the system of linear equations
mentioned above. Let us substitute the above observed probabilities
into the system:
\[
\small\begin{array}{ccc}
p_{1111}+p_{1112}+p_{1211}+p_{1212}=0.140, &  & p_{1111}+p_{1121}+p_{1211}+p_{1221}=0.198,\\
p_{1121}+p_{1122}+p_{1221}+p_{1222}=0.360, &  & p_{1112}+p_{1122}+p_{1212}+p_{1222}=0.302,\\
p_{2111}+p_{2112}+p_{2211}+p_{2212}=0.360, &  & p_{2111}+p_{2121}+p_{2211}+p_{2221}=0.302,\\
p_{2121}+p_{2122}+p_{2221}+p_{2222}=0.140, &  & p_{2112}+p_{2122}+p_{2212}+p_{2222}=0.198,\\
\\
p_{1111}+p_{1112}+p_{2111}+p_{2112}=0.189, &  & p_{1111}+p_{1121}+p_{2111}+p_{2121}=0.460,\\
p_{1121}+p_{1122}+p_{2121}+p_{2122}=0.311, &  & p_{1112}+p_{1122}+p_{2112}+p_{2122}=0.040,\\
p_{1211}+p_{1212}+p_{2211}+p_{2212}=0.311, &  & p_{1211}+p_{1221}+p_{2211}+p_{2221}=0.040,\\
p_{1221}+p_{1222}+p_{2221}+p_{2222}=0.189, &  & p_{1212}+p_{1222}+p_{2212}+p_{2222}=0.460.
\end{array}
\]
The values (found using the simplex linear programming algorithm)
\[
\small\begin{array}{lclclcl}
p_{1111}=0.067, &  & p_{1211}=0, &  & p_{2111}=0.122, &  & p_{2211}=0.04,\\
p_{1112}=0, &  & p_{1212}=0.073, &  & p_{2112}=0, &  & p_{2212}=0.198,\\
p_{1121}=0.131, &  & p_{1221}=0, &  & p_{2121}=0.14, &  & p_{2221}=0,\\
p_{1122}=0.04, &  & p_{1222}=0.189, &  & p_{2122}=0, &  & p_{2222}=0
\end{array}
\]
satisfy these equations, and as they are nonnegative and sum to one,
they represent a probability distribution. Thus, according to the
Joint Distribution Criterion, the observed joint distributions satisfy
selective influences.\hfill$\square$
\end{example}
To illustrate the second reason for the importance of Theorem \ref{TH: JDC},
we consider the following question. By the definition of selective
influences, the proposition $\left(A^{1},\ldots,A^{n}\right)\looparrowleft\left(\lambda^{1},\ldots,\lambda^{n}\right)$
means the existence of a random variable $R$ and functions $f_{1},\ldots,f_{n}$
such that
\[
\left(A_{\phi}^{1},\ldots,A_{\phi}^{n}\right)\sim\left(f_{1}\left(\lambda^{1},R\right),\ldots,f_{n}\left(\lambda^{n},R\right)\right)
\]
for every $\left(\lambda^{1},\ldots,\lambda^{n}\right)=\phi\in\Phi.$
This definition says nothing about the nature and complexity of $R$
and the functions involved, even for the simplest observable random
variables $\left(A^{1},\ldots,A^{n}\right)_{\phi}$. In most applications
$\left(A^{1},\ldots,A^{n}\right)_{\phi}$ are random variables in
the narrow sense (Section \ref{sec:Random-variables-narrow}). It
seems intuitive to expect that in such cases $R$, if it exists, is
also a random variable in the narrow sense. But this does not follow
from the definition of selective influences. Even if one manages to
prove that for a given family of random variables $\left(A^{1},\ldots,A^{n}\right)_{\phi}$
in the narrow sense this definition is satisfied by no random variable
$R$ in the narrow sense, we still do not know whether this means
that the selectiveness $\left(A^{1},\ldots,A^{n}\right)\looparrowleft\left(\lambda^{1},\ldots,\lambda^{n}\right)$
is ruled out. What if there is a random variable $R$ of a much greater
complexity (say, a random function or a random set) for which one
can find functions $f_{1},\ldots,f_{n}$ as required by the definition?

The Joint Distribution Criterion, however, allows one to rule out
such a possibility. Since the reduced coupling vector 
\[
H=\left(H_{1}^{1},\ldots,H_{m_{1}}^{1},\ldots,H_{1}^{n},\ldots,H_{m_{n}}^{n}\right),
\]
if it exists, should satisfy
\[
\left(A_{\phi}^{1},\ldots,A_{\phi}^{n}\right)\sim\left(H_{\lambda^{1}}^{1},\ldots,H_{\lambda^{n}}^{n}\right),
\]
it follows that, for any $k$ and $\lambda^{k}$,
\[
H_{\lambda^{k}}^{k}\sim A_{\phi}^{k},
\]
whenever the treatment $\phi$ contains $\lambda^{k}$. But this means
that each $H_{\lambda^{k}}^{k}$ is a random variable in a narrow
sense, and from Section \ref{sec:Random-variables-narrow} we know
then that $H$ is a random variable in the narrow sense. This constitutes
a proof of the following theorem, a simple corollary to the Joint
Distribution Criterion.
\begin{thm}
\label{TH: narrow sense}Let $\left(\lambda^{1},\ldots,\lambda^{n}\right)$,
$\Phi$, and $\left(A_{\phi}^{1},\ldots,A_{\phi}^{n}\right)$ be the
same as in Theorem \ref{TH: JDC} Let, in addition, $\left(A_{\phi}^{1},\ldots,A_{\phi}^{n}\right)$
be random variables in the narrow sense. Then 
\[
\left(A^{1},\ldots,A^{n}\right)\looparrowleft\left(\lambda^{1},\ldots,\lambda^{n}\right)
\]
if and only if there is a random variable $R$ in the narrow sense
and functions $f_{1},\ldots,f_{n}$ such that
\[
\left(A_{\phi}^{1},\ldots,A_{\phi}^{n}\right)\sim\left(f_{1}\left(\lambda^{1},R\right),\ldots,f_{n}\left(\lambda^{n},R\right)\right)
\]
for every $\left(\lambda^{1},\ldots,\lambda^{n}\right)=\phi\in\Phi.$ 
\end{thm}
If one feels dissatisfied with considering vectors of random variables
on a par with ``single'' random variables, this dissatisfaction
is not well-grounded. The fact is, the dimensionality of vectors of
random variables in the narrow sense is not essential. Consider, for
example, the reduced coupling vector
\[
H=\left(H_{1}^{1},H_{2}^{1},H_{1}^{2},H_{2}^{2}\right),
\]
constructed earlier for two binary random variables selectively influenced
by two binary inputs. Clearly, in all considerations this four-component
vector of binary random variables can be replaced with a single 16-valued
random variable, $H'$. Let these 16 values be $0,\ldots,15.$ The
two variables are equivalent if one puts
\begin{align*}
 & \Pr\left(H_{1}^{1}=i,H_{2}^{1}=i',H_{1}^{2}=j,H_{2}^{2}=j'\right)\\
 & =\Pr\left(H'=\left(i-1\right)2^{3}+\left(i'-1\right)2^{2}+\left(j-1\right)2+\left(j'-1\right)\right).
\end{align*}
In particular, any functions of $H$ can be presented as functions
of $H'$. 

In the case of continuous random variables the situation is, in a
sense, even simpler, although we will have to omit the underlying
justification. It follows from the theory of Borel-equivalent spaces
(which is part of descriptive set theory), that any vector of continuous
random variables
\[
R=\left(R^{1},\ldots,R^{k}\right),
\]
can be presented as a function of any continuous variable $R'$ with
an \emph{atomless} distribution on an interval of real numbers. The
``atomlessness'' means that the sigma-algebra of $R'$ contains
no null-set whose probability measure is not zero. Simple examples
are uniformly and normally distributed random variables. If the vector
is discrete, the previous statement applies with no modifications
(although we know that in this case one can also choose a discrete
$R'$). It follows that the statement also applies to mixed vectors,
containing both discrete and continuous random variables (or vectors
thereof, or vectors of vectors thereof, etc.)

We can complement, therefore, Theorem \ref{TH: narrow sense} with
the following statement.
\begin{thm}
\label{TH: Under-the-conditions}Under the conditions of Theorem \ref{TH: narrow sense},
the random variable $R$ can always be chosen to be any continuous
random variable with an atomless distribution on an interval of real
numbers. If all the random variables $A_{\phi}^{1},\ldots,A_{\phi}^{n}$
are discrete (in particular, have finite numbers of values), then
$R$ can be chosen to be discrete (respectively, have finite number
of values). 
\end{thm}
We have quite a bit more specificity now than based on the initial
definition of selective influences. And it is achieved due to the
Joint Distribution Criterion almost ``automatically.''

Theorem \ref{TH: JDC} is not restricted to finite-valued inputs.
Nor is it restricted to a finite number of inputs, or to outputs of
a specific kind. It is completely general. For the reader's convenience,
we formulate here the general version of this theorem, avoiding all
elaborations.
\begin{thm}[Joint Distribution Criterion (general version)]
Let $\left(\lambda^{k}:k\in K\right)$ be an indexed family of inputs,
with $\lambda^{k}\in\Lambda^{k}\not=\emptyset$, for all $k\in K$.
Let $\Phi\subset\prod_{k\in K}\Lambda^{k}$ be a nonempty set of allowable
treatments. Let $\left(A_{\phi}^{k}:k\in K\right)$ be a family of
random variables jointly distributed for every $\phi\in\Phi$. Then
\[
\left(A^{k}:k\in K\right)\looparrowleft\left(\lambda^{k}:k\in K\right)
\]
if and only if there exists an indexed family of jointly distributed
random variables
\[
H=\left(H_{\lambda^{k}}^{k}:\lambda^{k}\in\Lambda^{k},k\in K\right),
\]
(one variable per each value of each input) such that
\[
\left(A_{\phi}^{k}:k\in K\right)\sim\left(H_{\lambda^{k}}^{k}:k\in K\right)
\]
for every $\left(\lambda^{k}:k\in K\right)=\phi\in\Phi$.
\end{thm}

\section{\label{sec:Properties-of-selective}Properties of selective influences
and tests}

Certain properties of selective influences (in the canonical form)
are immediately obvious. 

The first one is \emph{nestedness with respect to input values}: if
random outputs $A^{1},\ldots,A^{n}$ are selectively influenced by
inputs $\lambda^{1},\ldots,\lambda^{n}$, with sets of possible values
$\Lambda^{1},\ldots,\Lambda^{n}$, then the same random outputs are
selectively influenced by inputs $\lambda'^{1},\ldots,\lambda'^{n}$
whose sets of possible values are $\Lambda'^{1}\subset\Lambda^{1},\ldots,\Lambda'^{n}\subset\Lambda^{n}$.
Every variable is essentially the set of its possible values. Inputs
are no exception. In fact, in a more rigorous development $\lambda$
would be reserved for input values, whereas input themselves, considered
as variables, would be identified by $\Lambda$ (see Remark \ref{REM: As-it-is}).
When a set of an input's values changes, the input is being replaced
by a new one. The nestedness property in question tells us that if
the change consists in removing some of the possible values of some
of the inputs, the selectiveness pattern established for the original
inputs cannot be violated. This does not, of course, work in the other
direction: if we augment $\Lambda^{1},\ldots,\Lambda^{n}$ by adding
to them new elements, then the initial pattern of selectiveness may
very well disappear.

The second property is \emph{nestedness with respect to inputs and
outputs} (in a canonical diagram they are in a one-to-one correspondence):
if a vector of random outputs is selectively influenced by a vector
of inputs, then any subvector of the random outputs is selectively
influenced by the corresponding subvector of the inputs. In symbols,
if 
\[
\left(A^{1},\ldots,A^{n}\right)\looparrowleft\left(\lambda^{1},\ldots,\lambda^{n}\right)
\]
and $i_{1},\ldots,i_{k}\in\left\{ 1,\ldots,n\right\} $, then
\[
\left(A^{i_{1}},\ldots,A^{i_{k}}\right)\looparrowleft\left(\lambda^{i_{1}},\ldots,\lambda^{i_{k}}\right).
\]

Note that the set of allowable treatments has to be redefined whether
we eliminate certain input-output pairs or certain input values. In
the latter case, the new set of allowable treatments is the largest
$\Phi'\subset\Lambda'^{1}\times\ldots\times\Lambda'^{n}$, such that
$\Phi'\subset\Phi$. In the case we drop input-output pairs, the new
set of allowable treatments is the largest $\Phi''\subset\Lambda^{i_{1}}\times\ldots\times\Lambda^{i_{k}}$,
such that every $\phi''\in\Phi''$ is a part of some $\phi\in\Phi$. 

Both these nestedness properties follow from the fact that any subset
of random variables that are components of a reduced coupling vector
\[
H=\left(H_{1}^{1},\ldots,H_{m_{1}}^{1},\ldots,H_{1}^{n},\ldots,H_{m_{n}}^{n}\right),
\]
are also jointly distributed. When we eliminate an $i$th value of
input $k$, we drop from this vector $H_{i}^{k}$. When we eliminate
an input $k$, we drop the subvector $H_{1}^{k},\ldots,H_{m_{k}}^{k}$.
In both cases the resulting $H'$ is easily checked to be a reduced
coupling vector for the redefined sets of treatments and outputs. 

By similar arguments one can establish that a pattern of selective
influences is well-behaved in response to all possible groupings of
the inputs, with or without a corresponding grouping of outputs: thus,
if 
\[
\left(A^{1},\ldots,A^{k},\ldots,A^{l},\ldots,A^{n}\right)\looparrowleft\left(\lambda^{1},\ldots,\lambda^{k},\ldots,\lambda^{l},\ldots,\lambda^{n}\right),
\]
then
\[
\left(A^{1},\ldots,A^{k},\ldots,A^{l},\ldots,A^{n}\right)\looparrowleft\left(\lambda^{1},\ldots,\left(\lambda^{k},\lambda^{l}\right),\ldots,\left(\lambda^{k},\lambda^{l}\right),\ldots,\lambda^{n}\right)
\]
and
\begin{align*}
 & \left(A^{1},\ldots,\left(A^{k},A^{l}\right),\ldots,\left(A^{k},A^{l}\right),\ldots,A^{n}\right)\\
 & \looparrowleft\left(\lambda^{1},\ldots,\left(\lambda^{k},\lambda^{l}\right),\ldots,\left(\lambda^{k},\lambda^{l}\right),\ldots,\lambda^{n}\right).
\end{align*}
We omit the details related to redefinitions of allowable treatments.

A simple consequence of the nestedness with respect to input-output
pairs turns out to be of a great importance for determining if a selectiveness
pattern is present. This consequence is called \emph{complete marginal
selectivity}: if $\left(A^{1},\ldots,A^{n}\right)\looparrowleft\left(\lambda^{1},\ldots,\lambda^{n}\right)$
and $i_{1},\ldots,i_{k}\in\left\{ 1,\ldots,n\right\} $, then the
distribution of $\left(A_{\phi}^{i_{1}},\ldots,A_{\phi}^{i_{k}}\right)$
depends only on $\left(\lambda^{i_{1}},\ldots,\lambda^{i_{k}}\right)$.
In other words, if $\phi$ and $\phi'$ include the same subset $\left(\lambda^{i_{1}},\ldots,\lambda^{i_{k}}\right)$,
\[
\left(A_{\phi}^{i_{1}},\ldots,A_{\phi}^{i_{k}}\right)\sim\left(A_{\phi'}^{i_{1}},\ldots,A_{\phi'}^{i_{k}}\right).
\]
In particular (\emph{simple marginal selectivity}),
\[
A_{\phi}^{i}\sim A_{\phi'}^{i}
\]
for any $\phi$ and $\phi'$ that share a value of $\lambda^{i}$
($i=1,\ldots,n$). The importance of marginal selectivity is that
it is easy to check, ruling out selective influences whenever it is
found violated.
\begin{example}
Let $A^{1},A^{2}$ have values in $\{1,2\}$ and depend on the external
factors $\lambda^{1}\in\Lambda^{1}=\{1,2\}$ and $\lambda^{2}\in\Lambda^{2}=\{1,2\}$.
Let the joint distribution of $A^{1},A^{2}$ for each treatment (all
four being allowable) be as follows:\[\small%
\begin{tabular}{r|cc|ccr|cc|c}
$\lambda^{1}=1$, $\lambda^{2}=1$ & $\!\! A^{2}=1\!\!$ & $\!\! A^{2}=2\!\!$ &  &  & $\lambda^{1}=1$, $\lambda^{2}=2$ & $\!\! A^{2}=1\!\!$ & $\!\! A^{2}=2\!\!$ & \tabularnewline
\cline{1-4} \cline{6-9} 
$A^{1}=1$ & .2 & .2 & .4 &  & $A^{1}=1$ & .3 & .1 & .4\tabularnewline
$A^{1}=2$ & .3 & .3 & .6 &  & $A^{1}=2$ & .2 & .4 & .6\tabularnewline
\cline{1-4} \cline{6-9} 
 & .5 & .5 &  &  &  & .5 & .5 & \tabularnewline
\multicolumn{1}{r}{} &  & \multicolumn{1}{c}{} &  &  & \multicolumn{1}{r}{} &  & \multicolumn{1}{c}{} & \tabularnewline
$\lambda^{1}=2$, $\lambda^{2}=1$ & $\!\! A^{2}=1\!\!$ & $\!\! A^{2}=2\!\!$ &  &  & $\lambda^{1}=2$, $\lambda^{2}=2$ & $\!\! A^{2}=1\!\!$ & $\!\! A^{2}=2\!\!$ & \tabularnewline
\cline{1-4} \cline{6-9} 
$A^{1}=1$ & .4 & .3 & .7 &  & $A^{1}=1$ & .3 & .4 & .7\tabularnewline
$A^{1}=2$ & .1 & .2 & .3 &  & $A^{1}=2$ & .1 & .2 & .3\tabularnewline
\cline{1-4} \cline{6-9} 
 & .5 & .5 &  &  &  & .4 & .6 & \tabularnewline
\end{tabular}\]Marginal selectivity here is violated because the marginal distribution
of $A^{2}$ changes when $\lambda^{2}=2$ and $\lambda^{1}$ changes
from $1$ to $2$.\hfill$\square$
\end{example}
Marginal selectivity is strictly weaker than selective influences.
The latter do imply marginal selectivity, but marginal selectivity
can very well hold in the absence of selective influences.
\begin{example}
\label{Example:marginal-selectivity-but}Consider the following joint
distributions:\[\small%
\begin{tabular}{r|cc|ccr|cc|c}
$\lambda^{1}=1$, $\lambda^{2}=1$ & $\!\! A^{2}=1\!\!$ & $\!\! A^{2}=2\!\!$ &  &  & $\lambda^{1}=1$, $\lambda^{2}=2$ & $\!\! A^{2}=1\!\!$ & $\!\! A^{2}=2\!\!$ & \tabularnewline
\cline{1-4} \cline{6-9} 
$A^{1}=1$ & .5 & 0 & .5 &  & $A^{1}=1$ & .5 & 0 & .5\tabularnewline
$A^{1}=2$ & 0 & .5 & .5 &  & $A^{1}=2$ & 0 & .5 & .5\tabularnewline
\cline{1-4} \cline{6-9} 
 & .5 & .5 &  &  &  & .5 & .5 & \tabularnewline
\multicolumn{1}{r}{} &  & \multicolumn{1}{c}{} &  &  & \multicolumn{1}{r}{} &  & \multicolumn{1}{c}{} & \tabularnewline
$\lambda^{1}=2$, $\lambda^{2}=1$ & $\!\! A^{2}=1\!\!$ & $\!\! A^{2}=2\!\!$ &  &  & $\lambda^{1}=2$, $\lambda^{2}=2$ & $\!\! A^{2}=1\!\!$ & $\!\! A^{2}=2\!\!$ & \tabularnewline
\cline{1-4} \cline{6-9} 
$A^{1}=1$ & .5 & 0 & .5 &  & $A^{1}=1$ & 0 & .5 & .5\tabularnewline
$A^{1}=2$ & 0 & .5 & .5 &  & $A^{1}=2$ & .5 & 0 & .5\tabularnewline
\cline{1-4} \cline{6-9} 
 & .5 & .5 &  &  &  & .5 & .5 & \tabularnewline
\end{tabular}\]Marginal selectivity is trivially satisfied as all marginals are
uniform. However, $\left(A^{1},A^{2}\right)\not\looparrowleft\left(\lambda^{1},\lambda^{2}\right)$
in this case. The joint distribution criterion would require the existence
of a jointly distributed vector $H$ whose components satisfy $\left(A_{ij}^{1},A_{ij}^{2}\right)\sim\left(H_{i}^{1},H_{j}^{2}\right)$
for $i,j\in\{1,2\}$. But combining this with the above joint distributions,
we obtain 
\[
H_{1}^{1}=H_{1}^{2},\qquad H_{1}^{1}=H_{2}^{2},\qquad H_{2}^{1}=H_{1}^{2},\qquad H_{2}^{1}=3-H_{2}^{2},
\]
which yields the contradiction 
\[
3-H_{2}^{2}=H_{2}^{2}.
\]
\hfill$\square$
\end{example}
Another property of selective influences is that if $\left(A^{1},\ldots,A^{n}\right)\looparrowleft\left(\lambda^{1},\ldots,\lambda^{n}\right)$,
and if, for all $\phi=\left(\lambda^{1},\ldots,\lambda^{n}\right)\in\Phi$,
\[
B_{\phi}^{1}=g_{1}\left(\lambda^{1},A_{\phi}^{1}\right),\ldots,B_{\phi}^{n}=g_{n}\left(\lambda^{n},A_{\phi}^{n}\right),
\]
then $\left(B^{1},\ldots,B^{n}\right)\looparrowleft\left(\lambda^{1},\ldots,\lambda^{n}\right)$.
The functions $g_{1},\ldots,g_{n}$ are referred to as \emph{input-value-specific
transformations of random outputs}. The property in question therefore
is the\emph{ }invariance of selective influences, if established,
with respect to such transformations. 

Let us make sure that this property is true. According to the general
definition, we have a random variable $R$ and functions $f_{1},\ldots,f_{n}$
such that
\[
\left(A_{\phi}^{1},\ldots,A_{\phi}^{n}\right)\sim\left(f_{1}\left(\lambda^{1},R\right),\ldots,f_{n}\left(\lambda^{n},R\right)\right),
\]
for every $\phi=\left(\lambda^{1},\ldots,\lambda^{n}\right)\in\Phi$.
But then
\[
\left(B_{\phi}^{1},\ldots,B_{\phi}^{n}\right)\sim\left(g_{1}\left(\lambda^{1},f_{1}\left(\lambda^{1},R\right)\right),\ldots,g_{n}\left(\lambda^{n},f_{n}\left(\lambda^{n},R\right)\right)\right),
\]
and every $g_{k}\left(\lambda^{k},f_{k}\left(\lambda^{k},R\right)\right)$
is some function $f{}_{k}^{*}\left(\lambda^{k},R\right)$. The vectors
$\left(B_{\phi}^{1},\ldots,B_{\phi}^{n}\right)$ therefore satisfy
the definition too.

As a special case, the transformation may not depend on input values,
\[
B_{\phi}^{1}=g_{1}\left(A_{\phi}^{1}\right),\ldots,B_{\phi}^{n}=g_{n}\left(A_{\phi}^{n}\right).
\]
This would include all possible \emph{renamings} and \emph{groupings}
of the values of the random outputs: a pattern of selective influences
is preserved under all such transformations. For instance, one can
rename values $1,2$ of a binary output into $\sqcup,\sqcap$, or
one can group values $1,2,3,4$ into ``cruder'' values, by means
of a transformation like
\[
1\mapsto\sqcup,2\mapsto\sqcup,3\mapsto\sqcap,4\mapsto\sqcap.
\]

The meaning of the input-value-specificity is this. We choose a $k\in\left\{ 1,\ldots,n\right\} $
and assume, for simplicity, that $\lambda^{k}$ has discrete values,
$1,2,\ldots$. Let $A_{\phi}^{k}$ be transformed into random variables
$B_{1,\phi}^{k}$, $B_{2,\phi}^{k}$, etc., all sharing the same set
of possible values and the same sigma-algebra. We know that one can
replace $A^{k}$ in
\[
\left(A^{1},\ldots A^{k},\ldots,A^{n}\right)\looparrowleft\left(\lambda^{1},\ldots,\lambda^{k},\ldots,\lambda^{n}\right)
\]
with any of these new random variables,
\begin{align*}
\left(A^{1},\ldots B_{1}^{k},\ldots,A^{n}\right) & \looparrowleft\left(\lambda^{1},\ldots,\lambda^{k},\ldots,\lambda^{n}\right),\\
\left(A^{1},\ldots B_{2}^{k},\ldots,A^{n}\right) & \looparrowleft\left(\lambda^{1},\ldots,\lambda^{k},\ldots,\lambda^{n}\right),\\
 & \textnormal{etc.}
\end{align*}
The input-value-specificity is involved if one forms a random variable
\[
B_{\phi}^{k}=\left\{ \begin{array}{ccc}
B_{1,\phi}^{k} & \textnormal{if} & \lambda^{k}=1\\
B_{2,\phi}^{k} & \textnormal{if} & \lambda^{k}=2\\
 & \textnormal{etc.}
\end{array}\right.
\]
The invariance property says that this random variable, too, can replace
$A^{k}$ in a pattern of selective influences,
\[
\left(A^{1},\ldots B^{k},\ldots,A^{n}\right)\looparrowleft\left(\lambda^{1},\ldots,\lambda^{k},\ldots,\lambda^{n}\right).
\]

Note that the property in question works in one direction only: if
$\left(A^{1},\ldots,A^{n}\right)\looparrowleft\left(\lambda^{1},\ldots,\lambda^{n}\right)$
then $\left(B^{1},\ldots,B^{n}\right)\looparrowleft\left(\lambda^{1},\ldots,\lambda^{n}\right)$.
It is perfectly possible (if we use grouping of values) that $\left(A^{1},\ldots,A^{n}\right)\not\looparrowleft\left(\lambda^{1},\ldots,\lambda^{n}\right)$
but following an input-value-specific transformation, $\left(B^{1},\ldots,B^{n}\right)\looparrowleft\left(\lambda^{1},\ldots,\lambda^{n}\right)$.
However, if the transformation $B_{\phi}^{1}=g_{1}(\lambda^{1},A_{\phi}^{1}),\ldots,B_{\phi}^{n}=g_{n}(\lambda^{n},A_{\phi}^{n})$,
is reversible, that is, there exist another transformation $A_{\phi}^{1}=h_{1}(\lambda^{1},B_{\phi}^{1}),\ldots,A_{\phi}^{n}=h_{n}(\lambda^{n},B_{\phi}^{n})$
back to the original variables, then $\left(A^{1},\ldots,A^{n}\right)\looparrowleft\left(\lambda^{1},\ldots,\lambda^{n}\right)$
if and only if $\left(B^{1},\ldots,B^{n}\right)\looparrowleft\left(\lambda^{1},\ldots,\lambda^{n}\right)$.
\begin{example}
\label{Example:input-value-specific}Consider the random variables
$A^{1},A^{2}$ with values in $\{1,2\}$, depending on the input factors
$\lambda^{1}\in\{1,2\}$, $\lambda^{2}\in\{1,2\}$, and having the
following joint distributions at the four possible treatments:\[\small%
\begin{tabular}{r|cc|ccr|cc|c}
$\lambda^{1}=1$, $\lambda^{2}=1$ & $\!\! A^{2}=1\!\!$ & $\!\! A^{2}=2\!\!$ &  &  & $\lambda^{1}=1$, $\lambda^{2}=2$ & $\!\! A^{2}=1\!\!$ & $\!\! A^{2}=2\!\!$ & \tabularnewline
\cline{1-4} \cline{6-9} 
$A^{1}=1$ & 0.3 & 0.4 & 0.7 &  & $A^{1}=1$ & 0.35 & 0.35 & 0.7\tabularnewline
$A^{1}=2$ & 0.1 & 0.2 & 0.3 &  & $A^{1}=2$ & 0.15 & 0.15 & 0.3\tabularnewline
\cline{1-4} \cline{6-9} 
 & 0.4 & 0.6 &  &  &  & 0.5 & 0.5 & \tabularnewline
\multicolumn{1}{r}{} &  & \multicolumn{1}{c}{} &  &  & \multicolumn{1}{r}{} &  & \multicolumn{1}{c}{} & \tabularnewline
$\lambda^{1}=2$, $\lambda^{2}=1$ & $\!\! A^{2}=1\!\!$ & $\!\! A^{2}=2\!\!$ &  &  & $\lambda^{1}=2$, $\lambda^{2}=2$ & $\!\! A^{2}=1\!\!$ & $\!\! A^{2}=2\!\!$ & \tabularnewline
\cline{1-4} \cline{6-9} 
$A^{1}=1$ & 0.32 & 0.48 & 0.8 &  & $A^{1}=1$ & 0.45 & 0.35 & 0.8\tabularnewline
$A^{1}=2$ & 0.08 & 0.12 & 0.2 &  & $A^{1}=2$ & 0.05 & 0.15 & 0.2\tabularnewline
\cline{1-4} \cline{6-9} 
 & 0.4 & 0.6 &  &  &  & 0.5 & 0.5 & \tabularnewline
\end{tabular}\]We will see in the next section that $\left(A^{1},A^{2}\right)\looparrowleft\left(\lambda^{1},\lambda^{2}\right)$
is satisfied in this case. Let us define the input value specific
transformations $B^{1}=g_{1}(\lambda^{1},A^{1})$ and $B^{2}=g_{2}(\lambda^{2},A^{2})$,
where
\begin{align*}
 & g_{1}(1,\{1,2\})=\{+1,-1\},\quad g_{2}(1,\{1,2\})=\{7,3\},\\
 & g_{1}(2,\{1,2\})=\{-1,+1\},\quad g_{2}(2,\{1,2\})=\{3,7\}.
\end{align*}
As we see, $A^{1}=1$ is mapped into $B^{1}=+1$ or $B^{1}=-1$ according
as $\lambda^{1}$ is 1 or 2, $A^{2}=1$ is mapped into $B^{2}=7$
or $B^{2}=3$ according as $\lambda^{2}$ is 1 or 2, etc. We obtain
the following joint distributions\[\small%
\begin{tabular}{r|cc|ccr|cc|c}
$\lambda^{1}=1$, $\lambda^{2}=1$ & $\!\! B^{2}=7\!\!$ & $\!\! B^{2}=3\!\!$ &  &  & $\lambda^{1}=1$, $\lambda^{2}=2$ & $\!\! B^{2}=7\!\!$ & $\!\! B^{2}=3\!\!$ & \tabularnewline
\cline{1-4} \cline{6-9} 
$B^{1}=\text{+1}$ & 0.3 & 0.4 & 0.7 &  & $B^{1}=\text{+1}$ & 0.35 & 0.35 & 0.7\tabularnewline
$B^{1}=-1$ & 0.1 & 0.2 & 0.3 &  & $B^{1}=-1$ & 0.15 & 0.15 & 0.3\tabularnewline
\cline{1-4} \cline{6-9} 
 & 0.4 & 0.6 &  &  &  & 0.5 & 0.5 & \tabularnewline
\multicolumn{1}{r}{} &  & \multicolumn{1}{c}{} &  &  & \multicolumn{1}{r}{} &  & \multicolumn{1}{c}{} & \tabularnewline
$\lambda^{1}=2$, $\lambda^{2}=1$ & $\!\! B^{2}=7\!\!$ & $\!\! B^{2}=3\!\!$ &  &  & $\lambda^{1}=2$, $\lambda^{2}=2$ & $\!\! B^{2}=7\!\!$ & $\!\! B^{2}=3\!\!$ & \tabularnewline
\cline{1-4} \cline{6-9} 
$B^{1}=\text{+1}$ & 0.08 & 0.12 & 0.2 &  & $B^{1}=\text{+1}$ & 0.15 & 0.05 & 0.2\tabularnewline
$B^{1}=-1$ & 0.32 & 0.48 & 0.8 &  & $B^{1}=-1$ & 0.35 & 0.45 & 0.8\tabularnewline
\cline{1-4} \cline{6-9} 
 & 0.4 & 0.6 &  &  &  & 0.5 & 0.5 & \tabularnewline
\end{tabular}\]We know that the transformed variables satisfy $\left(B^{1},B^{2}\right)\looparrowleft\left(\lambda^{1},\lambda^{2}\right)$
because $\left(A^{1},A^{2}\right)\looparrowleft\left(\lambda^{1},\lambda^{2}\right)$.\hfill$\square$
\end{example}
In the subsequent sections we will consider several \emph{tests of
selective influences}. Such a test is always a statement whose truth
value (whether it is true or false) determines whether a given pattern
of selective influences holds or does not hold. The truth value of
the test statement must be determinable from the distributions of
$\left(A_{\phi}^{1},\ldots,A_{\phi}^{n}\right)$ for all allowable
$\phi$. If its truth implies $\left(A^{1},\ldots,A^{n}\right)\looparrowleft\left(\lambda^{1},\ldots,\lambda^{n}\right)$,
then the test provides a sufficient condition for selective influences;
if its falsity implies $\left(A^{1},\ldots,A^{n}\right)\not\looparrowleft\left(\lambda^{1},\ldots,\lambda^{n}\right)$,
then the test provides a necessary condition for selective influences.
If the test provides both necessary and sufficient condition, it is
a criterion. 

The distribution of $\left(A_{\phi}^{1},\ldots,A_{\phi}^{n}\right)$,
if the random variables are known from their observed realizations,
cannot be known precisely, because probabilities are never observable.
All our tests require that the distributions of $\left(A_{\phi}^{1},\ldots,A_{\phi}^{n}\right)$,
or at least some parameters thereof, be known precisely. Therefore
they can only be applied to empirical observations if the latter are
replaced by theoretical distributions. This can be done based on statistical
considerations, outside the scope of the tests themselves. In particular,
if all sample sizes are sufficiently large, theoretical distributions
can be assumed to be so close to the empirical ones that their difference
cannot affect the outcome of a test.

As follows from the discussion above, the most basic and obvious test
of selective influences is the (complete) \emph{marginal selectivity
test}. This is a necessary condition for selective influences: if,
at least for one pair of distinct treatments $\phi$ and $\phi'$
that include one and the same subvector $\left(\lambda^{i_{1}},\ldots,\lambda^{i_{k}}\right)$,
the distributions of the $k$-marginal random variables $\left(A_{\phi}^{i_{1}},\ldots,A_{\phi}^{i_{k}}\right)$
and $\left(A_{\phi'}^{i_{1}},\ldots,A_{\phi'}^{i_{k}}\right)$ are
not the same, then $\left(A^{1},\ldots,A^{n}\right)\not\looparrowleft\left(\lambda^{1},\ldots,\lambda^{n}\right).$

\section{Linear Feasibility Test}

In this section we will discuss a test which is both a necessary and
sufficient condition for the selective influences in the case when
the number of input-output pairs, the set of values of each input,
and the set of possible values of each random output are all finite.
Let us enumerate, for $k=1,\ldots,n$, the values of each input $\lambda^{k}$
as $1,\ldots,m_{k}$, and the values of each random output $A^{k}$
as $1,\ldots,v_{k}$. In Section \ref{sec:Joint-Distribution-Criterion}
we discussed the case $n=2,$ $m_{1}=m_{2}=2$, and $v_{1}=v_{2}=2$.
We determined there that the question of whether $\left(A^{1},A^{2}\right)\looparrowleft\left(\lambda^{1},\lambda^{2}\right)$
translates into a question of whether certain linear equations have
a solution subject to certain constraints. We will see that this is
the case generally.

The observable distributions of $\left(A_{\phi}^{1},\ldots,A_{\phi}^{n}\right)$
are represented by the probabilities of the events that can be described
as
\[
\left(\overbrace{A^{1}=a_{1},\ldots,A^{k}=a_{k},\ldots,A^{n}=a_{n}};\overbrace{\lambda^{1}=l_{1},\ldots,\lambda^{k}=l_{k},\ldots,\lambda^{n}=l_{n}}\right),
\]
where $a_{k}\in\left\{ 1,\ldots,v_{k}\right\} $ (output values) and
$l_{k}\in\left\{ 1,\ldots,m_{k}\right\} $ (input values). Let us
form a matrix $M$ whose rows are enumerated (labeled) by all such
vectors. We only consider the vectors with allowable treatments, 
\[
\phi=\left(\lambda^{1}=l_{1},\ldots,\lambda^{k}=l_{k},\ldots,\lambda^{n}=l_{n}\right)\in\Phi.
\]
If the number of the allowable treatments is $t$ (between 1 and $m_{1}\times\ldots\times m_{n}$),
then the number of the rows in $M$ is $t\times v_{1}\times\ldots\times v_{n}$. 

The columns of the matrix $M$ are enumerated (labeled) by the vectors
of the form 
\[
\left(\overbrace{H_{1}^{1}=h_{1}^{1},\ldots H_{m_{1}}^{1}=h_{m_{1}}^{1}},\ldots,\overbrace{H_{1}^{n}=h_{1}^{n},\dots,H_{m_{n}}^{n}=h_{m_{n}}^{n}}\right),
\]
where $h_{i}^{k}\in\left\{ 1,\ldots,v_{k}\right\} $. Such vectors
represent events whose probabilities define the distribution of a
reduced coupling vector $H$ (if one exists). The number of such events,
hence the number of the columns in $M$ is $\left(v_{1}\right)^{m_{1}}\times\ldots\times\left(v_{n}\right)^{m_{n}}$
(where the superscripts represent conventional exponents). 

We also form a column vector $P$ whose elements are labeled in the
same way and in the same order as the rows of the matrix $M$, and
a column vector $Q$ whose elements are labeled in the same way and
in the same order as the columns of the matrix $M$. 

Let us now fill in the entries of the vectors $P,Q$, and the matrix
$M$. The matrix $M$ is Boolean: it is filled with 1's and 0's. Consider
a cell $\left(I,J\right)$ belonging to the column labeled
\[
J=\left(\overbrace{H_{1}^{1}=h_{1}^{1},\ldots H_{m_{1}}^{1}=h_{m_{1}}^{1}},\ldots,\overbrace{H_{1}^{n}=h_{1}^{n},\dots,H_{m_{n}}^{n}=h_{m_{n}}^{n}}\right)
\]
and to the row labeled
\[
I=\left(\overbrace{A^{1}=a_{1},\ldots,A^{k}=a_{k},\ldots,A^{n}=a_{n}};\overbrace{\lambda^{1}=l_{1},\ldots,\lambda^{k}=l_{k},\ldots,\lambda^{n}=l_{n}}\right).
\]
In the vector-label $J$ pick the entries
\[
H_{l_{1}}^{1}=h_{l_{1}}^{1},\ldots,H_{l_{k}}^{k}=h_{l_{k}}^{k},\ldots,H_{l_{n}}^{n}=h_{l_{n}}^{n}
\]
corresponding to the values of $\left(\lambda^{1},\ldots,\lambda^{n}\right)$
indicated in the vector-label $I$. If 
\[
\left(h_{l_{1}}^{1},\ldots,h_{l_{k}}^{k},\ldots,h_{l_{n}}^{n}\right)=\left(a_{1},\ldots,a_{k},\ldots,a_{n}\right)
\]
then the cell $\left(I,J\right)$ should be filled with 1; otherwise
its value is 0.

The vector $P$ is filled with the probabilities
\[
\Pr\left(A^{1}=a_{1},\ldots,A^{n}=a_{n}\right)_{\phi=\left(\lambda^{1}=l_{1},\ldots,\lambda^{n}=l_{n}\right)}.
\]
For any allowable $\phi$, the probabilities across all possible combinations
of $\left(a_{1},\ldots,a_{n}\right)$ sum to 1. These probabilities
are assumed to be known.

The vector $Q$ is filled with the probabilities
\[
\Pr\left(H_{1}^{1}=h_{1}^{1},\ldots,H_{m_{1}}^{1}=h_{m_{1}}^{1},\ldots,H_{1}^{n}=h_{1}^{n},\dots,H_{m_{n}}^{n}=h_{m_{n}}^{n}\right),
\]
which sum to 1 across all possible values of $\left(h_{1}^{1},\ldots,h_{m_{1}}^{1},\ldots,h_{1}^{n},\dots,h_{m_{n}}^{n}\right)$.
These probabilities are not known, they have to be found or determined
not to exist. 

\def\0{$\!\cdot\!\!$}\def\1{$\!1\!\!$}\def\2{$\!2\!\!$}
\begin{example}
\label{Example:PMQ-filled}Let us now apply these general definitions
to the simplest nontrivial case $n=2$, $m_{1}=m_{2}=2$, $v_{1}=v_{2}=2$
considered in Section \ref{sec:Joint-Distribution-Criterion}. The
matrix $M$ filled with binary values is (replacing 0 with ``$\cdot$''
for better legibility) \[\small%
\begin{tabular}{cr|cccccccccccccccc}
 & $H_{1}^{1}$ & \1 & \1 & \1 & \1 & \1 & \1 & \1 & \1 & \2 & \2 & \2 & \2 & \2 & \2 & \2 & \2\tabularnewline
 & $H_{2}^{1}$ & \1 & \1 & \1 & \1 & \2 & \2 & \2 & \2 & \1 & \1 & \1 & \1 & \2 & \2 & \2 & \2\tabularnewline
 & $H_{1}^{2}$ & \1 & \1 & \2 & \2 & \1 & \1 & \2 & \2 & \1 & \1 & \2 & \2 & \1 & \1 & \2 & \2\tabularnewline
 & $H_{2}^{2}$ & \1 & \2 & \1 & \2 & \1 & \2 & \1 & \2 & \1 & \2 & \1 & \2 & \1 & \2 & \1 & \2\tabularnewline
\hline 
\multirow{4}{*}{$\lambda^{1}=1$, $\lambda^{2}=1$} & $A^{1}=1$, $A^{2}=1$ & \1 & \1 & \0 & \0 & \1 & \1 & \0 & \0 & \0 & \0 & \0 & \0 & \0 & \0 & \0 & \0\tabularnewline
 & $A^{1}=1$, $A^{2}=2$ & \0 & \0 & \1 & \1 & \0 & \0 & \1 & \1 & \0 & \0 & \0 & \0 & \0 & \0 & \0 & \0\tabularnewline
 & $A^{1}=2$, $A^{2}=1$ & \0 & \0 & \0 & \0 & \0 & \0 & \0 & \0 & \1 & \1 & \0 & \0 & \1 & \1 & \0 & \0\tabularnewline
 & $A^{1}=2$, $A^{2}=2$ & \0 & \0 & \0 & \0 & \0 & \0 & \0 & \0 & \0 & \0 & \1 & \1 & \0 & \0 & \1 & \1\tabularnewline
\hline 
\multirow{4}{*}{$\lambda^{1}=1$, $\lambda^{2}=2$} & $A^{1}=1$, $A^{2}=1$ & \1 & \0 & \1 & \0 & \1 & \0 & \1 & \0 & \0 & \0 & \0 & \0 & \0 & \0 & \0 & \0\tabularnewline
 & $A^{1}=1$, $A^{2}=2$ & \0 & \1 & \0 & \1 & \0 & \1 & \0 & \1 & \0 & \0 & \0 & \0 & \0 & \0 & \0 & \0\tabularnewline
 & $A^{1}=2$, $A^{2}=1$ & \0 & \0 & \0 & \0 & \0 & \0 & \0 & \0 & \1 & \0 & \1 & \0 & \1 & \0 & \1 & \0\tabularnewline
 & $A^{1}=2$, $A^{2}=2$ & \0 & \0 & \0 & \0 & \0 & \0 & \0 & \0 & \0 & \1 & \0 & \1 & \0 & \1 & \0 & \1\tabularnewline
\hline 
\multirow{4}{*}{$\lambda^{1}=2$, $\lambda^{2}=1$} & $A^{1}=1$, $A^{2}=1$ & \1 & \1 & \0 & \0 & \0 & \0 & \0 & \0 & \1 & \1 & \0 & \0 & \0 & \0 & \0 & \0\tabularnewline
 & $A^{1}=1$, $A^{2}=2$ & \0 & \0 & \1 & \1 & \0 & \0 & \0 & \0 & \0 & \0 & \1 & \1 & \0 & \0 & \0 & \0\tabularnewline
 & $A^{1}=2$, $A^{2}=1$ & \0 & \0 & \0 & \0 & \1 & \1 & \0 & \0 & \0 & \0 & \0 & \0 & \1 & \1 & \0 & \0\tabularnewline
 & $A^{1}=2$, $A^{2}=2$ & \0 & \0 & \0 & \0 & \0 & \0 & \1 & \1 & \0 & \0 & \0 & \0 & \0 & \0 & \1 & \1\tabularnewline
\hline 
\multirow{4}{*}{$\lambda^{1}=2$, $\lambda^{2}=2$} & $A^{1}=1$, $A^{2}=1$ & \1 & \0 & \1 & \0 & \0 & \0 & \0 & \0 & \1 & \0 & \1 & \0 & \0 & \0 & \0 & \0\tabularnewline
 & $A^{1}=1$, $A^{2}=2$ & \0 & \1 & \0 & \1 & \0 & \0 & \0 & \0 & \0 & \1 & \0 & \1 & \0 & \0 & \0 & \0\tabularnewline
 & $A^{1}=2$, $A^{2}=1$ & \0 & \0 & \0 & \0 & \1 & \0 & \1 & \0 & \0 & \0 & \0 & \0 & \1 & \0 & \1 & \0\tabularnewline
 & $A^{1}=2$, $A^{2}=2$ & \0 & \0 & \0 & \0 & \0 & \1 & \0 & \1 & \0 & \0 & \0 & \0 & \0 & \1 & \0 & \1\tabularnewline
\end{tabular}\]The vector $P$ consists of the observed probabilities corresponding
to the row labels of the matrix, and the vector $Q$ consists of the
joint probabilities of the coupling vector $H=\left(H_{1}^{1},H_{2}^{1},H_{1}^{2},H_{2}^{2}\right)$
as indicated in the column labels of the matrix. Using the observed
probabilities of Example~\ref{Example:input-value-specific} we obtain
\[
P=[.3,.4,.1,.2,.35,.35,.15,.15,.08,.12,.32,.48,.15,.05,.35,.45]^{T}.
\]
\hfill$\square$\end{example}
\begin{thm}
If the sets of values for all $n$ inputs and all $n$ random outputs
are finite, then, using the notation of this section, 
\[
\left(A^{1},\ldots,A^{n}\right)\looparrowleft\left(\lambda^{1},\ldots,\lambda^{n}\right)
\]
holds if and only if the system of linear equations
\[
MQ=P
\]
has a solution $Q\geq0$ (the inequality meaning that the elements
of $Q$ are non-negative).
\end{thm}
Without the non-negativity constraint, the system $MQ=P$ always has
solutions, because the number of the unknowns (elements of $Q$) equals
or exceeds the rank of the matrix $M$, which can be shown to never
exceed
\[
\left(m_{1}\left(v_{1}-1\right)+1\right)\times\ldots\times\left(m_{n}\left(v_{n}-1\right)+1\right).
\]
Moreover, the structure of the matrix $M$ is such that that any solution
for $Q$ should automatically have its elements summing to 1. The
latter therefore is not a constraint. However, it is not guaranteed
that $Q\geq0$: it is possible that all solutions for $Q$ have some
of the elements negative, in which case our test establishes that
$\left(A^{1},\ldots,A^{n}\right)\not\looparrowleft\left(\lambda^{1},\ldots,\lambda^{n}\right)$.

Let us introduce a function 
\[
\mathrm{Sol}\left(M,P\right)
\]
that attains two values: ``True,'' if $MQ=P$ has a non-negative
solution, and ``False,'' if such a solution does not exist. Note
that $M$ is an argument that is determined uniquely by the format
of the problem: the number of input-output pairs and number of possible
values for inputs and outputs. The task of computing $\mathrm{Sol}\left(M,P\right)$
is a standard \emph{feasibility problem} of the area of linear algebra
called linear programming. Due to this term, the test in question
is called the \emph{linear feasibility test},
\[
\left(A^{1},\ldots,A^{n}\right)\looparrowleft\left(\lambda^{1},\ldots,\lambda^{n}\right)\textnormal{ if and only if }\mathrm{Sol}\left(M,P\right)=\textnormal{True}.
\]
It is known from linear programming that $\mathrm{Sol}\left(M,P\right)$
can always be computed.
\begin{example}
Let us apply the linear feasibility test to the matrix $M$ and vector
$P$ of Example~\ref{Example:PMQ-filled}. Using the simplex linear
programming algorithm, we obtain the solution
\[
Q=[.03,0,0,0,0,.27,.32,.08,0,.05,.12,0,0,.05,.03,.05]^{T}\geq0
\]
satisfying $MQ=P$. This means that $\mathrm{Sol}(M,P)=\text{\textquotedblleft True\textquotedblright}$,
hence $\left(A^{1},A^{2}\right)\looparrowleft\left(\lambda^{1},\lambda^{2}\right)$.
\end{example}
The (complete) marginal selectivity test mentioned in the previous
section is part of the linear feasibility test. If the former is violated,
so will also the latter. It follows from the structure of the matrix
$M$, as explained in the following example.
\begin{example}
Consider the matrix of Example~\ref{Example:PMQ-filled}. If $\left(A^{1},A^{2}\right)\looparrowleft\left(\lambda^{1},\lambda^{2}\right)$
is satisfied for a given vector $P$ of observed probabilities, then
we know that there exists a vector $Q\geq0$ such that $MQ=P$. The
marginal probabilities of $A^{1}$ and $A^{2}$ within each treatment
are obtained by summing certain elements of $P$. However, as $MQ=P$,
we can obtain these marginal probabilities also by summing certain
rows of $M$ and then multiplying these summed rows by $Q$. Thus,
if we sum the rows of $M$ corresponding to the same value of $A^{1}$
within each treatment, we obtain\[\small%
\begin{tabular}{cr|cccccccccccccccc}
\multirow{2}{*}{$\lambda^{1}=1$, $\lambda^{2}=1$} & $A^{1}=1$ & \1 & \1 & \1 & \1 & \1 & \1 & \1 & \1 & \0 & \0 & \0 & \0 & \0 & \0 & \0 & \0\tabularnewline
 & $A^{1}=2$ & \0 & \0 & \0 & \0 & \0 & \0 & \0 & \0 & \1 & \1 & \1 & \1 & \1 & \1 & \1 & \1\tabularnewline
\hline 
\multirow{2}{*}{$\lambda^{1}=1$, $\lambda^{2}=2$} & $A^{1}=1$ & \1 & \1 & \1 & \1 & \1 & \1 & \1 & \1 & \0 & \0 & \0 & \0 & \0 & \0 & \0 & \0\tabularnewline
 & $A^{1}=2$ & \0 & \0 & \0 & \0 & \0 & \0 & \0 & \0 & \1 & \1 & \1 & \1 & \1 & \1 & \1 & \1\tabularnewline
\hline 
\multirow{2}{*}{$\lambda^{1}=2$, $\lambda^{2}=1$} & $A^{1}=1$ & \1 & \1 & \1 & \1 & \0 & \0 & \0 & \0 & \1 & \1 & \1 & \1 & \0 & \0 & \0 & \0\tabularnewline
 & $A^{1}=2$ & \0 & \0 & \0 & \0 & \1 & \1 & \1 & \1 & \0 & \0 & \0 & \0 & \1 & \1 & \1 & \1\tabularnewline
\hline 
\multirow{2}{*}{$\lambda^{1}=2$, $\lambda^{2}=2$} & $A^{1}=1$ & \1 & \1 & \1 & \1 & \0 & \0 & \0 & \0 & \1 & \1 & \1 & \1 & \0 & \0 & \0 & \0\tabularnewline
 & $A^{1}=2$ & \0 & \0 & \0 & \0 & \1 & \1 & \1 & \1 & \0 & \0 & \0 & \0 & \1 & \1 & \1 & \1\tabularnewline
\end{tabular}\]As the rows corresponding to the marginal probabilities of $A^{1}$
are identical between the treatments with $\lambda^{1}=1$ and between
the treatments with $\lambda^{1}=2$, we see that the marginal distribution
of $A^{1}$ does not depend on $\lambda_{2}$. If we then sum the
rows of $M$ corresponding to the same value of $A^{2}$ within each
treatment, we obtain\[\small%
\begin{tabular}{cr|cccccccccccccccc}
\multirow{2}{*}{$\lambda^{1}=1$, $\lambda^{2}=1$} & $A^{2}=1$ & \1 & \1 & \0 & \0 & \1 & \1 & \0 & \0 & \1 & \1 & \0 & \0 & \1 & \1 & \0 & \0\tabularnewline
 & $A^{2}=2$ & \0 & \0 & \1 & \1 & \0 & \0 & \1 & \1 & \0 & \0 & \1 & \1 & \0 & \0 & \1 & \1\tabularnewline
\hline 
\multirow{2}{*}{$\lambda^{1}=1$, $\lambda^{2}=2$} & $A^{2}=1$ & \1 & \0 & \1 & \0 & \1 & \0 & \1 & \0 & \1 & \0 & \1 & \0 & \1 & \0 & \1 & \0\tabularnewline
 & $A^{2}=2$ & \0 & \1 & \0 & \1 & \0 & \1 & \0 & \1 & \0 & \1 & \0 & \1 & \0 & \1 & \0 & \1\tabularnewline
\hline 
\multirow{2}{*}{$\lambda^{1}=2$, $\lambda^{2}=1$} & $A^{2}=1$ & \1 & \1 & \0 & \0 & \1 & \1 & \0 & \0 & \1 & \1 & \0 & \0 & \1 & \1 & \0 & \0\tabularnewline
 & $A^{2}=2$ & \0 & \0 & \1 & \1 & \0 & \0 & \1 & \1 & \0 & \0 & \1 & \1 & \0 & \0 & \1 & \1\tabularnewline
\hline 
\multirow{2}{*}{$\lambda^{1}=2$, $\lambda^{2}=2$} & $A^{2}=1$ & \1 & \0 & \1 & \0 & \1 & \0 & \1 & \0 & \1 & \0 & \1 & \0 & \1 & \0 & \1 & \0\tabularnewline
 & $A^{2}=2$ & \0 & \1 & \0 & \1 & \0 & \1 & \0 & \1 & \0 & \1 & \0 & \1 & \0 & \1 & \0 & \1\tabularnewline
\end{tabular}\]and we can see that the marginal distribution of $A^{2}$ does
not depend on $\lambda^{1}$. Thus, linear feasibility test includes
the test for marginal selectivity, so if the latter is violated, the
former fails.\hfill$\square$
\end{example}
One may feel that $\mathrm{Sol}\left(M,P\right)$ is not a ``true''
function, as it requires a computer algorithm to be computed, and
it is not presented in an analytic form. Such a misgiving is not well-founded.
An analytic (or closed-form) solution is merely one that can be presented
in terms of familiar functions and operations. For example, if a solution
of a problem involves the standard normal integral
\[
N\left(t\right)=\frac{1}{\sqrt{2\pi}}\int_{-\infty}^{t}\exp\left(-z^{2}/2\right)\textnormal{d}z,
\]
the solution may or may not be called analytic depending on how familiar
and easily computable this function is. In the past, $N\left(t\right)$
could be viewed as ``less analytic'' than $\exp\left(x\right)$,
and in Napier's time $\exp\left(x\right)$ would be viewed as ``less
analytic'' than $x^{2}$. Familiarity is not a mathematical category,
and the existence of a rigorous definition of a function combined
with an algorithm allowing one to compute it to a desired level of
precision is all one needs to use it in a solution to a problem. The
computational complexity, of course, may be a concern. In our case,
however, it is known that as the size of the matrix $M$ increases,
the computational time required to compute $\mathrm{Sol}\left(M,P\right)$
increases only as a polynomial function of this size (rather than
exponentially or even faster). This makes the linear feasibility test
practical. 

It still may be of interest to see whether the linear feasibility
test could be formulated in terms of a system of equalities and inequalities
involving the entries of the vector $P$ alone. This can always be
achieved, with every linear feasibility problem. These equalities
and inequalities, in fact, can be generated by a computer algorithm
(called a \emph{facet enumeration algorithm}). \def\m{\!\!\text{-}1\!\!}\def\0{\!0\!\!}\def\1{\!1\!\!}
\begin{example}
Geometrically, the linear feasibility test checks if $P$ is within
the \emph{convex polytope} determined by points $MQ$ such that $Q\ge0$,
$\sum Q=1$. The columns of $M$ correspond to the vertices of this
polytope. A facet enumeration algorithm transforms this \emph{vertex
representation} of the polytope to the so-called \emph{half-plane
representation}, that is, to a representation of the form 
\[
M_{1}P\ge Q_{1},\; M_{2}P=Q_{2},
\]
where $M_{1},M_{2}$ are matrices and $Q_{1},Q_{2}$ are vectors. For
our $16\times16$ example matrix, this yields
\[
M_{1}=\left[\begin{array}{cccccccccccccccc}
\1 & \0 & \0 & \0 & \0 & \0 & \1 & \0 & \0 & \1 & \0 & \0 & \m & \0 & \0 & \0\\
\0 & \1 & \0 & \0 & \0 & \0 & \1 & \0 & \1 & \0 & \0 & \0 & \m & \0 & \0 & \0\\
\0 & \0 & \1 & \0 & \1 & \0 & \0 & \0 & \0 & \1 & \0 & \0 & \m & \0 & \0 & \0\\
\0 & \0 & \0 & \1 & \1 & \0 & \0 & \0 & \1 & \0 & \0 & \0 & \m & \0 & \0 & \0\\
\1 & \1 & \1 & \0 & \m & \0 & \0 & \0 & \m & \0 & \0 & \0 & \1 & \0 & \0 & \0\\
\1 & \1 & \0 & \1 & \m & \0 & \0 & \0 & \0 & \m & \0 & \0 & \1 & \0 & \0 & \0\\
\1 & \0 & \1 & \1 & \0 & \0 & \m & \0 & \m & \0 & \0 & \0 & \1 & \0 & \0 & \0\\
\0 & \1 & \1 & \1 & \0 & \0 & \m & \0 & \0 & \m & \0 & \0 & \1 & \0 & \0 & \0
\end{array}\right]\qquad Q_{1}=\left[\begin{array}{c}
0\\
0\\
0\\
0\\
0\\
0\\
0\\
0
\end{array}\right]
\]
and
\[
M_{2}=\left[\begin{array}{cccccccccccccccc}
\1 & \1 & \0 & \0 & \m & \m & \0 & \0 & \0 & \0 & \0 & \0 & \0 & \0 & \0 & \0\\
\1 & \0 & \1 & \0 & \0 & \0 & \0 & \0 & \m & \0 & \m & \0 & \0 & \0 & \0 & \0\\
\0 & \1 & \0 & \1 & \0 & \0 & \0 & \0 & \0 & \m & \0 & \m & \0 & \0 & \0 & \0\\
\0 & \0 & \1 & \1 & \0 & \0 & \m & \m & \0 & \0 & \0 & \0 & \0 & \0 & \0 & \0\\
\0 & \0 & \0 & \0 & \1 & \0 & \1 & \0 & \0 & \0 & \0 & \0 & \m & \0 & \m & \0\\
\0 & \0 & \0 & \0 & \0 & \0 & \0 & \0 & \1 & \1 & \0 & \0 & \m & \m & \0 & \0\\
\1 & \1 & \1 & \1 & \0 & \0 & \0 & \0 & \0 & \0 & \0 & \0 & \0 & \0 & \0 & \0\\
\1 & \1 & \1 & \1 & \m & \0 & \m & \0 & \m & \m & \0 & \0 & \1 & \0 & \0 & \m
\end{array}\right]\qquad Q_{2}=\left[\begin{array}{c}
0\\
0\\
0\\
0\\
0\\
0\\
1\\
0
\end{array}\right].
\]
The equations $M_{2}P=Q_{2}$ of this representation always correspond
to the marginal selectivity constraints. Thus, a vector $P$ of observed
probabilities satisfying marginal selectivity satisfies selective
influences if and only if $M_{1}P\ge Q_{1}.$ Assuming marginal selectivity,
$M_{1}P\ge Q_{1}$ can in this case also be simplified into the four
double-inequalities
\[
0\le p_{i\cdot}+p_{\cdot j}+p_{i'j'}-p_{ij'}-p_{i'j'}-p_{i'j}\le1,\qquad i\ne i',\, j\ne j',
\]
where we denote
\begin{align*}
p_{i\cdot} & =\Pr(A^{1}=1)_{\phi=(\lambda^{1}=i,\lambda^{2}=\cdot)},\\
p_{\cdot j} & =\Pr(A^{2}=1)_{\phi=(\lambda^{1}=\cdot,\lambda^{2}=j)},\\
p_{ij} & =\Pr(A^{1}=1,A^{2}=1)_{\phi=(\lambda^{1}=i,\lambda^{2}=j)}
\end{align*}
(the definition of $p_{i\cdot}$ and $p_{\cdot j}$ presupposes marginal
selectivity). These are known as the \emph{Bell/CHSH/Fine inequalities}
in quantum mechanics. \hfill$\square$
\end{example}
In the same way, the representation as inequalities can be obtained
for any linear feasibility test matrix $M$. It should be noted, however,
that the number of the inequalities increases explosively as the size
of the matrix $M$ increases. Thus, for three pairs of completely
crossed binary inputs and three binary random outputs, the number
of independent equalities representing marginal selectivity is 42,
and the number of inequalities is 53792. From a practical point of
view, therefore, computing $\mathrm{Sol}\left(M,P\right)$ directly
is a better approach in all but the simplest cases.

\section{Distance Tests}

Let us establish some general terminology. A \emph{pseudo-quasi-metric}
(or \emph{p.q.-metric}, for short) on a nonempty set $X$ is defined
as a function $d:X\times X\rightarrow\mathbb{R}^{+}$ (set of non-negative
real numbers), such that, for any $x,y,z\in X$, 
\begin{lyxlist}{00.00.0000}
\item [{(1)}] (zero$\:$property) $d\left(x,x\right)=0$,
\item [{(2)}] (triangle inequality) $d\left(x,y\right)+d\left(y,z\right)\geq d\left(x,z\right)$.
\end{lyxlist}
A p.q.-metric that satisfies, in addition, 
\begin{lyxlist}{00.00.0000}
\item [{(3)}] (symmetry) $d\left(x,y\right)=d\left(y,x\right)$,
\end{lyxlist}
is called a \emph{pseudo-metric.} A p.q.-metric that satisfies
\begin{lyxlist}{00.00.0000}
\item [{(4)}] (positivity) if $x\not=y$, then $d\left(x,y\right)>0$,
\end{lyxlist}
is called a \emph{quasi-metric.} Finally, a p.q.-metric that satisfies
both (3) and (4) is called a \emph{metric}. The terminology is not
well-established and varies from one area or application to another. 
\begin{rem}
To refer to the value $d\left(x,y\right)$ of a metric, pseudo-metrics,
quasi-metrics, or a p.q.-metric at a specific pair of points $\left(x,y\right)$,
one usually uses the generic term ``distance,'' adding the corresponding
prefixes (pseudo, quasi, or p.q.) only if it is required for disambiguation.
Thus, the value of a p.q.-metric for a specific pair $\left(x,y\right)$
can be called the distance from $x$ to $y$, or the p.q.-distance
from $x$ to $y$. (For pseudo-metrics, ``from $x$ to $y$'' can
be replaced with ``between $x$ and $y$.'') The term ``distance''
can also be used (with or without the prefixes) to refer to the functions
themselves. Therefore ``p.q.-metric tests'' below can also be referred
to as ``distance tests'' or ``p.q.-distance tests.'' 
\end{rem}
The nature of the set $X$ in the definition is entirely arbitrary.
We are interested in a set of jointly distributed random variables,
that is, those representable as functions of one and the same random
variable. A p.q.-metric on such a set is a function $d$ mapping pairs
of random variables into non-negative real numbers, such that $d\left(R,R\right)=0$
and $d\left(R^{1},R^{2}\right)+d\left(R^{2},R^{3}\right)\geq$d$\left(R^{1},R^{3}\right)$,
for any random variables $R^{1},R^{2},R^{3}$ in the set. We assume
that $d\left(R^{1},R^{2}\right)$ is entirely determined by the joint
distribution of $\left(R^{1},R^{2}\right)$. In other words, it does
not depend on the identifying label of the pair (or on how $R^{1}$
and $R^{2}$ are presented as functions of a common random variable). 

An immediate consequence (and generalization) of the triangle inequality
is the following \emph{chain inequality}: if $R^{1},\ldots,R^{l}$
are elements of $X$ ($l\geq3$), not necessarily distinct, then
\[
d\left(R^{1},R^{l}\right)\leq\sum_{i=2}^{l}d\left(R^{i-1},R^{i}\right).
\]
This inequality, as it turns out, can be utilized to construct tests
of selective influences. 

Suppose that the random outputs $\left(A_{\phi}^{1},\ldots,A_{\phi}^{n}\right)$
across all $\phi\in\Phi$ belong to a certain type, or class of random
variables (e.g., those in the narrow sense, or with a finite number
of values, etc.). We continue to consider, for simplicity, inputs
with finite number of values each. We know that $\left(A^{1},\ldots,A^{n}\right)\looparrowleft\left(\lambda^{1},\ldots,\lambda^{n}\right)$
if and only if there exists a reduced coupling vector $H$. Assuming
that it does exist, its elements are of the same type, or class, as
$\left(A_{\phi}^{1},\ldots,A_{\phi}^{n}\right)$, and any $l\geq3$
of these elements,
\[
H_{j_{1}}^{k_{1}},H_{j_{2}}^{k_{2}},\ldots,H_{j_{l}}^{k_{l}}
\]
can be used to form a chain inequality,
\[
d\left(H_{j_{1}}^{k_{1}},H_{j_{l}}^{k_{l}}\right)\leq\sum_{i=2}^{l}d\left(H_{j_{i-1}}^{k_{i-1}},H_{j_{i}}^{k_{i}}\right).
\]
Let us choose these elements of $H$ so that $\lambda^{k_{1}}=j_{1}$
and $\lambda^{k_{l}}=j_{l}$ belong to some allowable treatment $\phi_{1k}$,
and each pair $\lambda^{k_{i-1}}=j_{i-1},\lambda^{k_{i}}=j_{i}$ belongs
to some allowable treatment $\phi_{i-1,i}$ ($i=2,\ldots,l$). The
allowable treatments $\phi_{1k},\phi_{12},\ldots,\phi_{l-1,l}$ need
not be pairwise distinct. Such a sequence of input values,
\[
\lambda^{k_{1}}=j_{1},\lambda^{k_{2}}=j_{2},\ldots,\lambda^{k_{l}}=j_{l}
\]
is called \emph{treatment-realizable}. This choice ensures that
\[
\left(H_{j_{1}}^{k_{1}},H_{j_{l}}^{k_{l}}\right)\sim\left(A_{\phi_{1k}}^{k_{1}},A_{\phi_{1l}}^{k_{l}}\right)
\]
and
\[
\left(H_{j_{i-1}}^{k_{i-1}},H_{j_{i}}^{k_{i}}\right)\sim\left(A_{\phi_{i-1,i}}^{k_{i-1}},A_{\phi_{i-1,i}}^{k_{i}}\right),\textnormal{ for }i=2,\ldots,l.
\]
But then
\[
d\left(H_{j_{1}}^{k_{1}},H_{j_{l}}^{k_{l}}\right)=d\left(A_{\phi_{1k}}^{k_{1}},A_{\phi_{1l}}^{k_{l}}\right)
\]
and
\[
d\left(H_{j_{i-1}}^{k_{i-1}},H_{j_{i}}^{k_{i}}\right)=d\left(A_{\phi_{i-1,i}}^{k_{i-1}},A_{\phi_{i-1,i}}^{k_{i}}\right),\textnormal{ for }i=2,\ldots,l,
\]
whence the chain inequality can be rewritten using only observable
pairwise distributions,
\[
d\left(A_{\phi_{1k}}^{k_{1}},A_{\phi_{1l}}^{k_{l}}\right)\leq\sum_{i=2}^{l}d\left(A_{\phi_{i-1,i}}^{k_{i-1}},A_{\phi_{i-1,i}}^{k_{i}}\right).
\]
This inequality is a necessary condition for the existence of $H$.
If it is found violated for at least one treatment-realizable sequence
of input values, then the existence of $H$ is ruled out, and one
should conclude that $\left(A^{1},\ldots,A^{n}\right)\not\looparrowleft\left(\lambda^{1},\ldots,\lambda^{n}\right).$  

There are numerous ways of constructing p.q.-metrics for jointly distributed
random variables. We will confine our consideration to only two examples.

If all random outputs have one and the same set of possible values
$S$, then one way of creating a p.q.-metric on a set $X$ of such
random variables is to use any p.q.-metric $D$ on $S$ and put, for
any random variables $Q,R\in X$,
\[
d\left(Q,R\right)=\textnormal{E}\left[D\left(Q,R\right)\right].
\]
The right-hand expression is the expected value of the random variable
$D\left(Q,R\right)$. The underlying assumption is, of course, that
this random variable is well-defined (that is, $D$ is a measurable
function from $S\times S$ to non-negative real numbers), and that
its expectation is finite. It can easily be proved then that $d$
is a p.q.-metric on $X$.

As a simple example, consider the p.q.-metric
\[
D\left(x,y\right)=\left\{ \begin{array}{cc}
\left|x-y\right|^{p} & \textnormal{if }x<y\\
0 & \textnormal{otherwise}
\end{array}\right.
\]
on the set of real numbers, with $0\leq p\leq 1$ (a power exponent). It
is a p.q.-metric because $D\left(x,x\right)=0$, and 
\[
D\left(x,y\right)+D\left(y,z\right)\geq D\left(x,z\right),
\]
as one can prove by considering various arrangements of numbers $x,y,z$.
Using $D$ one can construct a p.q.-metric for any set $X$ of random
variables whose (common) set of possible values is a subset of reals.
Let this set be a subset of integers. Then the p.q.-metric on $X$
derived from $D$ is 
\[
d_{p}\left(Q,R\right)=\sum_{q<r}\left|q-r\right|^{p}p\left(q,r\right),
\]
where 
\[
p\left(q,r\right)=\Pr\left(Q=q,R=r\right).
\]
 
\begin{example}
\label{example:d1}Let the outputs $A^{1},A^{2}$ have the following
distributions for treatments in $\Lambda^{1}\times\Lambda^{2}=\{1,2\}\times\{1,2\}$:\[\small%
\begin{tabular}{r|ccclr|ccc}
\multicolumn{4}{l}{$\lambda^{1}=1$, $\lambda^{2}=1$} &  & \multicolumn{4}{l}{$\lambda^{1}=1$, $\lambda^{2}=2$}\tabularnewline
 & $A^{2}=0$ & $A^{2}=1$ & $A^{2}=2$ &  &  & $A^{2}=0$ & $A^{2}=1$ & $A^{2}=2$\tabularnewline
\cline{1-4} \cline{6-9} 
$A^{1}=0$ & .24 & .07 & 0 &  & $A^{1}=0$ & .24 & .07 & 0\tabularnewline
$A^{1}=2$ & .07 & .24 & .07 &  & $A^{1}=2$ & .07 & .24 & .07\tabularnewline
$A^{1}=4$ & 0 & .07 & .24 &  & $A^{1}=4$ & 0 & .07 & .24\tabularnewline
\multicolumn{1}{r}{} &  &  &  &  & \multicolumn{1}{r}{} &  &  & \tabularnewline
\multicolumn{4}{l}{$\lambda^{1}=2$, $\lambda^{2}=1$} &  & \multicolumn{4}{l}{$\lambda^{1}=2$, $\lambda^{2}=2$}\tabularnewline
 & $A^{2}=0$ & $A^{2}=1$ & $A^{2}=2$ &  &  & $A^{2}=0$ & $A^{2}=1$ & $A^{2}=2$\tabularnewline
\cline{1-4} \cline{6-9} 
$A^{1}=0$ & .24 & .07 & 0 &  & $A^{1}=0$ & 0 & .07 & .24\tabularnewline
$A^{1}=2$ & .07 & .24 & .07 &  & $A^{1}=2$ & .07 & .24 & .07\tabularnewline
$A^{1}=4$ & 0 & .07 & .24 &  & $A^{1}=4$ & .24 & .07 & 0\tabularnewline
\end{tabular}\]Let us put $p=1$ and compute the values of the $d_{1}$-p.q.-metric.
For any $\lambda^{1},\lambda^{2}$ here,
\[
\begin{array}{l}
d_{1}(A^{1},A^{2})=\sum_{a_{1}<a_{2}}\left|a_{1}-a_{2}\right|^{1}p\left(a_{1},a_{2}\right)=\left|1-0\right|^{1}p\left(0,1\right)+\left|2-0\right|^{1}p\left(0,2\right)\end{array}
\]
and
\begin{align*}
d_{1}(A^{2},A^{1}) & =\sum_{a_{2}<a_{1}}\left|a_{1}-a_{2}\right|^{1}p\left(a_{1},a_{2}\right)\\
 & =\left|2-0\right|^{1}p\left(2,0\right)+\left|4-0\right|^{1}p\left(4,0\right)+\ldots+\left|4-2\right|^{1}p\left(4,2\right).
\end{align*}
The calculations yield the following distances:\[\small%
\begin{tabular}{c|cccc}
 & $\!\lambda^{1}=1$, $\lambda^{2}=1\!$ & $\!\lambda^{1}=1$, $\lambda^{2}=2\!$ & $\!\lambda^{1}=2$, $\lambda^{2}=1\!$ & $\!\lambda^{1}=2$, $\lambda^{2}=2\!$\tabularnewline
\hline 
$d_{p=1}(A^{1},A^{2})$ & .07 & .07 & .07 & .55\tabularnewline
$d_{p=1}(A^{2},A^{1})$ & 1.07 & 1.07 & 1.07 & 1.55\tabularnewline
\end{tabular}\]Using this table, all possible distance test inequalities are of
the form $a\le b+c+d$, where $a$, $b$, and $d$ belong to one row
and $c$ to another, provided all four values are in distinct columns.
It is easy to see that all the inequalities are passed.\hfill$\square$
\end{example}
P.q.-metrics can be introduced directly in probabilistic terms rather
derived from ``deterministic'' metrics on sets of possible values.
Consider, as an example, the following construction. Let $\left(S^{1},\Sigma^{1}\right),\ldots,\left(S^{m},\Sigma^{m}\right)$
be the sets of possible values with sigma-algebras for random variables
$R^{1},\ldots,R^{m}$, respectively, and let us partition each $S^{k}$
into $l_{k}>1$ measurable subsets $S^{1k},\ldots,S^{l_{k}k}\in\Sigma^{k}$.
It follows that the joint probabilities of any pair $S^{ik},S^{i'k'}$,
\[
\Pr\left(R^{k}\in S^{ik},R^{k'}\in S^{i'k'}\right),
\]
are well defined. It can easily be proved that the function
\[
d_{class}\left(R^{k},R^{k'}\right)=\sum_{i<i'}\Pr\left(R^{k}\in S^{ik},R^{k'}\in S^{i'k'}\right)
\]
is a p.q.-metric. It is called a \emph{classification p.q.-metric},
and it can be applied to all types of random variables without restrictions. 
\begin{example}
Consider the case with two real-valued random variables $R^{1},R^{2}$
and define the partition of $S^{1}=\mathbb{R}$ and $S^{2}=\mathbb{R}$
as, respectively, 
\[
S^{11}=(-\infty,x),\; S^{21}=[x,\infty)
\]
and 
\[
S^{12}=(-\infty,y),\; S^{22}=[y,\infty).
\]
Then, the classification distance is simply
\[
d_{class}(R^{1},R^{2})=\Pr\left(R^{1}\in S^{11},R^{2}\in S^{22}\right)=\Pr\left(R^{1}<x,R^{2}\ge y\right).
\]
Different choices of $x,y$ give us different classification distances.\hfill$\square$\end{example}
\begin{rem}
\label{REM: A-classification-p.q.-metric}A classification p.q.-metric
can also be viewed as a limit case of the metric $d_{p}$ introduced
above, provided we first map by a measurable function $f_{k}$ each
$S^{k}$ into a set $\left\{ 1,\ldots,l_{k}\right\} $, and then define
all the transformed random variables $f_{k}\left(R^{k}\right)$ as
distributed on $\left\{ 1,\ldots,l\right\} $, with $l=\max\left(l_{1},\ldots,l_{m}\right)$.
The latter is always possible by assigning to the ``redundant''
integers probability zero. Following this transformation and equalization
of domains, $d_{class}$ is obtained as $d_{p=0}$. Another way of
introducing the classification metric is as a special case of an \emph{order-distance}.
Without elaborating, the latter involves a relation of strict order
$\prec$ between values of one random variable and values of another.
The order-distance is defined as
\[
d_{ord}\left(Q,R\right)=\Pr\left(Q\prec R\right).
\]

\end{rem}
Recall that a sequence $\lambda^{k_{1}}=j_{1},\lambda^{k_{2}}=j_{2},\ldots,\lambda^{k_{l}}=j_{l}$
of input values is treatment-realizable if $\left\{ \lambda^{k_{1}}=j_{1},\lambda^{k_{k}}=j_{k}\right\} $
and $\left\{ \lambda^{k_{i-1}}=j_{i-1},\lambda^{k_{i}}=j_{i}\right\} $
for $i=2,\ldots,l$ belong to allowable treatments. If the elements
of all these pairs are distinct, and if these pairs are the only subsequences
of more than one element that have the property of being a subset
of an allowable treatment, then the sequence is called \emph{irreducible}.
It turns out that one only has to check the chain inequalities for
irreducible sequences: these inequalities are satisfied for all treatment-realizable
sequences if and only if they are satisfied for all irreducible ones. 

The set of irreducible sequences may be significantly smaller than
the set of all treatment-realizable sequences. Thus, it can be shown
that if the set $\Phi$ consists of all possible combinations of input
values, then the only irreducible sequences are quadruples of the
form
\[
\lambda^{k}=j_{1},\lambda^{k'}=j_{2},\lambda^{k}=j_{3},\lambda^{k'}=j_{4},
\]
with $k\not=k'$, $j_{1}\not=j_{3}$ and $j_{2}\not=j_{4}$. The only
inequalities to check then are of the form, 
\[
d\left(A_{\phi_{14}}^{k},A_{\phi_{14}}^{k'}\right)\leq d\left(A_{\phi_{12}}^{k},A_{\phi_{12}}^{k'}\right)+d\left(A_{\phi_{23}}^{k'},A_{\phi_{23}}^{k}\right)+d\left(A_{\phi_{34}}^{k},A_{\phi_{34}}^{k'}\right),
\]
where $\phi_{14},\phi_{12},\phi_{23},\phi_{34}$ are any allowable
treatments that contain, respectively,
\begin{align*}
 & \left\{ \lambda^{k}=j_{1},\lambda^{k'}=j_{4}\right\} ,\ \left\{ \lambda^{k}=j_{1},\lambda^{k'}=j_{2}\right\} ,\\
 & \left\{ \lambda^{k'}=j_{2},\lambda^{k}=j_{3}\right\} ,\ \left\{ \lambda^{k}=j_{3},\lambda^{k'}=j_{4}\right\} .
\end{align*}

\section{(Non)Invariance of tests with respect to transformations }

In this section we introduce another class of tests of selective influences,
called \emph{cosphericity tests}. Prior to introducing them, however,
we should discuss an important issue. 

We know from Section \ref{sec:Properties-of-selective} that if $\left(A^{1},\ldots,A^{n}\right)\looparrowleft\left(\lambda^{1},\ldots,\lambda^{n}\right)$,
then $\left(B^{1},\ldots,B^{n}\right)\looparrowleft\left(\lambda^{1},\ldots,\lambda^{n}\right)$,
where the $B$'s are input-value-specific transformations of the $A'$s,
that is,

\[
B_{\phi}^{1}=g_{1}\left(\lambda^{1},A_{\phi}^{1}\right),\ldots,B_{\phi}^{n}=g_{n}\left(\lambda^{n},A_{\phi}^{n}\right),
\]
for all $\phi=\left(\lambda^{1},\ldots,\lambda^{n}\right)\in\Phi$.
It follows that if a test provides a necessary condition for selective
influences, then its failure for any of the input-value-specific transformations
of $\left(A^{1},\ldots,A^{n}\right)_{\phi}$ establishes $\left(A^{1},\ldots,A^{n}\right)\not\looparrowleft\left(\lambda^{1},\ldots,\lambda^{n}\right)$.
If the outcome of a test is not invariant with respect to some of
such transformations, this consideration automatically expands this
test into a multitude of tests, one for each of these transformations.
This may enormously increase the ability of a test to detect violations
of selective influences. This might sound paradoxical, or at least
unexpected, but this is generally true for any test that provides
a necessary but not sufficient condition for a tested proposition:
the lack of invariance in the test's outcome with respect to transformations
that preserve the tested proposition is an advantage rather than a
drawback. 
\begin{rem}
If a test provides a sufficient condition for $\left(A^{1},\ldots,A^{n}\right)\looparrowleft\left(\lambda^{1},\ldots,\lambda^{n}\right)$,
and it is not invariant with respect to input-value-specific transformations,
then one should apply it to a variety of $\left(B^{1},\ldots,B^{n}\right)_{\phi}$
from which $\left(A^{1},\ldots,A^{n}\right)_{\phi}$ can be obtained
by such a transformation. At the time this is written (end of 2012),
we do not have nontrivial tests that provide sufficient but not necessary
conditions. If a test is a criterion when applied to input-output
pairs of a particular type, then its (non)invariance with respect
to transformations is immaterial for establishing or rejecting selective
influences for original random variables (although transformed ones
may be of interest for their own sake).
\end{rem}
Of the two distance tests considered in the previous section, $d_{p}$-test
is not invariant (for any fixed $p$) with respect to numerical transformations
of the random outputs.
\begin{example}
\label{example:distance-test-transformations}

Continuing Example \ref{example:d1}, let us transform the outputs
$A^{1},A^{2}$ as $B^{1}=g_{1}(A^{1})$, $B^{2}=g_{2}(A^{2})$, where
$g_{1}$ is given by $0\mapsto2$, $2\mapsto1$, $4\mapsto1$ and
$g_{2}$ is given by $0\mapsto2$, $1\mapsto1$, $2\mapsto1$. We
get the joint distributions\[\small%
\begin{tabular}{r|cccr|cc}
$\lambda^{1}=1,\lambda^{2}=1$ & $B^{2}=1$ & $B^{2}=2$ &  & $\lambda^{1}=1,\lambda^{2}=2$ & $B^{2}=1$ & $B^{2}=2$\tabularnewline
\cline{1-3} \cline{5-7} 
$B^{1}=1$ & .62 & .07 &  & $B^{1}=1$ & .62 & .07\tabularnewline
$B^{1}=2$ & .07 & .24 &  & $B^{1}=2$ & .07 & .24\tabularnewline
\multicolumn{1}{r}{} &  &  &  & \multicolumn{1}{r}{} &  & \tabularnewline
$\lambda^{1}=2,\lambda^{2}=1$ & $B^{2}=1$ & $B^{2}=2$ &  & $\lambda^{1}=2,\lambda^{2}=2$ & $B^{2}=1$ & $B^{2}=2$\tabularnewline
\cline{1-3} \cline{5-7} 
$B^{1}=1$ & .62 & .07 &  & $B^{1}=1$ & .38 & .31\tabularnewline
$B^{1}=2$ & .07 & .24 &  & $B^{1}=2$ & .31 & 0\tabularnewline
\end{tabular}\]and the corresponding $d_{p=1}$ distances are\[\small%
\begin{tabular}{c|cccc}
 & $\!\lambda^{1}=1$, $\lambda^{2}=1\!$ & $\!\lambda^{1}=1$, $\lambda^{2}=2\!$ & $\!\lambda^{1}=2$, $\lambda^{2}=1\!$ & $\!\lambda^{1}=2$, $\lambda^{2}=2\!$\tabularnewline
\hline 
$d_{p=1}(B^{1},B^{2})$ & .07 & .07 & .07 & .31\tabularnewline
$d_{p=1}(B^{2},B^{1})$ & .07 & .07 & .07 & .31\tabularnewline
\end{tabular}\]Now the distance test inequality $.31\le.07+.07+.07=.21$ fails
implying $\left(B^{1},B^{2}\right)\not\looparrowleft\left(\lambda^{1},\lambda^{2}\right)$
which in turn implies $\left(A^{1},A^{2}\right)\not\looparrowleft\left(\lambda^{1},\lambda^{2}\right)$.
Thus, the $d_{p}$-test is not invariant with respect to transformations
of the variables.\hfill$\square$
\end{example}
The second distance test considered in the previous section, $d_{class}$-test,
is invariant (for any given partition scheme) with respect to any
transformations of the possible values of random outputs. The obvious
proviso for this statement is that a transformed value is always classified
into a partition with the same number as the original value. If this
proviso is violated, it would amount to changing the partition scheme
for the original outputs. The power of the $d_{class}$-test to detect
violations of selective influences does not come from different transformations.
Rather it comes from complete flexibility in the partitioning scheme.
Another way of looking at this test (see Remark \ref{REM: A-classification-p.q.-metric})
is that a transformation of the random outputs (different mappings
into natural numbers) is built into the identity of the test. If the
transformation changes, we apply a different test.
\begin{example}
Consider the system $(A^{1},A^{2})$ of Example~\ref{example:d1}.
Let us partition $S^{1}$ into $S^{11}=\{0\}$, $S^{21}=\{2,4\}$,
and $S^{2}$ into $S^{12}=\{0,1\}$, $S^{22}=\{$2\}. We obtain the
following joint probabilities for the partition memberships $A^{k}\in S^{ik}$:\[\small%
\begin{tabular}{r|cccr|cc}
$\lambda^{1}=1,\lambda^{2}=1$ & $\! A^{2}\in S^{21}\!$ & $\! A^{2}\in S^{22}\!$ &  & $\lambda^{1}=1,\lambda^{2}=2$ & $\! A^{2}\in S^{21}\!$ & $\! A^{2}\in S^{22}\!$\tabularnewline
\cline{1-3} \cline{5-7} 
$A^{1}\in S^{11}$ & .31 & 0 &  & $A^{1}\in S^{11}$ & .31 & 0\tabularnewline
$A^{1}\in S^{12}$ & .38 & .31 &  & $A^{1}\in S^{12}$ & .38 & .31\tabularnewline
\multicolumn{1}{r}{} &  &  &  & \multicolumn{1}{r}{} &  & \tabularnewline
$\lambda^{1}=2,\lambda^{2}=1$ & $\! A^{2}\in S^{21}\!$ & $\! A^{2}\in S^{22}\!$ &  & $\lambda^{1}=2,\lambda^{2}=2$ & $\! A^{2}\in S^{21}\!$ & $\! A^{2}\in S^{22}\!$\tabularnewline
\cline{1-3} \cline{5-7} 
$A^{1}\in S^{11}$ & .31 & 0 &  & $A^{1}\in S^{11}$ & .07 & .24\tabularnewline
$A^{1}\in S^{12}$ & .38 & .31 &  & $A^{1}\in S^{12}$ & .62 & .07\tabularnewline
\end{tabular}\]This yields the classification distances\[\small%
\begin{tabular}{c|cccc}
 & $\!\lambda^{1}=1$, $\lambda^{2}=1\!$ & $\!\lambda^{1}=1$, $\lambda^{2}=2\!$ & $\!\lambda^{1}=2$, $\lambda^{2}=1\!$ & $\!\lambda^{1}=2$, $\lambda^{2}=2\!$\tabularnewline
\hline 
$d_{class}(A^{1},A^{2})$ & 0 & 0 & 0 & .24\tabularnewline
$d_{class}(A^{2},A^{1})$ & .38 & .38 & .38 & .62\tabularnewline
\end{tabular}\] which can be seen to satisfy all distance test inequalities, as
in Example~\ref{example:d1}. 

Consider now the partitioning of $S^{1}$ into $S^{11}=\{0,2\}$,
$S^{21}=\{4\}$, and of $S^{2}$ into $S^{12}=\{0,1\},$ $S^{22}=\{2\}$.
The partition membership indicator $B^{k}$ (given by $B^{k}=i$ when
$A^{k}\in S^{ik}$) corresponds to the transformed variables $B^{k}$
of Example~\ref{example:distance-test-transformations}. As a result,
we get the same joint distribution tables as there. We know that $d_{class}$
corresponds to $d_{p=0}$ (see Remark \ref{REM: A-classification-p.q.-metric}),
and it is easy to see that $d_{p=0}$ is identical to $d_{p=1}$ when
the sets are partitioned into only two classes each. Therefore $d_{class}$
distance table we obtain is identical to the $d_{p=1}$ table shown
in Example~\ref{example:distance-test-transformations}, and we conclude
that the $d_{class}$-distance test fails, implying $\left(A^{1},A^{2}\right)\not\looparrowleft\left(\lambda^{1},\lambda^{2}\right)$.\hfill$\square$
\end{example}
We conclude this section by presenting a test based on pairwise correlation
between random outputs. It is called the \emph{cosphericity test},
and confined to random variables for which conventional correlations
can be computed. These are all variables that are defined (or can
be redefined) on the set of real numbers with the Lebesgue sigma-algebra.
Discrete random variables can always be redefined to fall within this
category.

The primary application of the cosphericity test is to two input-output
pairs, with two values per input, and all four treatments allowable.
That is, we test the assumption $\left(A^{1},A^{2}\right)\looparrowleft\left(\lambda^{1},\lambda^{2}\right)$,
with $\Lambda^{1}=\left\{ 1,2\right\} $, $\Lambda^{2}=\left\{ 1,2\right\} $,
and allowable treatments $\phi_{11}=\left(\lambda_{1}^{1},\lambda_{1}^{2}\right)$,
$\phi_{12}=\left(\lambda_{1}^{1},\lambda_{2}^{2}\right)$, etc. The
use of the test for larger designs will be discussed later. 

Denote the correlation between $A_{\phi_{ij}}^{1}$ and $A_{\phi_{ij}}^{2}$
(as the two are jointly distributed) by $\rho_{ij}$, $i,j\in\left\{ 1,2\right\} $.
The cosphericity test is the proposition: if $\left(A^{1},A^{2}\right)\looparrowleft\left(\lambda^{1},\lambda^{2}\right)$
, then
\[
\left\vert \rho_{11}\rho_{12}-\rho_{21}\rho_{22}\right\vert \leq\sqrt{1-\left(\rho_{11}\right)^{2}}\sqrt{1-\left(\rho_{12}\right)^{2}}+\sqrt{1-\left(\rho_{21}\right)^{2}}\sqrt{1-\left(\rho_{22}\right)^{2}}.
\]
Superscript $2$ here indicates squaring. If this inequality is violated,
then the initial assumption $\left(A^{1},A^{2}\right)\looparrowleft\left(\lambda^{1},\lambda^{2}\right)$
should be rejected. 

The explanation for the name ``cosphericity'' is this: the inequality
above holds if and only if one can place four points, $\mathbf{A}_{1},\mathbf{A}_{2},\mathbf{B}_{1},\mathbf{B}_{2},$
on the surface of a unit sphere (in the Euclidean three-dimensional
space) centered at point $\mathbf{O}$, so that
\begin{align*}
 & \cos\angle\mathbf{A}_{1}\mathbf{O}\mathbf{B}_{1}=\rho_{11},\quad\cos\angle\mathbf{A}_{1}\mathbf{O}\mathbf{B}_{2}=\rho_{12},\\
 & \cos\angle\mathbf{A}_{2}\mathbf{O}\mathbf{B}_{1}=\rho_{21},\quad\cos\angle\mathbf{A}_{2}\mathbf{O}\mathbf{B}_{2}=\rho_{22}.
\end{align*}

\begin{center}\includegraphics[width=7cm]{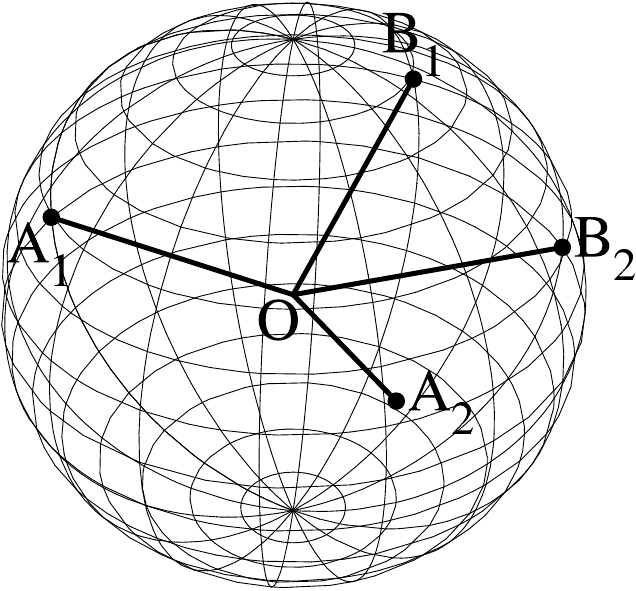}\end{center}
\begin{example}
\label{example:cosphericity-test-passed}Consider the following output
distributions of $A^{1},A^{2}$ for the treatments in $\Lambda^{1}\times\Lambda^{2}=\{1,2\}\times\{1,2\}$:\[\small%
\begin{tabular}{r|ccccr|ccc}
\multicolumn{4}{l}{$\lambda^{1}=1$, $\lambda^{2}=1$} &  & \multicolumn{4}{l}{$\lambda^{1}=1$, $\lambda^{2}=2$}\tabularnewline
 & $A^{2}=0$ & $A^{2}=1$ & $A^{2}=5$ &  &  & $A^{2}=0$ & $A^{2}=1$ & $A^{2}=5$\tabularnewline
\cline{1-4} \cline{6-9} 
$A^{1}=0$ & .24 & .07 & 0 &  & $A^{1}=0$ & .24 & .07 & 0\tabularnewline
$A^{1}=1$ & .07 & .24 & .07 &  & $A^{1}=1$ & .07 & .24 & .07\tabularnewline
$A^{1}=5$ & 0 & .07 & .24 &  & $A^{1}=5$ & 0 & .07 & .24\tabularnewline
\multicolumn{1}{r}{} &  &  &  &  & \multicolumn{1}{r}{} &  &  & \tabularnewline
\multicolumn{4}{l}{$\lambda^{1}=2$, $\lambda^{2}=1$} &  & \multicolumn{4}{l}{$\lambda^{1}=2$, $\lambda^{2}=2$}\tabularnewline
 & $A^{2}=0$ & $A^{2}=1$ & $A^{2}=5$ &  &  & $A^{2}=0$ & $A^{2}=1$ & $A^{2}=5$\tabularnewline
\cline{1-4} \cline{6-9} 
$A^{1}=0$ & .24 & .07 & 0 &  & $A^{1}=0$ & 0 & .07 & .24\tabularnewline
$A^{1}=1$ & .07 & .24 & .07 &  & $A^{1}=1$ & .07 & .24 & .07\tabularnewline
$A^{1}=5$ & 0 & .07 & .24 &  & $A^{1}=5$ & .24 & .07 & 0\tabularnewline
\end{tabular}\]The correlations coefficients of the four distributions are $\rho_{11}=\rho_{12}=\rho_{21}\approx.7299$
and $\rho_{22}\approx-.6322$. Substituting these in the cosphericity
test, we obtain
\begin{align*}
.9942 & \approx|.7299\cdot.7299-.7299(-.6322)|\\
 & \le\sqrt{1-.7299^{2}}\sqrt{1-.7299^{2}}+\sqrt{1-.7299^{2}}\sqrt{1-.6322^{2}}\approx.9969,
\end{align*}
so the test is passed. \hfill$\square$
\end{example}
Correlation between two random variables is not invariant with respect
to any but affine transformations of the random variables. This allows
one to expand the single cosphericity test into a potential infinity
of tests, corresponding to different nonlinear input-value-specific
transformations $g_{1}\left(\lambda^{1},A_{\phi}^{1}\right)$ and
$g_{2}\left(\lambda^{2},A_{\phi}^{2}\right)$. An interesting fact
is that if, by means of some \emph{reversible }transformations $g_{1},g_{2}$
the random variables $\left(A^{1},A^{2}\right)$$_{\phi}$ can be
made bivariate-normally distributed at all four treatments, then the
cosphericity test performed on thus transformed random outputs provides
both a necessary and sufficient condition for $\left(A^{1},A^{2}\right)\looparrowleft\left(\lambda^{1},\lambda^{2}\right)$.
\begin{example}
The system of Example~\ref{example:cosphericity-test-passed} passed
the coshpericity test. However, if we apply the nonlinear transformation
$B^{1}=g(A^{1})$, $B^{2}=g(A^{2})$, where $g$ is given by $0\mapsto0$,
$1\mapsto1$, $5\mapsto2$, we get\[\small%
\begin{tabular}{r|ccclr|ccc}
\multicolumn{4}{l}{$\lambda^{1}=1$, $\lambda^{2}=1$} &  & \multicolumn{4}{l}{$\lambda^{1}=1$, $\lambda^{2}=2$}\tabularnewline
 & $B^{2}=0$ & $B^{2}=1$ & $B^{2}=2$ &  &  & $B^{2}=0$ & $B^{2}=1$ & $B^{2}=2$\tabularnewline
\cline{1-4} \cline{6-9} 
$B^{1}=0$ & .24 & .07 & 0 &  & $B^{1}=0$ & .24 & .07 & 0\tabularnewline
$B^{1}=1$ & .07 & .24 & .07 &  & $B^{1}=1$ & .07 & .24 & .07\tabularnewline
$B^{2}=2$ & 0 & .07 & .24 &  & $B^{2}=2$ & 0 & .07 & .24\tabularnewline
\multicolumn{1}{r}{} &  &  &  &  & \multicolumn{1}{r}{} &  &  & \tabularnewline
\multicolumn{1}{r}{} &  &  &  &  & \multicolumn{1}{r}{} &  &  & \tabularnewline
\multicolumn{4}{l}{$\lambda^{1}=2$, $\lambda^{2}=1$} &  & \multicolumn{4}{l}{$\lambda^{1}=2$, $\lambda^{2}=2$}\tabularnewline
 & $B^{2}=0$ & $B^{2}=1$ & $B^{2}=2$ &  &  & $B^{2}=0$ & $B^{2}=1$ & $B^{2}=2$\tabularnewline
\cline{1-4} \cline{6-9} 
$B^{1}=0$ & .24 & .07 & 0 &  & $B^{1}=0$ & 0 & .07 & .24\tabularnewline
$B^{1}=1$ & .07 & .24 & .07 &  & $B^{1}=1$ & .07 & .24 & .07\tabularnewline
$B^{2}=2$ & 0 & .07 & .24 &  & $B^{2}=2$ & .24 & .07 & 0\tabularnewline
\end{tabular}\]and the correlations for these joint distributions are $\rho_{11}=\rho_{12}=\rho_{21}\approx.7742$
and $\rho_{22}\approx-.7742$. Substituting these in the cosphericity
test, we obtain
\begin{align*}
1.1988 & \approx|.7742\cdot.7742-.7742(-.7742)|\\
 & \le\sqrt{1-.7742^{2}}\sqrt{1-.7742^{2}}+\sqrt{1-.7742^{2}}\sqrt{1-.7742^{2}}\approx.8012.
\end{align*}
We see that the cosphericity test is not passed for the transformed
variables. As a result selective influences are ruled out for the
original variables as well.\hfill$\square$
\end{example}
The cosphericity test can also be applied to more than two input-output
pairs. If we assume that $\left(A^{1},\ldots,A^{n}\right)\looparrowleft\left(\lambda^{1},\ldots,\lambda^{n}\right)$,
then, by the nestedness property for input-output pairs, for any two
of them, $\left(A^{k},\lambda^{k}\right)$ and $\left(A^{k'},\lambda^{k'}\right)$,
we should have $\left(A^{k},A^{k'}\right)\looparrowleft\left(\lambda^{k},\lambda^{k'}\right)$.
The test only applies if there are two values $i$ and $i'$ of $\lambda^{k}$
and two values $j$ and $j'$ of $\lambda^{k'}$ such that, for some
allowable treatments $\phi_{ij},\phi_{ij'},\phi_{i'j},\phi_{i'j'}$,
\[
\lambda^{k}=i,\lambda^{k'}=j\in\phi_{ij},\lambda^{k}=i,\lambda^{k'}=j'\in\phi_{ij'},\textnormal{ etc}.
\]
In other words, the inputs and their values should be chosen so that
$\left\{ \lambda^{k}=i,\lambda^{k}=i'\right\} $ and $\left\{ \lambda^{k'}=j,\lambda^{k'}=j'\right\} $
form a completely crossed subdesign within the set of allowable treatments.
By the nestedness property for input values, we have $\left(A^{k},A^{k'}\right)\looparrowleft\left(\lambda^{k},\lambda^{k'}\right)$
with the input values restricted to $\left\{ i,i'\right\} $ and $\left\{ j,j'\right\} $
and the new set of allowable treatments consisting of all four possible
combinations. If this cosphericity inequality is violated for all
least one combination of $k,k',i,i',j,j'$, then the initial assumption
$\left(A^{1},\ldots,A^{n}\right)\looparrowleft\left(\lambda^{1},\ldots,\lambda^{n}\right)$
should be rejected.

\section{\label{Conditional determinism}Conditional determinism and conditional independence of outcomes}

The definition of selective influences (in the canonical form) requires
the existence of a random variable $R$ and functions $f_{1},\ldots,f_{n}$
such that, for all allowable treatments $\phi$, 
\[
\left(A_{\phi}^{1},\ldots,A_{\phi}^{n}\right)\sim\left(f_{1}\left(\lambda^{1},R\right),\ldots,f_{n}\left(\lambda^{n},R\right)\right).
\]
One obvious consequence of this definition is that, conditioned on
any value $r$ of $R$, the outputs become (equal to) deterministic
functions of the corresponding factors,
\[
f_{1}\left(\lambda^{1},r\right),\ldots,f_{n}\left(\lambda^{n},r\right).
\]
It is sometimes easy to deal with these deterministic quantities,
derive certain inequalities that hold for every value of $r$, and
then show that they are preserved as $R$ randomly varies. It is an
especially useful approach if the distributions of $\left(A_{\phi}^{1},\ldots,A_{\phi}^{n}\right)$
at allowable treatments $\phi$ are not known, and instead we know
distributions of certain functions of these random variables, such
as their sums or maxima. 

Let us discuss this on an example from studies of mental architectures.
This is a traditional area of psychology dealing with decomposing performance of a task into a network of subprocesses when we only know the distributions of the overall performance time (referred to as response time) at different treatments. Let us assume that we observe response times $T$ in an experiment
with two factors, $\lambda^{1},\lambda^{2}$, manipulated at two levels
each, denoted in both cases by 1 and 2. All four treatments are allowable.
Let us postulate that there are two processes involved, with their
durations $A^{1}$ and $A^{2}$ being random variables, and that $\left(A^{1},A^{2}\right)\looparrowleft\left(\lambda^{1},\lambda^{2}\right)$.
We want to determine which of the three ``architectures,'' or composition
schemes, is being employed: 
\begin{enumerate}
\item serial, $T_{\phi}=A_{\phi}^{1}+A_{\phi}^{2}$
\item parallel-OR, $T_{\phi}=\min\left(A_{\phi}^{1},A_{\phi}^{2}\right)$,
or
\item parallel-AND, $T_{\phi}=\max\left(A_{\phi}^{1},A_{\phi}^{2}\right)$. 
\end{enumerate}
One tool traditionally used for this purpose is the \emph{interaction
contrast},
\[
c\left(t\right)=\Pr\left(T_{11}\leq t\right)+\Pr\left(T_{22}\leq t\right)-\Pr\left(T_{12}\leq t\right)-\Pr\left(T_{21}\leq t\right),
\]
where $t$ is any non-negative number, and $T_{ij}$ abbreviates $T_{\phi=\left(i,j\right)}$. 

We do not know the joint distribution of $A_{\phi}^{1},A_{\phi}^{2}$
at any of the four treatments, but we can write
\[
\left(A_{ij}^{1},A_{ij}^{2}\right)\sim\left(f_{1}\left(\lambda^{1}=i,R\right),f_{2}\left(\lambda^{2}=j,R\right)\right)=\left(g_{i}^{1}\left(R\right),g_{j}^{2}\left(R\right)\right),\; i,j\in\left\{ 1,2\right\} .
\]
We need one additional assumption: that $R$ can be chosen in such
a way that, for any of its possible values $r$,
\[
g_{1}^{1}\left(r\right)\leq g_{2}^{1}\left(r\right),\quad g_{1}^{2}\left(r\right)\leq g_{2}^{2}\left(r\right).
\]
In other words, switching either factor from level 1 to level 2 prolongs
the corresponding processing time. We call this assumption \emph{prolongation
constraints}. Various analogues of this assumption are common in studies
of mental architectures.

Deterministic real-valued quantities can be viewed as random variables
with Heaviside distribution functions:
\[
\Pr\left(g_{l}^{k}\left(r\right)\leq t\right)=\left\{ \begin{array}{ccc}
0 & \textnormal{if} & t<g_{l}^{k}\left(r\right),\\
1 & \textnormal{if} & t\geq g_{l}^{k}\left(r\right).
\end{array}\right.
\]
Analogously,
\[
\Pr\left(\textnormal{comp}\left(g_{i}^{1}\left(r\right),g_{j}^{2}\left(r\right)\right)\leq t\right)=\left\{ \begin{array}{ccc}
0 & \textnormal{if} & t<\textnormal{comp}\left(g_{i}^{1}\left(r\right),g_{j}^{2}\left(r\right)\right),\\
1 & \textnormal{if} & t\geq\textnormal{comp}\left(g_{i}^{1}\left(r\right),g_{j}^{2}\left(r\right)\right),
\end{array}\right.
\]
where $\textnormal{comp}$ stands for one of the three composition
rules of interest, plus, maximum, or minimum. This allows us to form
the \emph{conditional interaction contrast},
\[
c^{*}\left(t,r\right)=\Pr\left(t_{11}\left(r\right)\leq t\right)+\Pr\left(t_{22}\left(r\right)\leq t\right)-\Pr\left(t_{12}\left(r\right)\leq t\right)-\Pr\left(t_{21}\left(r\right)\leq t\right),
\]
where
\[
t_{ij}=\textnormal{comp}\left(g_{i}^{1}\left(r\right),g_{j}^{2}\left(r\right)\right).
\]

It is easy to see that
\[
\Pr\left(T_{ij}\leq t\right)=\int_{S_{R}}\Pr\left(t_{ij}\left(r\right)\leq t\right)\textnormal{d}p_{R}\left(r\right)
\]
and 
\[
c\left(t\right)=\int_{S_{R}}c^{*}\left(t,r\right)\textnormal{d}p_{R}\left(r\right),
\]
where the Lebesgue integral is taken over the entire domain $S_{R}$
of $R$, and $p_{R}$ is the probability measure in the distribution
of $R$. (The reader not familiar with Lebesgue integrals can think
of $\textnormal{d}p_{R}\left(r\right)$ above as a generalized version
of $f_{R}\left(r\right)\textnormal{d}r$, where $f_{R}$ is the density
function of $R$ over the set of real numbers.) 

Using this observation we can easily establish that if the composition
rule is $\min$ (parallel-OR architecture), then $c\left(t\right)\leq0$,
for all $t$, because $c^{*}\left(t,r\right)\leq0$ at any $t$ and
any fixed $r$. Indeed, consider all possible arrangements of $g_{1}^{1}\left(r\right),g_{2}^{1}\left(r\right),g_{1}^{2}\left(r\right),g_{2}^{2}\left(r\right)$
keeping in mind the prolongation constraints and assuming, with no
loss of generality, that $g_{1}^{1}\left(r\right)\leq g_{1}^{2}\left(r\right)$.
These possible arrangements are
\[
\begin{array}{cc}
\textnormal{(i)} & g_{1}^{1}\left(r\right)\leq g_{2}^{1}\left(r\right)\leq g_{1}^{2}\left(r\right)\leq g_{2}^{2}\left(r\right),\\
\textnormal{(ii)} & g_{1}^{1}\left(r\right)\leq g_{1}^{2}\left(r\right)\leq g_{2}^{1}\left(r\right)\leq g_{2}^{2}\left(r\right),\\
\textnormal{(iii)} & g_{1}^{1}\left(r\right)\leq g_{1}^{2}\left(r\right)\leq g_{2}^{2}\left(r\right)\leq g_{2}^{1}\left(r\right).
\end{array}
\]
Thus, for (ii), we have
\[
\begin{array}{c}
t_{11}\left(r\right)=\min\left(g_{1}^{1}\left(r\right),g_{1}^{2}\left(r\right)\right)=g_{1}^{1}\left(r\right),\\
t_{12}\left(r\right)=\min\left(g_{1}^{1}\left(r\right),g_{2}^{2}\left(r\right)\right)=g_{1}^{1}\left(r\right),\\
t_{21}\left(r\right)=\min\left(g_{2}^{1}\left(r\right),g_{1}^{2}\left(r\right)\right)=g_{1}^{2}\left(r\right),\\
t_{22}\left(r\right)=\min\left(g_{2}^{1}\left(r\right),g_{2}^{2}\left(r\right)\right)=g_{2}^{1}\left(r\right).
\end{array}
\]
Then, substituting for the numerical values
\begin{align*}
c^{*}\left(t,r\right) & =\Pr\left(t_{11}\left(r\right)\leq t\right)+\Pr\left(t_{22}\left(r\right)\leq t\right)-\Pr\left(t_{12}\left(r\right)\leq t\right)-\Pr\left(t_{21}\left(r\right)\leq t\right)\\
 & =\begin{cases}
0+0-0-0=0 & \textnormal{if }t<g_{1}^{1}\left(r\right),\\
1+0-1-0=0 & \textnormal{if }g_{1}^{1}\left(r\right)\leq t<g_{1}^{2}\left(r\right),\\
1+0-1-1<0 & \textnormal{if }g_{1}^{2}\left(r\right)\leq t<g_{2}^{1}\left(r\right),\\
1+1-1-1=0 & \textnormal{if }g_{2}^{1}\left(r\right)\leq t<g_{2}^{2}\left(r\right),\\
1+1-1-1=0 & \textnormal{if }t\geq g_{2}^{2}\left(r\right).
\end{cases}
\end{align*}
In the same way one proves that $c^{*}\left(t,r\right)$ is never
positive in cases (i) and (iii). 

By analogous reasoning we can show that if the composition rule is
$\max$ (parallel-AND architecture), then $c\left(t\right)\geq0$,
for all $t$, because $c^{*}\left(t,r\right)\geq0$ at any $t$ and
any fixed $r$. 

For the serial architecture (the composition rule $+$) $c^{*}\left(t,r\right)$
does not preserve its sign, but the analysis of the arrangements shows
that, for any $t$ and $r$,
\[
\int_{0}^{t}c^{*}\left(t,r\right)\textnormal{d}t\geq0
\]
and
\[
\int_{0}^{\infty}c^{*}\left(t,r\right)\textnormal{d}t=0.
\]
Then the same properties should hold for $c\left(t\right)$, because
\[
\int_{0}^{t}c\left(t\right)\textnormal{d}t=\int_{0}^{t}\left(\int_{S_{R}}c^{*}\left(t,r\right)\textnormal{d}p_{R}\left(r\right)\right)\textnormal{d}t=\int_{S_{R}}\left(\int_{0}^{t}c^{*}\left(t,r\right)\textnormal{d}t\right)\textnormal{d}p_{R}\left(r\right).
\]

However, dealing with deterministic quantities is not always convenient.
If a deterministic quantity changes as a function of $r$, the probability
with which it falls within a given measurable subset may jump from
0 to 1 or vice versa. In some cases it may be desirable to deal with
``well-behaved'' distributions only, with associated probabilities
that change continuously or even sufficiently smoothly. (The term
``smooth'' refers to the highest order of continuous derivative
a function possesses.) To make this desideratum achievable in the
context of selective influences, we begin by stating the following
equivalence. 
\begin{thm}
\label{TH: stoch. determinism}$\left(A^{1},\ldots,A^{n}\right)\looparrowleft\left(\lambda^{1},\ldots,\lambda^{n}\right)$
if and only if one can find stochastically independent random variables
$R,R^{1},\ldots,R^{n}$ and functions $w_{1},\ldots,w_{n}$, such
that
\[
\left(A_{\phi}^{1},\ldots,A_{\phi}^{n}\right)\sim\left(w_{1}\left(\lambda^{1},R,R^{1}\right),\ldots,w_{n}\left(\lambda^{n},R,R^{n}\right)\right)
\]
for all allowable treatments $\phi=\left(\lambda^{1},\ldots,\lambda^{n}\right)$. 
\end{thm}
By analogy with factor analysis, we can call $R^{1},\ldots,R^{n}$
\emph{specific sources of variability}, and call $R$ a \emph{common
source of variability}. The proof of the theorem is very simple. If
a representation 
\[
\left(A_{\phi}^{1},\ldots,A_{\phi}^{n}\right)\sim\left(f_{1}\left(\lambda^{1},R\right),\ldots,f_{n}\left(\lambda^{n},R\right)\right)
\]
exists, one can choose arbitrary $R^{1},\ldots,R^{n}$ (combined together
and with $R$ by an independent coupling) and put $w_{k}\left(\lambda,r,r'\right)=f_{k}\left(\lambda,r\right)$,
$k=1,\ldots,n$. If a representation stated in the theorem exists,
then define $R^{*}=\left(R,R^{1},\ldots,R^{n}\right)$ and put $f_{k}\left(\lambda,\left(r,r^{1},\ldots,r^{n}\right)\right)=w_{k}\left(\lambda,r,\textnormal{Proj}_{k}\left(r^{1},\ldots,r^{n}\right)\right)$. 

The consequences of this simple theorem are significant. Once the
possibility of splitting a single source of randomness into a common
and specific components has been established, it becomes possible
that in certain situations this split can be more than a formal redefinition
of a single source. It follows from the theorem that conditioned upon
any value $r$ of $R$, the random variables $w_{1}\left(\lambda^{1},r,R^{1}\right),\ldots,w_{n}\left(\lambda^{n},r,R^{n}\right)$
are stochastically independent. One can hypothesize now, that these
independent random variables have distributions with desired properties.
For example, if all random variables $\left(A_{\phi}^{1},\ldots,A_{\phi}^{n}\right)$
are real-valued and continuous, $w_{1}\left(\lambda^{1},r,R^{1}\right),\ldots,w_{n}\left(\lambda^{n},r,R^{n}\right)$
may be assumed to possess densities, or have the property that the
probability with which $w_{k}\left(\lambda^{k},r,R^{k}\right)$ falls
within any interval of reals is a continuously differentiable function
of $r$. Such assumptions may be important in studying mental architectures
or random variables underlying comparisons of stimuli.

\section{Related literature}

There are many textbooks treating measure theory and probability (e.g.,
\citealp{Chung1974}). However, the reader should be aware that (a)
older textbooks usually deal with random variables in the narrow sense
only; (2) in most textbooks random variables are defined as measurable
functions on a sample space, restricting thereby the consideration
to jointly distributed random variables. For random variables that
need not be jointly distributed and the associated theory of coupling
them into jointly distributed entities, see \citet{Thorisson2000}.
The earliest explicit discussions of selective influences in psychology
can be found in \citet{Sternberg1969} and \citet{Townsend1984}.
Marginal selectivity for two random variables was first mentioned
in \citet{TownsendSchweickert1989}. Other historical details and
relations can be found in \citet{Dzhafarov2003c}, where the theory
of selective influences presented in this chapter was first proposed.
In this earlier work (and its elaboration in \citet{DzhafarovGluhovsky2006}
the ``is distributed as'' relation in the defining representation
for selective influences,
\[
\left(A_{\phi}^{1},\ldots,A_{\phi}^{n}\right)\sim\left(f_{1}\left(\lambda^{1},R\right),\ldots,f_{n}\left(\lambda^{n},R\right)\right)
\]
was somewhat carelessly replaced with equality. For a mathematically
rigorous and maximally general version of the definition and Joint
Distribution Criterion, see \citet{DzhafarovKujala2010}. The tests
of selective influences were first introduced in \citet{KujalaDzhafarov2008b}.
They included the cosphericity tests and a special form of distance
tests. A general version of distance tests (p.q.-metric tests) was
introduced in \citet{orderdist}. The linear feasibility test is described
in \citet{DzhafarovKujala2012a}. For applications of the theory of
selective influences to discrimination judgments and to mental processing
architectures, see \citet{Dzhafarov2003a,Dzhafarov2003b} and \citet*{DzhafarovSchweickertSung2004}.
The parallels between the theory of selective influences and the analysis
of determinism in the so-called Bohmian version of the Einsten-Podolsky-Rosen
entanglement paradigm of quantum physics are described in \citet{DzhafarovKujala2012a,DzhafarovKujala2012b}.
The history there dates back to \citeauthor{Bell1964}'s \citeyearpar{Bell1964}
epoch-making inequalities, and then to their elborations in \citet*{ClauserHorneShimonyHolt1969}
and \citet{Fine1982a}. Mathematically, this line of work is subsumed
by the linear feasibility test, whose most general version in quantum
physics is described in \citet{BasoaltoPercival2003}. 

\section{Acknowledgments}
This work was supported by NSF grant SES-1155956. We are grateful to Jing Chen, Shree Frazier, Nicole Murchison,  Alison Schroeder, and Ru Zhang for pointing out numerous typos and imprecisions in the original draft of the chapter.
\bibliography{si}

\end{document}